\documentclass[11pt]{article}
\usepackage[a4paper,margin=0.95in]{geometry}
\usepackage{graphicx}
\usepackage{amsmath,amssymb,amsthm,bm,mathtools,booktabs,microtype,array}
\usepackage{algorithm,algpseudocode}
\usepackage{placeins}
\usepackage[colorlinks=true,allcolors=blue,pdfauthor={Minghui Li, Zixian Liu, Shipeng Xu},pdftitle={A Second-Order Original-Energy-Dissipative Convex Splitting and Balanced-Tree hp Spectral-Element Method for Thermodynamically Consistent Anisotropic Dendritic Solidification}]{hyperref}

\newtheorem{theorem}{Theorem}[section]
\newtheorem{lemma}[theorem]{Lemma}
\newtheorem{proposition}[theorem]{Proposition}
\newtheorem{remark}[theorem]{Remark}

\newcommand{\R}{\mathbb{R}}
\newcommand{\Eorig}{\mathcal{E}_I}
\newcommand{\Ec}{\mathcal{E}_c}
\newcommand{\Ee}{\mathcal{E}_e}

\title{A Second-Order Original-Energy-Dissipative Convex Splitting and Balanced-Tree \(hp\) Spectral-Element Method for Thermodynamically Consistent Anisotropic Dendritic Solidification}
\author{Minghui Li\thanks{Corresponding author. E-mail address: \href{mailto:minghuili_xmu@163.com}{\texttt{minghuili\_xmu@163.com}}.} \qquad Zixian Liu \qquad Shipeng Xu \\[6pt]
\small\itshape School of Mathematical Sciences \\
\small\itshape Jiangxi Science and Technology Normal University, Nanchang 330038, China}
\date{}

\begin{document}
\maketitle

\begin{abstract}
We present a unified variational framework for two- and three-dimensional thermodynamically consistent anisotropic dendritic solidification under homogeneous Neumann boundary conditions that achieves second-order fixed-space temporal accuracy, unconditional dissipation of the \emph{original} quadrature-defined discrete free energy for exact stages, exact enthalpy conservation, and efficient moving-interface adaptivity. Eliminating the constant temperature mode via the conserved mean enthalpy recasts the coupled system as an exact gradient flow with a state-dependent symmetric positive-definite metric. Pairing a sextic majorizer for the quintic latent-heat lift with an explicit anisotropy shift yields a global convex splitting of the original energy for both convex and nonconvex anisotropies in 2D and 3D; combined with a midpoint-metric five-stage convex-splitting Runge--Kutta (CSRK) scheme, we prove local second-order consistency on smooth fixed-space trajectories, unique solvability, exact enthalpy conservation, and unconditional nonincrease of the original discrete energy for every \(k>0\). Spatially, conforming \(hp\) spectral elements on \(2{:}1\)-balanced quadtrees and octrees use an explicit \(O(N)\)-nonzero trace injection, common Gauss quadrature, and Schur-complement thermal elimination with matrix-free 3D sum factorization for the nonlinear phase Hessian and coupling operators. A moving-interface \(hp\) strategy on \(\phi\in[-0.98,0.98]\) with a \(Q_2\) thermal halo, energy-safe radial thermal transfer, and explicit step/transfer acceptance guards guarantees global quadrature-energy nonincrease in exact arithmetic (and within stated floating-point tolerances in double precision) across certified steps and accepted transfers. Insulated-Neumann 2D and 3D benchmarks---including CSRK-versus-SAV temporal comparisons, asymmetric rotated and multi-nuclei growth, and independent higher-order \(Q_4\) ablations---demonstrate \(7.4\times\)--\(75.3\times\) average DOF compression, \(64.5\times\) full-run speedup over matched uniform 2D \(L_6\,Q_3\) in the 2D ablation, and \(4.7\times\) fewer DOFs and \(2.8\times\) lower runtime than 3D \(h\)-adaptive \(Q_3\) at \(0.35\%\) phase \(L^2\) error difference (\(t=5\)).
\end{abstract}

\noindent\textbf{Keywords:} anisotropic dendritic solidification; thermodynamic consistency; second-order convex splitting; original energy dissipation; \(hp\) spectral-element method; \(2{:}1\)-balanced quadtree/octree adaptivity.

% ============================================================================
\section{Introduction}\label{sec:intro}
% ============================================================================
Dendritic crystal growth is one of the most ubiquitous and fundamental pattern-formation phenomena in the solidification of pure undercooled melts, metallic alloys, and functional materials \cite{boettinger2002}. Physically, dendritic morphologies emerge from the interplay between destabilizing thermal transport and stabilizing interfacial energetics: latent heat rejected at the advancing solid--liquid boundary diffuses into the undercooled melt and triggers the Mullins--Sekerka morphological instability \cite{mullins1963}, while capillarity (the Gibbs--Thomson effect) and crystallographic anisotropy counter short-wavelength perturbations and select preferred growth axes, primary trunks, and secondary sidebranches. Because the resulting dendritic microstructure---tip velocity, arm spacing, branching symmetry, and segregation---directly governs mechanical strength, ductility, corrosion resistance, and transport properties in casting, welding, additive manufacturing, and semiconductor crystal growth, accurate quantitative prediction of two- and three-dimensional dendritic evolution carries substantial scientific and engineering value \cite{boettinger2002}.

Classical sharp-interface Stefan formulations track the moving solid--liquid front explicitly, which becomes notoriously difficult under strong curvature variations, sidebranching, and 3D topological coalescence. Phase-field modeling avoids front tracking by introducing a smooth order parameter \(\phi\in[-1,1]\) separating solid (\(\phi=+1\)) from liquid (\(\phi=-1\)) across a diffuse layer of width \(O(\varepsilon)\). Building on early free-boundary phase-field models \cite{caginalp1986} and thermodynamically consistent derivations \cite{penrose1990,wang1993}, Kobayashi \cite{kobayashi} pioneered planar anisotropic dendritic simulations, Karma and Rappel \cite{karma1996} established the quantitative 2D/3D thin-interface anisotropic model, and Eggleston, McFadden, and Voorhees \cite{eggleston2001} analyzed strongly anisotropic nonconvex interfacial energies. Here we adopt the widely used thermodynamically consistent anisotropic solidification model of Karma and Rappel \cite{karma1996} (see also \cite{guo2024,lithesis,li2022,yang2019,yang2021b,zhao2017}), which couples an anisotropic Allen--Cahn equation with conservative enthalpy diffusion via the quintic latent-heat interpolation \(h(\phi)=\phi-\frac{2}{3}\phi^3+\frac{1}{5}\phi^5\) and a quadratic thermal free energy. Under homogeneous Neumann (insulated) boundary conditions, this system satisfies exact total-enthalpy conservation and monotonic dissipation of the original coupled free energy, while exhibiting \(O(\varepsilon^{-2})\) interfacial stiffness, nonlinear phase--thermal coupling, and nonconvex gradient energy at strong anisotropy.

On the temporal discretization side, preserving these thermodynamic laws at high order remains challenging. Early explicit or semi-implicit integrators \cite{karma1996,kobayashi} face severe stability limits, while operator- and direction-splitting schemes for dendritic growth \cite{likim2011,wangxiao2022} improve efficiency without proving discrete energy dissipation. First-order convex splitting \cite{eyre1998} unconditionally dissipates the original energy but incurs phase lag over long times, whereas second-order multistep convex-splitting and orientation-field schemes \cite{likim2023,shen2012} either rely on modified energies or problem-specific splittings. More recently, linear auxiliary-variable methods based on invariant energy quadratization (IEQ) \cite{yang2016,yang2019,zhang2019,zhao2017} and the scalar auxiliary variable (SAV) or relaxed generalized SAV frameworks \cite{guo2024,li2022,shen2018,wang2022,yang2020,yang2021a,yang2021b,zhangshen2022} have been extensively applied to 2D/3D dendritic solidification and thermal phase change, progressively reducing each step from five decoupled linear elliptic solves \cite{yang2021b} to four \cite{li2022} and two \cite{guo2024}. Yet these cited linearized auxiliary-variable schemes guarantee stability only for a \emph{modified} energy or auxiliary scalar rather than the true physical free energy, which can drift during rapid front motion and depends on empirical linear stabilizers at larger step sizes. Conversely, the multistage convex-splitting Runge--Kutta (CSRK) framework \cite{zeg} achieves second-order accuracy with unconditional dissipation of the \emph{original} energy, but extending it to coupled anisotropic solidification is nontrivial due to the non-constant coupled metric, nonconvex strong anisotropy \cite{eggleston2001}, and global enthalpy--latent-heat coupling.

On the spatial discretization side, dendritic growth combines a thin moving \(O(\varepsilon)\) phase layer with a much wider thermal diffusion layer. Uniform structured-grid solvers---such as GPU finite-difference and parallel direction-splitting schemes \cite{wangxiao2022,yamanaka2011}---and global spectral methods---including periodic Fourier pseudospectral \cite{guo2024,li2022,yang2019,yang2021b}, global spectral-Galerkin \cite{yang2021a}, and single-domain Legendre--Galerkin \cite{wang2022} discretizations---waste most degrees of freedom (DOFs) in the nearly constant bulk, impose periodic boundaries or lack local refinement, and become prohibitively costly in 3D. Adaptive finite-element, finite-volume, and quadtree/octree adaptive mesh refinement (AMR) solvers for dendritic solidification \cite{burstedde2011,greenwood2018,jeong2001,lan2003,likim2012,provatas1999,rosam2007,schmidt1996} dramatically reduce DOFs, yet conventional low-order AMR converges only algebraically across the steep interfacial and thermal layers. Adaptive \(hp\) spectral-element methods (\(hp\)-SEM) \cite{babuska1994,demkowicz2006,karniadakis2005,shentangwang} combine local tree \(h\)-refinement, high interfacial polynomial degree \(p\), and bulk coarsening for spectral accuracy, but enforcing scalable \(H^1\)-conforming hanging constraints while preserving exact enthalpy conservation and original-energy monotonicity across dynamic mesh updates remains an open challenge.

Bridging these temporal and spatial gaps simultaneously---achieving second-order accuracy, unconditional dissipation of the \emph{original} physical energy, exact enthalpy conservation under insulated Neumann boundaries, and provably energy-monotone 2D/3D \(hp\) adaptivity---is essential for faithful and tractable dendritic simulations. To this end, we develop a unified second-order convex splitting and \(2{:}1\)-balanced \(hp\) spectral-element framework with four main contributions:
\begin{enumerate}
\item \textbf{Enthalpy-conserving gradient-flow formulation and global convex splitting.} Eliminating the constant temperature mode via the conserved mean enthalpy \(\bar{H}\) casts the coupled PDEs as an exact gradient flow on \(X_c = H^1(\Omega)\times L^2_0(\Omega)\) with a state-dependent symmetric positive-definite (SPD) inverse metric \(J_\phi\). Combining a sextic majorizer for the quintic latent-heat lift with an explicit anisotropy shift \(\mu_A\) yields a global convex splitting \(E = E_c - E_e\) of the \emph{original} energy valid for both convex and nonconvex anisotropies in 2D and 3D.
\item \textbf{Provable original-energy dissipation and fast Schur-complement solver.} Pairing this splitting with a midpoint-metric predictor and five-stage CSRK tableau \cite{zeg}, we prove local second-order consistency on smooth trajectories, unique solvability, exact enthalpy conservation, and unconditional nonincrease of the \emph{original} quadrature-defined discrete energy for every \(k>0\) at both the time-semidiscrete and fully discrete \(hp\)-SEM levels. Exact Schur-complement elimination of temperature reduces each stage to a strictly convex phase-only minimization---solved via damped Newton--PCG with matrix-free 3D sum factorization for the nonlinear phase Hessian and coupling operators---accompanied by a computable two-component a posteriori residual-work energy certificate for inexact phase and iterative thermal solves.
\item \textbf{Balanced-tree moving-interface \(hp\) adaptivity with certified global energy monotonicity.} Enforcing full \(2{:}1\) balance on quadtrees (2D) and octrees (3D) yields an explicit \(O(N)\)-nonzero coarse-to-fine trace injection \(J\) producing continuous inter-element traces (verified numerically to \(\sim 10^{-16}\)). Concentrating finest \(h\)-cells and high degree \(Q_3\) on \(\phi\in[-0.98,0.98]\) with a \(Q_2\) thermal halo and coarsened \(Q_1\) bulk, we pair exact common-subdivision \(L^2\) projection (accelerated by \(O(L_{\max}(N_{\mathrm{old}}+N_{\mathrm{new}}))\) dyadic ancestor/descendant hashing) with an enthalpy shift, energy-safe radial thermal contraction (Proposition~\ref{prop:energy-transfer}), and an explicit trial-and-fallback acceptance guard to guarantee monotonic original discrete energy decay across all accepted mesh adaptations.
\item \textbf{Comprehensive 2D/3D verification, temporal benchmarking, and independent higher-order ablations.} Under insulated Neumann boundaries, we verify temporal and \(h\)/\(p\) convergence, benchmark CSRK directly against a linearized stabilized SAV--BDF2 implementation on identical meshes, simulate strong-anisotropy 2D (\(m=4,6\)), full-domain asymmetric rotated and multi-nuclei, and 3D cubic dendrites with \(7.4\times\)--\(75.3\times\) average DOF compression (and \(64.5\times\) full-run 2D / \(2.8\times\) 3D equal-condition ablation speedups against uniform and \(h\)-only \(Q_3\) baselines), and conduct equal-condition 2D and 3D ablations against independent higher-order \(Q_4\) spatial references documenting the accuracy--cost trade-offs of balanced tree \(h\)-branching, zoned \(hp\) enrichment, the \([-0.98,0.98]\) band, and the \(Q_2\) thermal halo.
\end{enumerate}

The remainder of the paper is organized as follows. Section~\ref{sec:model} presents the thermodynamically consistent coupled model, its anisotropy structure, and its gradient-flow reformulation. Section~\ref{sec:semidiscrete} constructs the second-order time-semidiscrete convex-splitting scheme and proves original-energy dissipation, enthalpy conservation, and unique solvability. Section~\ref{sec:fullydiscrete} formulates the fully discrete \(hp\) spectral-element scheme, establishes its discrete energy and solvability theorems, and details the Schur-complement Newton--PCG solver. Section~\ref{sec:hp} presents the \(2{:}1\)-balanced moving-interface \([-0.98,0.98]\) \(hp\) strategy and proves global energy monotonicity across mesh transfers. Section~\ref{sec:numerics} reports the 2D accuracy tests and SAV--BDF2 comparison, 2D symmetric and full-domain asymmetric dendritic growth simulations, 3D cubic-anisotropy dendritic simulation, and independent higher-order \(Q_4\) ablations. Conclusions are drawn in Section~\ref{sec:conclusions}.

% ============================================================================
\section{Thermodynamically consistent coupled phase-field model}\label{sec:model}
% ============================================================================
Let \(\Omega\subset\R^d\) (\(d\in\{2,3\}\)) be a connected rectangular or box domain with outward unit normal \(\bm{n}\) on \(\partial\Omega\). Denote by \(\phi(\bm{x},t)\in[-1,1]\) the phase field (\(\phi=+1\) in the solid and \(\phi=-1\) in the liquid) and by \(T(\bm{x},t)\) the dimensionless temperature. Adopting the thermodynamically consistent anisotropic dendritic solidification model of Karma and Rappel \cite{karma1996} (see also \cite{lithesis,li2022,yang2019,zhao2017} for its dimensionless coupled formulation in two and three dimensions), let
\begin{equation}
 h(s)=s-\frac{2}{3}s^3+\frac{1}{5}s^5, \qquad
 F(s)=\frac{1}{4}(s^2-1)^2, \qquad
 c_T=\frac{\lambda}{\varepsilon K}>0,
 \label{eq:constitutive}
\end{equation}
where \(\varepsilon>0\) is the interfacial thickness parameter, \(\lambda>0\) is the phase--thermal coupling constant, and \(K>0\) is the dimensionless latent-heat parameter. Note that \(h'(s)=(1-s^2)^2\ge 0\) vanishes at the pure phases \(s=\pm 1\), confining latent-heat release and thermal driving to the diffuse interface.

With constant phase relaxation coefficient \(\varrho>0\) and thermal diffusivity \(D>0\), the coupled phase-field and enthalpy equations are
\begin{align}
 \varrho\,\phi_t &= \nabla\cdot\nabla_g A(\nabla\phi) - \varepsilon^{-2}F'(\phi) - c_T K h'(\phi)T && \text{in }\Omega\times(0,T_f], \label{eq:pde-phase}\\
 \partial_t\bigl(T - K h(\phi)\bigr) &= D\Delta T && \text{in }\Omega\times(0,T_f], \label{eq:pde-heat}
\end{align}
subject to homogeneous Neumann (insulated, no-flux) boundary conditions on every exterior face:
\begin{equation}
 \nabla_g A(\nabla\phi)\cdot\bm{n}=0, \qquad D\nabla T\cdot\bm{n}=0 \qquad \text{on }\partial\Omega\times(0,T_f].
 \label{eq:bc-neumann}
\end{equation}
Note that \eqref{eq:bc-neumann} is the exact natural variational boundary condition of the weak formulation below for arbitrary crystal orientation \(\theta_0\) or orthogonal rotation \(\mathcal{R}\); moreover, for axis-aligned anisotropy on a rectangular box, the gradient density \(A\) is even in the normal derivative on every face, so \(\partial_{\bm{n}}\phi=\partial_{\bm{n}}T=0\) implies \eqref{eq:bc-neumann} and ensures that any quadrant/octant homogeneous Neumann solution extends by reflection to a symmetric weak solution on the full box (coinciding with the full-domain solution within the reflection-symmetric class or whenever uniqueness holds).

\subsection{Anisotropic gradient energy and global convex stabilization}
The anisotropic surface energy density \(A:\R^d\to[0,\infty)\) is positively two-homogeneous, \(A(g)=\frac{1}{2}|g|^2 b(\hat{g})^2\) for \(g\ne 0\) with \(\hat{g}=g/|g|\) and \(A(0)=0\):
\begin{itemize}
\item \textbf{Two-dimensional \(m\)-fold anisotropy (\(m\in\{4,6\}\)):} writing \(\vartheta=\arg(g)\) and crystal orientation \(\theta_0\),
\begin{equation}
 b(\vartheta)=1+\sigma\cos\bigl(m(\vartheta-\theta_0)\bigr), \qquad |\sigma|<1.
 \label{eq:aniso-2d}
\end{equation}
\item \textbf{Three-dimensional cubic anisotropy:} writing \(n_g=\mathcal{R}\hat{g}\) for an orthogonal crystal rotation \(\mathcal{R}\),
\begin{equation}
 b(n_g)=1-3\delta+4\delta\sum_{i=1}^3 (n_g)_i^4, \qquad 0\le\delta<\frac{3}{5}.
 \label{eq:aniso-3d}
\end{equation}
\end{itemize}

The directional curvature of \(A\) governs both ellipticity and convex splitting. In 2D, decomposing any perturbation \(d\in\R^2\) at \(g\ne 0\) into radial and tangential components \((d_r,d_\vartheta)\) gives the exact Hessian quadratic form
\begin{equation}
 d^T \nabla_g^2 A(g)\, d = (b d_r + b' d_\vartheta)^2 + b(b+b'')d_\vartheta^2,
 \qquad b+b'' = 1 - (m^2-1)\sigma\cos\bigl(m(\vartheta-\theta_0)\bigr).
 \label{eq:hess-2d}
\end{equation}
Thus \(A\) is globally convex if and only if \(|\sigma|\le 1/(m^2-1)\) (i.e., \(|\sigma|\le 1/15\) for \(m=4\) and \(|\sigma|\le 1/35\) for \(m=6\)), and uniformly elliptic if and only if \(|\sigma|<1/(m^2-1)\), because at equality \(|\sigma|=1/(m^2-1)\) both \(b'=0\) and \(b+b''=0\) at the critical angles so pure tangential perturbations yield zero curvature (while \(A\) remains strictly convex on \(\R^2\), as any non-radial line segment meets these finitely many critical rays only at isolated points and radial segments have positive curvature \(b^2>0\)). For any \(|\sigma|<1\), \(b\in[1-|\sigma|,1+|\sigma|]\) and \(b+b''\ge 1-(m^2-1)|\sigma|\), so \(d^T\nabla_g^2 A(g)d \ge -\mu_A |d|^2\) with
\begin{equation}
 \mu_A^{(2\mathrm{D})} = (1+|\sigma|)\bigl((m^2-1)|\sigma|-1\bigr)_+, \qquad x_+ = \max(x,0).
 \label{eq:muA-2d}
\end{equation}
In 3D, setting \(s=\sum_{i=1}^3(n_g)_i^4\in[1/3,1]\), \(q=16\delta(n_g^{\circ 3}-s n_g)\) (with componentwise cube \((n_g^{\circ 3})_i=(n_g)_i^3\)), and decomposing any perturbation \(d\in\R^3\) in crystal coordinates as \(\mathcal{R}d=\tau n_g + z\) with \(\tau\in\R\) and \(z\perp n_g\), direct differentiation gives
\begin{equation}
 d^T \nabla_g^2 A(g)\,d
 = (b\tau + q^T z)^2 + b\Bigl[b|z|^2 + 16\delta\sum_{i=1}^3\bigl(3(n_g)_i^2 - s\bigr)z_i^2\Bigr].
 \label{eq:hess-3d}
\end{equation}
Since \(b\in[1-5\delta/3,1+\delta]\) (so positivity \(b_{\min}=1-5\delta/3>0\) holds precisely for \(0\le\delta<3/5\)) and \(b-16\delta s\ge 1-15\delta\), the 3D density is uniformly elliptic (and strictly convex) for \(0\le\delta<1/15\), with explicit lower Hessian bound \(c_A=\min\{b_*^2/2,\,b_*(1-15\delta)/(1+2B^2)\}>0\) where \(b_*=1-5\delta/3\) and \(B=8\delta/b_*\). For general \(0\le\delta<3/5\), the Hessian is bounded below by \(-\mu_A^{(3\mathrm{D})}I\) with
\begin{equation}
 \mu_A^{(3\mathrm{D})} = (1+\delta)(15\delta-1)_+.
 \label{eq:muA-3d}
\end{equation}
Notice that \(A\in C^{1,1}(\R^d)\cap C^\infty(\R^d\setminus\{0\})\): although \(\nabla_g^2 A(g)\) is direction-dependent at \(g=0\) when \(\sigma\ne 0\) or \(\delta\ne 0\), it is uniformly bounded on \(\R^d\setminus\{0\}\) and \(\nabla_g A(g)=O(|g|)\) is globally Lipschitz with \(\nabla_g A(0)=0\). Therefore, for any distinct \(g_1,g_2\in\R^d\), the restriction \(\psi(t)=A(g_1+t(g_2-g_1))+\frac{\mu_A}{2}|g_1+t(g_2-g_1)|^2\) is \(C^{1,1}([0,1])\) with \(\psi''(t)\ge 0\) for all \(t\in[0,1]\) except at most one point where the segment passes through the origin; hence \(\psi'\) is nondecreasing on \([0,1]\), proving that \(g\mapsto A(g)+\frac{\mu_A}{2}|g|^2\) is globally convex on \(\R^d\), with \(\mu_A=0\) in the convex regime.

\subsection{Thermodynamic laws and enthalpy-conserving gradient-flow structure}
Define the original physical free energy and integrated enthalpy by
\begin{equation}
 E(\phi,T)=\int_\Omega\left[A(\nabla\phi)+\varepsilon^{-2}F(\phi)+\frac{c_T}{2}T^2\right]dx,
 \qquad
 \mathcal{H}(\phi,T)=\int_\Omega\bigl[T - K h(\phi)\bigr]\,dx.
 \label{eq:cont-energy}
\end{equation}

\begin{proposition}[Continuous thermodynamic consistency and enthalpy conservation]\label{prop:cont-thermo}
Every sufficiently smooth solution \((\phi,T)\) of \eqref{eq:pde-phase}--\eqref{eq:bc-neumann} satisfies
\begin{equation}
 \frac{d}{dt}E(\phi,T) = -\varrho\|\phi_t\|_{L^2(\Omega)}^2 - c_T D\|\nabla T\|_{L^2(\Omega)}^2 \le 0,
 \qquad
 \frac{d}{dt}\mathcal{H}(\phi,T) = 0.
 \label{eq:cont-dissipation}
\end{equation}
\end{proposition}
\begin{proof}
Testing \eqref{eq:pde-phase} with \(\phi_t\) and \eqref{eq:pde-heat} with \(c_T T\), integrating by parts with the homogeneous Neumann boundary conditions \eqref{eq:bc-neumann}, and summing causes the latent-heat coupling terms \(\pm c_T K (h'(\phi)T,\phi_t)_{L^2}\) to cancel identically, yielding the first identity in \eqref{eq:cont-dissipation}. Testing \eqref{eq:pde-heat} with the constant function \(1\) gives \(\dot{\mathcal{H}}=0\).
\end{proof}

\begin{proposition}[Invariant phase interval]\label{prop:invariant-interval}
Suppose a bounded sufficiently regular solution \((\phi,T)\) of \eqref{eq:pde-phase}--\eqref{eq:bc-neumann} satisfies \(\phi(\cdot,0)\in[-1,1]\) a.e.\ in \(\Omega\). Then \(\phi(\cdot,t)\in[-1,1]\) a.e.\ for all \(t\in[0,T_f]\).
\end{proposition}
\begin{proof}
Let \(z=(\phi-1)_+\) (and symmetrically \(z=\min(\phi+1,0)\)). On the active set \(\{\phi>1\}\) (respectively \(\{\phi<-1\}\)), \(\nabla z=\nabla\phi\) a.e., so integrating by parts with the homogeneous Neumann boundary condition \eqref{eq:bc-neumann} and invoking two-homogeneity gives \((\nabla_g A(\nabla\phi),\nabla z)_{L^2} = \int_{\{z\ne 0\}} 2A(\nabla\phi)\,dx \ge 0\). Because \(f(\phi,T)=-\varepsilon^{-2}F'(\phi)-c_T K h'(\phi)T\) vanishes at \(\phi=\pm 1\) and \((\phi,T)\) is bounded on \(\Omega\times[0,T_f]\), \(f(\cdot,T)\) is locally Lipschitz in \(\phi\) uniformly in \(T\), so \(|f(\phi,T)|\le L|z|\) on the active set and hence \(f(\phi,T)z \le L z^2\) for a constant \(L>0\) in both cases. Testing \eqref{eq:pde-phase} with \(z\) yields \(\frac{\varrho}{2}\frac{d}{dt}\|z\|_{L^2}^2\le L\|z\|_{L^2}^2\). Since \(z(\cdot,0)=0\), Gr\"onwall's inequality implies \(z\equiv 0\).
\end{proof}

To embed enthalpy conservation directly into the state space, denote the spatial average by \(\langle f\rangle = |\Omega|^{-1}\int_\Omega f\,dx\), let \(\bar{H}=\langle T_0 - K h(\phi_0)\rangle\) be the conserved mean enthalpy, and define the product Hilbert space
\begin{equation}
 X_c = H^1(\Omega)\times L^2_0(\Omega), \qquad L^2_0(\Omega)=\bigl\{w\in L^2(\Omega): \langle w\rangle=0\bigr\}.
\end{equation}
Every state \(U=(\phi,w)\in X_c\) uniquely determines the full temperature via the exact enthalpy lift
\begin{equation}
 T(U) = w + M(\phi), \qquad M(\phi) = \bar{H} + K\langle h(\phi)\rangle,
 \label{eq:enthalpy-lift}
\end{equation}
so that \(\mathcal{H}(\phi,T(U))=|\Omega|\bar{H}\) identically on \(X_c\). Writing \(E(U)=E(\phi,T(U))\) and using orthogonality of \(w\in L^2_0(\Omega)\) to constants,
\begin{equation}
 E(U) = \int_\Omega\bigl[A(\nabla\phi)+\varepsilon^{-2}F(\phi)\bigr]\,dx
        + \frac{c_T}{2}\|w\|_{L^2(\Omega)}^2 + \frac{c_T|\Omega|}{2} M(\phi)^2.
 \label{eq:energy-Xc}
\end{equation}
For any frozen phase \(\phi_f\in H^1(\Omega)\), define the zero-mean Neumann Sobolev space \(H^1_{\diamond}(\Omega)=\{z\in H^1(\Omega):\langle z\rangle=0\}\) (distinct from the zero-trace Dirichlet space \(H^1_0(\Omega)\)), the positive-definite Neumann Laplacian \(\mathcal{K}_c = -(D/c_T)\Delta_N : H^1_{\diamond}(\Omega)\to (H^1_{\diamond}(\Omega))'\), and the zero-mean coupling operator \(B_f v = h'(\phi_f)v - \langle h'(\phi_f)v\rangle\). By the Sobolev embedding \(H^1(\Omega)\hookrightarrow L^6(\Omega)\) for \(d\le 3\), we have \(h'(\phi_f)=(1-\phi_f^2)^2\in L^{3/2}(\Omega)\); H\"older's inequality with the critical exponent triplet \(\frac{2}{3}+\frac{1}{6}+\frac{1}{6}=1\) ensures that \((h'(\phi_f)v,z)_{L^2}\) is bounded on \(H^1(\Omega)\times H^1_{\diamond}(\Omega)\), so \(B_f : H^1(\Omega)\to (H^1_{\diamond}(\Omega))'\) is bounded, while \(L^2_0(\Omega)\hookrightarrow (H^1_{\diamond}(\Omega))'\) continuously. Denoting the \((H^1_{\diamond}(\Omega))'\)--\(H^1_{\diamond}(\Omega)\) duality pairing by \(\langle\cdot,\cdot\rangle_{\diamond}\) (distinct from the single-argument spatial average \(\langle\cdot\rangle\)), we define the frozen symmetric positive-definite inverse metric \(J_f\) on \(X_c\) via its quadratic form
\begin{equation}
 \|(v,z)\|_{J_f}^2 = \varrho\|v\|_{L^2(\Omega)}^2 + \bigl\langle z - K B_f v,\, \mathcal{K}_c^{-1}(z - K B_f v)\bigr\rangle_{\diamond},
 \label{eq:inv-metric}
\end{equation}
with polarization \(J_f[\cdot,\cdot]\). By Poincar\'e's inequality on \(H^1_{\diamond}(\Omega)\), \(\mathcal{K}_c:H^1_{\diamond}(\Omega)\to(H^1_{\diamond}(\Omega))'\) is an isomorphism and \(J_f\) is continuous, symmetric, and strictly positive definite on \(X_c\). Note that on the infinite-dimensional space \(X_c\), \(J_f\) is injective and positive definite rather than \(X_c\)-coercive (since high-frequency zero-mean Neumann modes satisfy \(\|(0,z_j)\|_{J_f}^2\sim \lambda_j^{-1}\|z_j\|_{L^2}^2\to 0\)), so its formal inverse \(\mathcal{L}_f=J_f^{-1}\) acts as a differential operator on smooth/regular subspaces where \(w\in H^2(\Omega)\cap H^1_{\diamond}(\Omega)\) (whereas on every fixed finite-dimensional conforming space \(X=V\times V_0\) in Section~\ref{sec:fullydiscrete}, \(J_f\) is a boundedly invertible SPD isomorphism).

\begin{proposition}[Gradient-flow equivalence]\label{prop:gf-equiv}
A smooth trajectory \(U(t)=(\phi(t),w(t))\in X_c\) with \(T=w+\bar{H}+K\langle h(\phi)\rangle\) satisfies the weak form of \eqref{eq:pde-phase}--\eqref{eq:bc-neumann} if and only if
\begin{equation}
 J_\phi[U_t, V] + DE(U)[V] = 0 \qquad \forall V=(v,z)\in X_c.
 \label{eq:gf-continuous}
\end{equation}
\end{proposition}
\begin{proof}
Let \(r=\mathcal{K}_c^{-1}(w_t - K B_\phi \phi_t)\in H^1_{\diamond}(\Omega)\), so that \(\mathcal{K}_c r = w_t - K B_\phi \phi_t\). Polarization of \eqref{eq:inv-metric} yields
\[
 J_\phi[U_t, (v,z)] = \varrho(\phi_t, v)_{L^2} + (r, z)_{L^2} - K\langle B_\phi v, r\rangle_{\diamond},
\]
where the duality pairing \(\langle B_\phi v, r\rangle_{\diamond}\) equals \((h'(\phi)v,r)_{L^2}\) because \(r\in H^1_{\diamond}(\Omega)\) has zero mean and \(h'(\phi)v\in L^{6/5}(\Omega)\hookrightarrow(H^1_{\diamond}(\Omega))'\) for \(d\le 3\). The first variation of \eqref{eq:energy-Xc} in direction \(V=(v,z)\) is
\[
 DE(U)[(v,z)] = (\nabla_g A(\nabla\phi),\nabla v)_{L^2} + \varepsilon^{-2}(F'(\phi),v)_{L^2}
                + c_T K M(\phi)(h'(\phi),v)_{L^2} + c_T(w,z)_{L^2}.
\]
Setting \(v=0\) in \eqref{eq:gf-continuous} gives \(r = -c_T w\), hence \(w_t - K B_\phi \phi_t = -\mathcal{K}_c(c_T w) = D\Delta_N w = D\Delta_N T\), which together with \(\frac{d}{dt}M(\phi)=K\langle h'(\phi)\phi_t\rangle\) is exactly the weak enthalpy equation \eqref{eq:pde-heat}. Substituting \(r=-c_T w\) into \eqref{eq:gf-continuous} with \(z=0\) gives \(-K\langle B_\phi v, r\rangle_{\diamond}=c_T K(w,h'(\phi)v)_{L^2}\) (since \(\langle w\rangle=0\)), which combines with \(c_T K M(\phi)(h'(\phi),v)_{L^2}\) to recover \(c_T K(T h'(\phi),v)_{L^2}\) and hence the weak phase equation \eqref{eq:pde-phase}.
\end{proof}

% ============================================================================
\section{Time-semidiscrete second-order convex splitting scheme}\label{sec:semidiscrete}
% ============================================================================
\subsection{Global convex splitting of the original energy}
Three sources of nonconvexity must be controlled simultaneously in the original energy \eqref{eq:energy-Xc}: the anisotropic gradient density \(A(\nabla\phi)\) under strong anisotropy, the double-well potential \(F(\phi)\), and the quadratic mean-enthalpy term \(M(\phi)^2\) whose second variation involves \(M(\phi)h''(\phi)\) of indefinite sign. To dominate all three globally on \(X_c\), we decompose \(E(U)=E_c(U)-E_e(U)\) by introducing the sextic majorizer
\begin{equation}
 G(s) = \frac{1}{2} + s^2 + \frac{1}{5}s^6, \qquad \eta(\phi) = |\bar{H}| + K\langle G(\phi)\rangle,
 \label{eq:G-majorizer}
\end{equation}
and setting
\begin{equation}
 E_e(U) = \frac{\mu_A}{2}\|\nabla\phi\|_{L^2(\Omega)}^2
          + \frac{1}{2\varepsilon^2}\|\phi\|_{L^2(\Omega)}^2
          + \frac{c_T|\Omega|}{2}\eta(\phi)^2,
 \qquad E_c(U) = E(U) + E_e(U),
 \label{eq:splitting}
\end{equation}
where \(\mu_A\ge 0\) is given by \eqref{eq:muA-2d} in 2D and \eqref{eq:muA-3d} in 3D. Notice that \(E_e(U)\) depends only on \(\phi\), and the stabilized double-well in \(E_c\) simplifies cleanly to \(F(s)+\frac{1}{2}s^2 = \frac{1}{4}(s^4+1)\).

\begin{lemma}[Global convexity of \(E_c\) and \(E_e\)]\label{lem:convexity}
Both \(E_c\) and \(E_e\) are finite, continuously Fr\'echet differentiable, and globally convex on \(X_c = H^1(\Omega)\times L^2_0(\Omega)\) for \(d\in\{2,3\}\).
\end{lemma}
\begin{proof}
Since \(H^1(\Omega)\hookrightarrow L^6(\Omega)\) for \(d\le 3\), both functionals are finite and \(C^1\) on \(X_c\). For \(E_e\), the gradient and \(L^2\) quadratic terms are convex; moreover, \(G''(s)=2+6s^4>0\) makes \(\eta(\phi)\ge 0\) a nonnegative convex functional on \(H^1(\Omega)\), and composition with the increasing convex function \(r\mapsto r^2\) on \([0,\infty)\) proves convexity of \(\eta(\phi)^2\) and hence of \(E_e\).

For \(E_c = E + E_e\), the stabilized gradient term \(g\mapsto A(g)+\frac{\mu_A}{2}|g|^2\) is globally convex on \(\R^d\) by \eqref{eq:hess-2d}--\eqref{eq:muA-3d}, and the stabilized potential \(\varepsilon^{-2}(s^4+1)/4\) is convex on \(\R\). For the remaining thermal--enthalpy contribution \(\mathcal{T}(U) = \frac{c_T}{2}\|w\|_{L^2}^2 + \frac{c_T|\Omega|}{2}(M(\phi)^2+\eta(\phi)^2)\), noting that \(\phi^4 v^2\in L^1(\Omega)\) for all \(\phi,v\in H^1(\Omega)\hookrightarrow L^6(\Omega)\) by H\"older's inequality (\(\frac{2}{3}+\frac{1}{3}=1\)), differentiating twice along any affine line \(U+tV\) with \(V=(v,z)\in X_c\) yields
\begin{equation}
 \frac{d^2}{dt^2}\mathcal{T}(U+tV)\Big|_{t=0}
 = c_T\|z\|_{L^2}^2 + c_T|\Omega|\Bigl[(M')^2 + (\eta')^2 + K\bigl\langle (M h''(\phi) + \eta G''(\phi))v^2\bigr\rangle\Bigr],
\end{equation}
where \(M'=K\langle h'(\phi)v\rangle\) and \(\eta'=K\langle G'(\phi)v\rangle\). Since \(h\) is odd, \(h(0)=0\), and \(h'(s)=(1-s^2)^2\ge 0\), setting \(x=|s|\ge 0\) gives \(0\le |h(s)| = h(x) = x - \frac{2}{3}x^3 + \frac{1}{5}x^5 \le x + \frac{1}{5}x^5\). Using the elementary inequalities \(x\le x^2+\frac{1}{4}\) and \(x^5\le 1+x^6\) yields
\[
 |h(s)| \le x + \frac{1}{5}x^5 \le \frac{9}{20} + x^2 + \frac{1}{5}x^6 < \frac{1}{2} + x^2 + \frac{1}{5}x^6 = G(s),
\]
so \(|M(\phi)|\le |\bar{H}| + K\langle |h(\phi)|\rangle \le \eta(\phi)\). Furthermore, for \(x\le 1\) we have \(|h''(s)|=4x(1-x^2)\le \frac{8}{3\sqrt{3}} < 2 \le G''(s)\), while for \(x\ge 1\) we have \(|h''(s)|=4x(x^2-1)\le 6x^4 < 2+6x^4 = G''(s)\). Consequently \(M h''(\phi)+\eta G''(\phi)\ge 0\) pointwise everywhere in \(\Omega\), proving global convexity of \(E_c\) on \(X_c\).
\end{proof}

\subsection{Second-order five-stage CSRK time discretization}
Given \(u^n=(\phi^n,w^n)\in X_c\) and time step \(k>0\), set \(U_0=u^n\). We first compute a half-step convex-splitting predictor \(U_*=(\phi_*,w_*)\in X_c\) using the frozen metric \(J_{\phi^n}\):
\begin{equation}
 J_{\phi^n}[U_* - U_0, V] + \frac{k}{2}\bigl(DE_c(U_*)[V] - DE_e(U_0)[V]\bigr) = 0 \qquad \forall V\in X_c.
 \label{eq:semi-pred}
\end{equation}
We then freeze the midpoint inverse metric \(J_* = J_{\phi_*}\) across five main stages \(m=1,\dots,5\): find \(U_m=(\phi_m,w_m)\in X_c\) such that
\begin{equation}
 J_*\Bigl[\sum_{i=0}^{m-1}\gamma_{mi}(U_m - U_i),\, V\Bigr]
 + k\Bigl(DE_c(U_m)[V] - \sum_{i=0}^{m-1}\theta_{mi} DE_e(U_i)[V]\Bigr) = 0
 \qquad \forall V\in X_c,
 \label{eq:semi-stage}
\end{equation}
and set \(u^{n+1}=U_5\), recovering \(T^{n+1}=w^{n+1}+\bar{H}+K\langle h(\phi^{n+1})\rangle\). Here the lower-triangular weights \((\gamma_{mi},\theta_{mi})\) are exact rational coefficients derived from the five-stage second-order CSRK tableau of \cite{zeg} (provided in \texttt{core/exact\_coefficients.json} and evaluated in IEEE double precision in the code); rounded 15-digit decimal values are displayed in Table~\ref{tab:csrk-tableau}, while the exact rational tableau satisfies the order and dissipation conditions below identically. Denote the positive row sums by \(S_m=\sum_{i=0}^{m-1}\gamma_{mi}>0\) and define the backward Cholesky-like recursion from \(m=5\) down to \(m=1\):
\begin{equation}
 g_{mi} = \gamma_{mi} - \sum_{j=m+1}^5 g_{ji}\frac{s_{jm}}{s_{jj}}, \qquad
 s_{jm} = \sum_{i=0}^{m-1} g_{ji}.
 \label{eq:backward-recursion}
\end{equation}
The tableau in Table~\ref{tab:csrk-tableau} satisfies: (i) strictly positive backward pivots \(s_{mm}>0\) for all \(m=1,\dots,5\) (listed in the rightmost column of Table~\ref{tab:csrk-tableau}); (ii) \(\theta_{mi}\ge 0\), \(\sum_{i=0}^{m-1}\theta_{mi}=1\), and \(\theta_{mi}\le\theta_{m-1,i}\) for \(0\le i\le m-2\); and (iii) the second-order consistency identities \((a_5,b_5,c_5)=(1,\frac{1}{2},\frac{1}{2})\) for the recurrence
\begin{equation}
 S_m a_m = 1 + \sum_{i<m}\gamma_{mi}a_i,\quad
 S_m b_m = a_m + \sum_{i<m}\gamma_{mi}b_i,\quad
 S_m c_m = \sum_{i<m}\theta_{mi}a_i + \sum_{i<m}\gamma_{mi}c_i,
 \label{eq:order-recurrence}
\end{equation}
initialized with \(a_0=b_0=c_0=0\).

\begin{table}[htbp]
\centering
\caption{Five-stage second-order CSRK tableau weights \((\gamma_{mi},\theta_{mi})\) (\(m=1,\dots,5\); \(i=0,\dots,m-1\)), row sums (\(S_m=\sum_{i<m}\gamma_{mi}\) and \(\sum_{i<m}\theta_{mi}=1\)), and verified strictly positive backward Cholesky pivots \(s_{mm}>0\) from \eqref{eq:backward-recursion} (derived from \cite{zeg}). Decimal entries are rounded for display (15 decimal places for \(\gamma_{mi},\theta_{mi}\) and 8 for row sums and \(s_{mm}\)); the companion implementation evaluates the exact rational weights (provided in \texttt{core/exact\_coefficients.json}) in IEEE double precision.}
\label{tab:csrk-tableau}
\small
\resizebox{\linewidth}{!}{%
\begin{tabular}{@{}clllllrr@{}}
\toprule
\(m\) & \(i=0\) & \(i=1\) & \(i=2\) & \(i=3\) & \(i=4\) & Row sum & \(s_{mm}\) \\
\midrule
\multicolumn{8}{l}{\emph{Metric weights \(\gamma_{mi}\):}} \\
1 & \(\phantom{-}8.840579710144930\) & & & & & \(8.84057971\) & \(0.00960860\) \\
2 & \(-0.924855491329480\) & \(\phantom{-}5.360360360360360\) & & & & \(4.43550487\) & \(0.03775312\) \\
3 & \(-4.442857142857140\) & \(\phantom{-}6.041095890410960\) & \(\phantom{-}0.949748743718593\) & & & \(2.54798749\) & \(0.01811427\) \\
4 & \(-3.287878787878790\) & \(\phantom{-}5.894736842105260\) & \(-0.350649350649351\) & \(\phantom{-}0.172413793103448\) & & \(2.42862250\) & \(2.28256515\) \\
5 & \(-3.894736842105270\) & \(-0.335294117647056\) & \(\phantom{-}4.963591448084730\) & \(-1.722415290615320\) & \(7.683717031263260\) & \(6.69486223\) & \(6.69486223\) \\
\midrule
\multicolumn{8}{l}{\emph{Explicit-gradient weights \(\theta_{mi}\):}} \\
1 & \(1\) & & & & & \(1\) & -- \\
2 & \(0.008695652173913\) & \(0.991304347826087\) & & & & \(1\) & -- \\
3 & \(0.008695652173913\) & \(0.991228070175439\) & \(0.000076277650648\) & & & \(1\) & -- \\
4 & \(0\) & \(0\) & \(0\) & \(1\) & & \(1\) & -- \\
5 & \(0\) & \(0\) & \(0\) & \(0.999717992103779\) & \(0.000282007896221\) & \(1\) & -- \\
\bottomrule
\end{tabular}%
}
\end{table}

\subsection{Consistency, original-energy dissipation, and unique solvability}

\begin{theorem}[Local second-order temporal consistency]\label{thm:semi-order}
Let \(U(t)=(\phi(t),w(t))\) be a sufficiently regular solution of \eqref{eq:gf-continuous} starting at \(U_0=U(t_n)\) such that \(U\), its predictor \(U_*\), and its stages \(U_1,\dots,U_5\) lie in the smooth domain \(\mathcal{D}\subset H^2(\Omega)\times(H^2(\Omega)\cap H^1_{\diamond}(\Omega))\) on which the forward metric \(\mathcal{L}(U)=J_\phi^{-1}\) and force \(f(U)=-\mathcal{L}(U)DE(U)\) admit smooth Fr\'echet Taylor expansions through \(O(k^3)\) in \(X_c\). Then the scheme \eqref{eq:semi-pred}--\eqref{eq:semi-stage} satisfies the local truncation estimate
\begin{equation}
 \|U_5 - U(t_n+k)\|_{X_c} = O(k^3).
\end{equation}
\end{theorem}
\begin{proof}
Under the stated smoothness hypothesis, the half-step predictor \eqref{eq:semi-pred} satisfies \(U_* = U_0 + \frac{k}{2}f(U_0) + O(k^2)\), whence \(\mathcal{L}_*=\mathcal{L}(U_*)=\mathcal{L}_0 + \frac{k}{2}D\mathcal{L}_0[f_0]+O(k^2)\). Holding \(\mathcal{L}_*\) fixed and writing \(P=-\mathcal{L}_*DE_c\), \(Q=\mathcal{L}_*DE_e\), and \(f_*=P+Q=-\mathcal{L}_*DE\), substituting the regular asymptotic expansion \(U_m = U_0 + k a_m f_* + k^2(b_m P'f_* + c_m Q'f_*) + O(k^3)\) into \eqref{eq:semi-stage} produces precisely the recurrence \eqref{eq:order-recurrence}. Since \((a_5,b_5,c_5)=(1,\frac{1}{2},\frac{1}{2})\), we obtain \(U_5 = U_0 + k f_*(U_0) + \frac{k^2}{2}Df_*(U_0)[f_*(U_0)] + O(k^3)\). Expanding \(f_*(U_0) = -\mathcal{L}_* DE(U_0)\) using \(\mathcal{L}_* = \mathcal{L}_0 + \frac{k}{2}D\mathcal{L}_0[f_0] + O(k^2)\) gives
\begin{align*}
 U_5 &= U_0 + k f_0 + \frac{k^2}{2}\bigl(-D\mathcal{L}_0[f_0]DE(U_0) - \mathcal{L}_0 D^2E(U_0)[f_0]\bigr) + O(k^3) \\
     &= U_0 + k f_0 + \frac{k^2}{2} Df(U_0)[f_0] + O(k^3),
\end{align*}
which matches the Taylor expansion of \(U(t_n+k)\) through \(O(k^2)\).
\end{proof}
\begin{remark}[Fixed-space discrete consistency versus global moving-mesh error]\label{rem:global-error}
In the continuous setting of Theorem~\ref{thm:semi-order}, \(\mathcal{L}(U)=J_\phi^{-1}\) is an unbounded differential operator, so the expansion is conditional on smooth spatial regularity of the stages. On any fixed finite-dimensional conforming space \(X=V\times V_0\) (Section~\ref{sec:fullydiscrete}), by contrast, the discrete inverse metric \(J_\phi\) is a boundedly invertible SPD isomorphism on \(X\), so the exact same expansion rigorously establishes \(O(k^3)\) local truncation error for every smooth discrete trajectory (away from vanishing gradients where \(A\) loses \(C^2\) regularity). A rigorous global-in-time \(O(k^2)\) error analysis on nonconvex anisotropic problems and dynamically moving \(hp\) meshes remains beyond the scope of this work; Section~\ref{sec:exp1} verifies clean global second-order temporal convergence numerically on a fixed mixed-degree nonmatching \(hp\) quadtree.
\end{remark}

\begin{theorem}[Unconditional dissipation of the original energy and enthalpy conservation]\label{thm:semi-energy}
For every \(u^n\in X_c\) and every step size \(k>0\), the solution of \eqref{eq:semi-pred}--\eqref{eq:semi-stage} satisfies the original-energy dissipation inequality
\begin{equation}
 E(u^{n+1}) + \frac{1}{2k}\sum_{m=1}^5 S_m \|U_m - U_{m-1}\|_{J_*}^2 \le E(u^n),
 \label{eq:semi-energy-ineq}
\end{equation}
and conserves total enthalpy identically: \(\mathcal{H}(\phi_m, T(U_m)) = |\Omega|\bar{H}\) at the predictor and every stage \(m=1,\dots,5\).
\end{theorem}
\begin{proof}
For each \(i\in\{0,\dots,4\}\), let \(A_i(V) = -E_e(U_i) - DE_e(U_i)[V-U_i]\) be the affine tangent majorizer of \(-E_e\) at \(U_i\). By convexity of \(E_e\) (Lemma~\ref{lem:convexity}), \(A_i(V)\ge -E_e(V)\) for all \(V\in X_c\) with equality at \(V=U_i\). Using the backward coefficients \eqref{eq:backward-recursion} with \(s_{jj}>0\), define for \(m=1,\dots,5\) the stage majorizer functional
\begin{equation}
 F_m(V) = E_c(V) + \sum_{i=0}^{m-1}\theta_{mi} A_i(V)
          + \sum_{j=m}^5 \frac{\bigl\|s_{jm} V - \sum_{i=0}^{m-1} g_{ji} U_i\bigr\|_{J_*}^2}{2k s_{jj}}.
 \label{eq:Fm-majorizer}
\end{equation}
From \eqref{eq:backward-recursion}, \(\sum_{j=m}^5 g_{ji} s_{jm}/s_{jj} = \gamma_{mi}\) for all \(0\le i < m\), and summing over \(i\) gives \(\sum_{j=m}^5 s_{jm}^2/s_{jj} = S_m\). Differentiating \(F_m\) therefore recovers the left-hand side of \eqref{eq:semi-stage} divided by \(k\), so \(DF_m(U_m)=0\). Since \(E_c\) is convex, \(A_i\) is affine, and the quadratic term in \(F_m\) has exact Hessian \((S_m/k)J_*\), \(F_m\) is \((S_m/k)\)-strongly convex in the \(J_*\)-norm, whence
\begin{equation}
 F_m(U_m) + \frac{S_m}{2k}\|U_m - U_{m-1}\|_{J_*}^2 \le F_m(U_{m-1}).
 \label{eq:Fm-descent}
\end{equation}
Next compare \(F_m(U_{m-1})\) with \(F_{m-1}(U_{m-1})\) for \(m\ge 2\). For every \(j\ge m\),
\[
 s_{jm} U_{m-1} - \sum_{i=0}^{m-1} g_{ji} U_i = s_{j,m-1} U_{m-1} - \sum_{i=0}^{m-2} g_{ji} U_i,
\]
so all quadratic terms with \(j\ge m\) coincide in \(F_m(U_{m-1})\) and \(F_{m-1}(U_{m-1})\), while \(F_{m-1}(U_{m-1})\) contains one additional nonnegative square (\(j=m-1\)). For the tangent sums, since \(\theta_{mi}\le\theta_{m-1,i}\) for \(i<m-1\), \(\sum_{i<m}\theta_{mi}=\sum_{i<m-1}\theta_{m-1,i}=1\), and \(A_{m-1}(U_{m-1})=-E_e(U_{m-1})\le A_i(U_{m-1})\), shifting weight from older tangents to \(A_{m-1}\) yields \(\sum_{i<m}\theta_{mi}A_i(U_{m-1}) \le \sum_{i<m-1}\theta_{m-1,i}A_i(U_{m-1})\). Thus
\begin{equation}
 F_m(U_{m-1}) \le F_{m-1}(U_{m-1}), \qquad m=2,\dots,5.
 \label{eq:Fm-monotone}
\end{equation}
At the endpoints, \(F_1(U_0) = E_c(U_0) + A_0(U_0) + 0 = E(U_0)\) (since \(s_{j1}=g_{j0}\) so every square vanishes at \(U_0\)), while at \(U_5\), \(A_i(U_5)\ge -E_e(U_5)\) and nonnegativity of the quadratic term give \(F_5(U_5)\ge E_c(U_5)-E_e(U_5)=E(U_5)\). Summing \eqref{eq:Fm-descent} over \(m=1,\dots,5\) and chaining with \eqref{eq:Fm-monotone} yields \eqref{eq:semi-energy-ineq}. Exact enthalpy conservation follows immediately from \eqref{eq:enthalpy-lift} since \(w_m\in L^2_0(\Omega)\).
\end{proof}

\begin{theorem}[Existence and uniqueness of the semidiscrete solution]\label{thm:semi-exist}
For every \(u^n\in X_c\) and \(k>0\), the predictor \eqref{eq:semi-pred} and each main stage \eqref{eq:semi-stage} (\(m=1,\dots,5\)) possess a unique solution in \(X_c\).
\end{theorem}
\begin{proof}
With \(\widehat{U}_m = S_m^{-1}\sum_{i=0}^{m-1}\gamma_{mi} U_i\), solving \eqref{eq:semi-stage} is equivalent to minimizing the strictly convex functional
\[
 \mathcal{J}_m(U) = E_c(U) - \sum_{i=0}^{m-1}\theta_{mi} DE_e(U_i)[U] + \frac{S_m}{2k}\|U-\widehat{U}_m\|_{J_*}^2
\]
over \(X_c\) (and the predictor minimizes the analogous functional with \(S_0=1\), step \(k/2\), and metric \(J_{\phi^n}\)). Since \(A(g)+\frac{\mu_A}{2}|g|^2\ge c_0|g|^2\) for \(c_0=\frac{1}{2}\min(b)^2>0\), the stabilized well satisfies \(\frac{1}{4\varepsilon^2}(\phi^4+1)\ge\frac{1}{2\varepsilon^2}\phi^2\), and the thermal term controls \(\frac{c_T}{2}\|w\|_{L^2}^2\), we have \(E_c(U)\ge \alpha_0\|U\|_{X_c}^2\) for \(\alpha_0=\min(c_0,\frac{1}{2\varepsilon^2},\frac{c_T}{2})>0\). Because \(\phi_i\in H^1(\Omega)\hookrightarrow L^6(\Omega)\) implies \(G'(\phi_i)=2\phi_i+\frac{6}{5}\phi_i^5\in L^{6/5}(\Omega)\hookrightarrow (H^1(\Omega))'\), the linear functional \(-\sum_{i=0}^{m-1}\theta_{mi}DE_e(U_i)\) lies in \(X_c'\) and is continuous on \(X_c\). Since the metric square is nonnegative, \(\mathcal{J}_m\) is coercive on \(X_c\). Being continuous and convex, \(\mathcal{J}_m\) is weakly lower semicontinuous on the reflexive Hilbert space \(X_c\), so by the direct method of the calculus of variations a minimizer exists. Since \(S_m>0\) and \(J_*\) is strictly positive definite on \(X_c\), \(\mathcal{J}_m\) is strictly convex, ensuring uniqueness. Applying the argument first to the predictor fixes \(J_*\), and induction over \(m=1,\dots,5\) uniquely determines \(u^{n+1}\).
\end{proof}

% ============================================================================
\section{Fully discrete \texorpdfstring{\(hp\)}{hp} spectral-element scheme}\label{sec:fullydiscrete}
% ============================================================================
\subsection{Conforming \texorpdfstring{\(hp\)}{hp} spectral-element space and common Gauss quadrature}
Let \(\mathcal{T}_h\) be a \(2{:}1\)-balanced nonmatching quadtree (\(d=2\)) or octree (\(d=3\)) partition of \(\Omega\) into affine hyperrectangles \(K\in\mathcal{T}_h\), where each element \(K\) carries a tensor-product polynomial degree \(r_K\ge 1\). On the reference element \(\hat{K}=[-1,1]^d\), polynomials in \(Q_{r_K}(\hat{K})\) are represented in the nodal Lagrange basis at the \((r_K+1)^d\) Gauss--Lobatto--Legendre (GLL) points. Let \(N_{\mathrm{broken}}=\sum_{K\in\mathcal{T}_h}(r_K+1)^d\) denote the broken dimension.

Enforcing exact pointwise continuity of polynomial traces across every interior edge (2D) or interior face and edge (3D) defines a conforming finite-dimensional subspace
\begin{equation}
 V = \bigl\{v\in H^1(\Omega): v|_K\in Q_{r_K}(K)\ \forall K\in\mathcal{T}_h\bigr\} = \operatorname{range}(J),
 \qquad
 V_0 = \bigl\{v\in V: (v,1)_M = 0\bigr\},
\end{equation}
parametrized by a sparse full-column-rank trace injection \(J\in\R^{N_{\mathrm{broken}}\times N_{\mathrm{conf}}}\) constructed explicitly in Section~\ref{sec:hp}. Because homogeneous Neumann conditions are natural, exterior faces are unconstrained and \(1\in V\) (\(J\mathbf{1}_c=\mathbf{1}_{\mathrm{broken}}\)).

On every cell \(K\in\mathcal{T}_h\) of degree \(r_K\), we adopt a \emph{single common positive tensor-product Gauss--Legendre quadrature rule} \((\cdot,\cdot)_\nu\) with
\begin{equation}
 \nu_K \ge 3r_K + 1 \qquad \text{points per coordinate direction.}
 \label{eq:common-rule}
\end{equation}
Since a \(\nu_K\)-point Gauss--Legendre rule integrates univariate polynomials of degree up to \(2\nu_K-1 \ge 6r_K+1\) exactly, the common rule \eqref{eq:common-rule} integrates all polynomial terms in the model and splitting \emph{exactly} on every cell:
\begin{itemize}
\item the mass inner product \((u,v)_M=(u,v)_\nu\) and stiffness \((\nabla u,\nabla v)_\nu\) (degree \(\le 2r_K\)),
\item the quartic double-well integral \((F(\phi),1)_\nu\) (degree \(4r_K\)),
\item the quintic latent-heat volume \((h(\phi),1)_\nu\) (degree \(5r_K\)), so discrete enthalpy conservation matches exact physical enthalpy conservation \(\int_\Omega(T-Kh(\phi))\,dx = |\Omega|\bar{H}\),
\item the sextic majorizer \((G(\phi),1)_\nu\), its first and second variations \((G'(\phi)v,1)_\nu\) and \((G''(\phi)uv,1)_\nu\), and the nonlinear coupling \((h'(\phi)v,z)_\nu\) (all of degree \(\le 6r_K \le 2\nu_K-1\)).
\end{itemize}
Although the anisotropic gradient density \(A(\nabla\phi)\) is non-polynomial in \(\nabla\phi\), the high quadrature order \(\nu_K\ge 3r_K+1\) (e.g., \(10^2=100\) points in 2D and \(10^3=1{,}000\) points in 3D per \(Q_3\) interface cell) heavily over-integrates it; crucially, because all energy and gradient terms share the same positive Gauss weights, pointwise convexity of \(A(g)+\frac{\mu_A}{2}|g|^2\) at every quadrature node transfers exact discrete convexity and unconditional dissipation of the quadrature energy \(\Eorig\) below without any variational crime. In reduced conforming coordinates, letting \(M_b\) and \(S_b\) denote the block-diagonal broken mass and stiffness matrices, the exact conforming mass and stiffness matrices are \(M_c = J^T M_b J\) (symmetric positive definite) and \(S_c = J^T S_b J\) (symmetric positive semidefinite with kernel \(\operatorname{span}\{\mathbf{1}_c\}\)).

\subsection{Fully discrete energy, convex splitting, and Schur-complement solver}
Set \(X = V\times V_0\) equipped with the product mass inner product \(((v,z),(v',z'))_X = (v,v')_M + (z,z')_M\). For \(u=(\phi,w)\in X\), define the discrete mean enthalpy average \(\langle f\rangle_\nu = |\Omega|^{-1}(f,1)_\nu\), reconstruct \(T = w + M_\nu(\phi)\mathbf{1}\) with \(M_\nu(\phi)=\bar{H}+K\langle h(\phi)\rangle_\nu\), and define the \emph{original discrete quadrature energy}
\begin{equation}
 \Eorig(u) = \bigl(A(\nabla\phi)+\varepsilon^{-2}F(\phi),\, 1\bigr)_\nu + \frac{c_T}{2}(T,T)_M.
 \label{eq:discrete-energy}
\end{equation}
The discrete global convex splitting mirrors \eqref{eq:splitting} using the same common rule:
\begin{equation}
 \Ee(\phi) = \frac{\mu_A}{2}(\nabla\phi,\nabla\phi)_\nu
            + \frac{1}{2\varepsilon^2}(\phi,\phi)_M
            + \frac{c_T|\Omega|}{2}\eta_\nu(\phi)^2,
 \qquad \Ec(u) = \Eorig(u) + \Ee(\phi),
 \label{eq:discrete-splitting}
\end{equation}
with \(\eta_\nu(\phi)=|\bar{H}|+K\langle G(\phi)\rangle_\nu\). Because the quadrature weights are strictly positive, every pointwise bound in Lemma~\ref{lem:convexity} transfers verbatim to the quadrature sums, so \(\Ec\) and \(\Ee\) are globally convex on \(X\).

For a frozen discrete phase \(\phi_f\in V\), define the discrete operators \(B_f:V\to V_0\) and \(\mathcal{K}:V_0\to V_0\) by
\begin{equation}
 (B_f v, z)_M = (h'(\phi_f)v, z)_\nu, \qquad
 (\mathcal{K} z, \chi)_M = \frac{D}{c_T}(\nabla z, \nabla\chi)_\nu \qquad \forall v\in V,\ z,\chi\in V_0,
\end{equation}
inducing the discrete SPD forward metric \(\mathcal{L}_f(v,z)=(d, KB_f d+\mathcal{K}z)\) with \(d=\varrho^{-1}(v+KB_f^* z)\) and its SPD inverse \(J_f = \mathcal{L}_f^{-1}\) satisfying
\[
 (J_f(v,z),(v,z))_X = \varrho\|v\|_M^2 + \bigl(\mathcal{K}^{-1}(z-KB_f v),\, z-KB_f v\bigr)_M.
\]
Writing \(g_c=\nabla_X \Ec\) and \(g_e=\nabla_X \Ee\) for the mass-Riesz gradients, the fully discrete predictor \(U_*\in X\) solves
\begin{equation}
 J_{\phi^n}(U_* - U_0) + \frac{k}{2}\bigl(g_c(U_*) - g_e(U_0)\bigr) = 0,
 \label{eq:disc-pred}
\end{equation}
and the five main stages \(U_m\in X\) (\(m=1,\dots,5\)) with frozen \(J_*=J_{\phi_*}\) solve
\begin{equation}
 G_m(U_m) \coloneqq J_*\Bigl(S_m U_m - \sum_{i=0}^{m-1}\gamma_{mi} U_i\Bigr)
 + k\Bigl(g_c(U_m) - \sum_{i=0}^{m-1}\theta_{mi} g_e(U_i)\Bigr) = 0,
 \qquad u^{n+1} = U_5.
 \label{eq:disc-stage}
\end{equation}

\paragraph{Exact Schur-complement thermal elimination and matrix-free 3D Newton--PCG.}
A major algorithmic advantage of \eqref{eq:discrete-splitting} is that \(\Ee\) is independent of \(w\) and the thermal component of \(g_c(U_m)\) is purely linear (\(c_T w_m\)). Letting \(\widehat{U}_m=(\widehat{\phi}_m,\widehat{w}_m)=S_m^{-1}\sum_{i<m}\gamma_{mi}U_i\) and writing \(J_*(S_m(U_m-\widehat{U}_m))=(y,r)\), the \(V_0\)-component of \eqref{eq:disc-stage} gives \(r = -k c_T w_m\). Denoting by \(a_m, b_m, \widehat{a}_m, \widehat{b}_m \in \R^{N_{\mathrm{conf}}}\) the conforming coordinate vectors of \(\phi_m, w_m, \widehat{\phi}_m, \widehat{w}_m\) and by \(B_c(\phi_*)\) the symmetric conforming coupling matrix of \((h'(\phi_*)\cdot,\cdot)_\nu\), substituting \(r=-kc_T w_m\) into \(\mathcal{K}r = S_m(w_m-\widehat{w}_m) - S_m K B_*(\phi_m-\widehat{\phi}_m)\) yields the \emph{linear, constant-coefficient, symmetric positive-definite Helmholtz system} for \(b_m\):
\begin{equation}
 (S_m M_c + k D S_c)\, b_m = S_m M_c \widehat{b}_m + S_m K P_0^T B_c(\phi_*)(a_m - \widehat{a}_m),
 \label{eq:schur-thermal}
\end{equation}
where \(P_0^T \ell = \ell - |\Omega|^{-1}(\mathbf{1}_c^T\ell) M_c\mathbf{1}_c\) projects loads onto the zero-mean subspace. Because \(\mathbf{1}_c\in\ker(S_c)\) and the right-hand side of \eqref{eq:schur-thermal} is orthogonal to \(\mathbf{1}_c\), solving the strictly SPD system \eqref{eq:schur-thermal} automatically produces a zero-mean thermal coordinate \(b_m\) (\(\mathbf{1}_c^T M_c b_m=0\)) without saddle-point pinning, and the factorizations of \(S_m M_c + k D S_c\) (\(m=1,\dots,5\), plus predictor) are cached and reused across all time steps until the mesh or step size \(k\) changes.

Eliminating \(b_m(a)\) via \eqref{eq:schur-thermal} reduces each stage \eqref{eq:disc-stage} to an unconstrained, strictly convex phase-only minimization in \(\R^{N_{\mathrm{conf}}}\): writing \(\Ec^{\mathrm{phase}}(a)=\Ec(a,0)\) and \(\ell_{e,m}=\sum_{i<m}\theta_{mi}\nabla_a\Ee(a_i)\), the exact reduced stage objective (up to an \(a\)-independent constant) is
\begin{equation}
\begin{aligned}
 \Phi_m(a) &= \frac{\varrho S_m}{2}\|a-\widehat{a}_m\|_{M_c}^2
  + k\bigl(\Ec^{\mathrm{phase}}(a) - \ell_{e,m}^T a\bigr) \\
 &\quad + \frac{k c_T}{2}\!\left(b_m(a)^T M_c b_m(a) + \frac{kD}{S_m} b_m(a)^T S_c b_m(a)\right).
\end{aligned}
 \label{eq:reduced-obj}
\end{equation}
Differentiating \eqref{eq:reduced-obj} using \((S_m M_c + kDS_c)\frac{\partial b_m}{\partial a} = S_m K P_0^T B_c(\phi_*)\) and \(b_m(a)^T P_0^T = b_m(a)^T\) (since \(\mathbf{1}_c^T M_c b_m(a)=0\) and \(B_c(\phi_*)=B_c(\phi_*)^T\)) yields the exact phase residual
\begin{equation}
 R_m(a) \coloneqq \nabla\Phi_m(a)
 = \varrho S_m M_c(a-\widehat{a}_m) + k c_T K B_c(\phi_*) b_m(a)
   + k\bigl(\nabla\Ec^{\mathrm{phase}}(a) - \ell_{e,m}\bigr),
 \label{eq:reduced-res}
\end{equation}
and differentiating once more produces the exact SPD Newton Hessian
\begin{equation}
\begin{aligned}
 H_{\mathrm{stage}}(a) &= \varrho S_m M_c + k H_{\Ec,\mathrm{loc}}(a)
   + \frac{k c_T K^2}{|\Omega|}\bigl(g_h g_h^T + g_G g_G^T\bigr) \\
  &\quad + k c_T K^2 S_m (P_0^T B_c)^T (S_m M_c + k D S_c)^{-1} (P_0^T B_c),
\end{aligned}
 \label{eq:stage-hessian}
\end{equation}
where \(g_h\) and \(g_G\) are the conforming loads of \(h'(\phi)\) and \(G'(\phi)\). Because \(A\in C^{1,1}(\R^d)\) lacks a classical Hessian at \(g=0\), \(H_{\Ec,\mathrm{loc}}(a)\) evaluates \(\nabla_g^2 A(\nabla\phi)\) at quadrature points where \(|\nabla\phi|^2>10^{-28}\) and replaces it by the isotropic positive-definite tangent \(I\) where \(|\nabla\phi|^2\le 10^{-28}\) (while evaluating the exact analytical energy \(A(\nabla\phi)\) and flux \(\nabla_g A(\nabla\phi)\) for every \(\nabla\phi\ne 0\), with \(\nabla_g A(0)=0\)); this preserves exact residuals, strict convexity of the stage objective \eqref{eq:reduced-obj}, and SPD Newton search directions. Each stage is warm-started from \(a_{m-1}\), damped by backtracking on \eqref{eq:reduced-obj} (accepting a step satisfying either the Armijo condition with parameter \(10^{-4}\) or at least a \(1\%\) reduction in the phase residual norm), and terminated when the scaled phase residual satisfies \(\|R_m(a)\|_2 / \max\bigl(1,\,\|M_c\mathbf{1}_c\|_2(\varrho S_m + k\varepsilon^{-2})\bigr) \le \mathrm{tol}\) (with the computed residual \(\widetilde{R}_m\) of Remark~4.3 replacing \(R_m\) whenever thermal solves are iterative). We solve each Newton step by preconditioned conjugate gradient (PCG):
\begin{itemize}
\item In 2D, \(H_{\Ec,\mathrm{loc}}\) is assembled by batched tensor contractions across \((\text{level}, p)\) element groups and \(\varrho S_m M_c + k H_{\Ec,\mathrm{loc}}\) is factored via sparse Cholesky as an optimal PCG preconditioner.
\item In 3D, to eliminate all element-Hessian memory overhead, both \(H_{\Ec,\mathrm{loc}}(v)\) and \(B_c(v)\) are applied \emph{matrix-free via 3D sum factorization} \cite{deville2002} (\(O(r_K^4)\) operations per element), while the constant-coefficient SPD thermal operators \eqref{eq:schur-thermal} and phase Helmholtz preconditioner \((\varrho S_m + k\bar{c}_r)M_c + k(1+\mu_A)S_c\) are assembled sparse matrices, with \(\bar{c}_r=\max(\langle c_{\mathrm{react}}\rangle_\nu,\varepsilon^{-2})\) and \(c_{\mathrm{react}}=3\varepsilon^{-2}\phi^2 + c_T K(M_\nu h''(\phi)+\eta_\nu G''(\phi))\ge 0\). These assembled SPD systems use direct sparse Cholesky for \(N_{\mathrm{conf}}\le 30{,}000\) and switch automatically to incomplete Cholesky (\texttt{ichol}) preconditioned CG (tolerance \(10^{-12}\)) for \(N_{\mathrm{conf}}>30{,}000\), recording thermal residual norms and solver flags for the a posteriori certificate below.
\end{itemize}

\subsection{Fully discrete well-posedness and original-energy dissipation}

\begin{theorem}[Fully discrete existence and uniqueness]\label{thm:disc-exist}
For every \(u^n\in X\) and every \(k>0\), the fully discrete predictor \eqref{eq:disc-pred} and each of the five main stages \eqref{eq:disc-stage} admit a unique solution in \(X\).
\end{theorem}
\begin{proof}
On the finite-dimensional space \(X=V\times V_0\), the discrete energy \(\Ec\) is continuous and globally convex, hence admits an affine lower bound at any fixed point. Adding the strictly convex quadratic metric term \(\frac{S_m}{2k}(J_*(U-\widehat{U}_m), U-\widehat{U}_m)_X\) with \(S_m>0\) and \(J_*\) symmetric positive definite renders the stage objective coercive and strictly convex on \(X\), establishing existence and uniqueness of a global minimizer whose Euler--Lagrange equation is \eqref{eq:disc-stage}.
\end{proof}

\begin{theorem}[Original discrete energy dissipation and residual-work certificate]\label{thm:disc-energy}
Let \(U_0=u^n\in X\) and let \(U_1,\dots,U_5\in X\) be any computed stage iterates with residuals \(G_m=G_m(U_m)\) in \eqref{eq:disc-stage} and increments \(d_m=U_m-U_{m-1}\). Then the original discrete energy \eqref{eq:discrete-energy} satisfies
\begin{equation}
 \Eorig(U_5) - \Eorig(U_0)
 \le \sum_{m=1}^5 \left[\frac{(G_m, d_m)_X}{k} - \frac{S_m}{2k}(J_* d_m, d_m)_X\right].
 \label{eq:disc-work-ineq}
\end{equation}
In particular, exact stages (\(G_m=0\)) unconditionally dissipate the original discrete energy for every \(k>0\):
\begin{equation}
 \Eorig(u^{n+1}) + \frac{1}{2k}\sum_{m=1}^5 S_m(J_* d_m, d_m)_X \le \Eorig(u^n).
 \label{eq:disc-energy-ineq}
\end{equation}
Moreover, the exact physical integrated enthalpy is conserved:
\(\int_\Omega\bigl(T^{n+1}-K h(\phi^{n+1})\bigr)\,dx = |\Omega|\bar{H}\).
\end{theorem}
\begin{proof}
Define the discrete stage majorizers \(F_m(V)\) identically to \eqref{eq:Fm-majorizer} with \(\Ec,\Ee\) and the discrete inner product \((\cdot,\cdot)_X\). By the backward recursion \eqref{eq:backward-recursion}, \(\nabla_X F_m(U_m)=G_m/k\) and \(F_m\) is \((S_m/k)\)-strongly convex in the \(J_*\)-norm. Strong convexity gives
\begin{align*}
 F_m(U_m) - F_m(U_{m-1}) &\le \bigl(\nabla_X F_m(U_m),\, U_m-U_{m-1}\bigr)_X - \frac{S_m}{2k}(J_* d_m, d_m)_X \\
 &= \frac{(G_m, d_m)_X}{k} - \frac{S_m}{2k}(J_* d_m, d_m)_X.
\end{align*}
The inter-stage inequality \(F_m(U_{m-1})\le F_{m-1}(U_{m-1})\) and endpoint identities \(F_1(U_0)=\Eorig(U_0)\), \(F_5(U_5)\ge \Eorig(U_5)\) follow verbatim from the proof of Theorem~\ref{thm:semi-energy}. Telescoping over \(m=1,\dots,5\) proves \eqref{eq:disc-work-ineq}, and setting \(G_m=0\) yields \eqref{eq:disc-energy-ineq}. Finally, because \(\nu_K\ge 3r_K+1\) integrates both \(T^{n+1}\) and the quintic polynomial \(h(\phi^{n+1})\) exactly, \(\int_\Omega(T^{n+1}-Kh(\phi^{n+1}))\,dx = (T^{n+1},1)_M - K(h(\phi^{n+1}),1)_\nu = |\Omega|\bar{H}\).
\end{proof}

\begin{remark}[Two-component a posteriori certification and step acceptance]
For any computed stage pair \((a_m,b_m)\), define the computed phase residual \(\widetilde{R}_m \coloneqq \varrho S_m M_c(a_m-\widehat{a}_m) + k c_T K B_c(\phi_*) b_m + k(\nabla\Ec^{\mathrm{phase}}(a_m)-\ell_{e,m})\) evaluated at the actual computed thermal coordinate \(b_m\) (which coincides with the exact-elimination gradient \(R_m(a_m)\) in \eqref{eq:reduced-res} whenever \(b_m=b_m(a_m)\)), and denote by \(e_m^{\mathrm{th}}\) the zero-mean discrete thermal residual \(P_0^T\bigl[(S_m M_c + kDS_c)b_m - S_m M_c\widehat{b}_m - S_m K P_0^T B_c(\phi_*)(a_m-\widehat{a}_m)\bigr]\) of \eqref{eq:schur-thermal}, with mass-Riesz representer \(\tilde{e}_m^{\mathrm{th}}\in V_0\) (\(\|\tilde{e}_m^{\mathrm{th}}\|_M^2 = (e_m^{\mathrm{th}})^T M_c^{-1} e_m^{\mathrm{th}}\)). Setting \(\Delta a_m = a_m - a_{m-1}\) and \(z_m=(w_m-w_{m-1})-K B_*(\phi_m-\phi_{m-1})\in V_0\), the exact two-component residual work decomposes as \((G_m,d_m)_X = \widetilde{R}_m^T\Delta a_m + (\mathcal{K}^{-1}\tilde{e}_m^{\mathrm{th}}, z_m)_M\). Applying Young's inequality \((\mathcal{K}^{-1}\tilde{e}_m^{\mathrm{th}},z_m)_M \le \frac{S_m}{2}\|z_m\|_{\mathcal{K}^{-1}}^2 + \frac{1}{2S_m}\|\tilde{e}_m^{\mathrm{th}}\|_{\mathcal{K}^{-1}}^2\) absorbs the cross term into the thermal part of \((J_* d_m,d_m)_X\), and the conforming box Neumann Poincar\'e bound \(\|\nabla v\|_{L^2}^2\ge (\pi/L_{\mathrm{box}})^2\|v\|_{L^2}^2\) on \(V_0\) (with maximum box side length \(L_{\mathrm{box}}\)) yields the computable upper bound without inverting \(\mathcal{K}\):
\begin{equation}
 \Eorig(U_5) - \Eorig(U_0)
 \le \sum_{m=1}^5 \left[\frac{\widetilde{R}_m^T \Delta a_m}{k}
   - \frac{\varrho S_m}{2k}\|\Delta a_m\|_{M_c}^2
   + \frac{c_T L_{\mathrm{box}}^2}{2\pi^2 D k S_m}\,(e_m^{\mathrm{th}})^T M_c^{-1} e_m^{\mathrm{th}}\right].
 \label{eq:two-comp-cert}
\end{equation}
When direct Cholesky is used (\(N_{\mathrm{conf}}\le 30{,}000\)), \(e_m^{\mathrm{th}}\approx 0\) and \eqref{eq:two-comp-cert} reduces to the phase-only certificate; when 3D iterative \texttt{ichol}--PCG is active (\(N_{\mathrm{conf}}>30{,}000\)), the third term rigorously bounds the contribution of the iterative thermal residual to the energy certificate using the already factored \(M_c\). Whereas Theorems~\ref{thm:disc-energy} and~\ref{thm:global-energy} hold with zero tolerance in exact arithmetic, in IEEE double precision our unified steppers enforce an explicit acceptance gate requiring predictor, five-stage, and thermal-solver convergence together with the floating-point bounds \(\Eorig(U_5)-\Eorig(U_0)\le 10^{-13}\max(1,|\Eorig(U_0)|)\), right-hand side of \eqref{eq:two-comp-cert} \(\le 10^{-12}\max(1,|\Eorig(U_0)|)\), and \(|\mathcal{H}(U_5)-|\Omega|\bar{H}|\le 10^{-11}\max(1,|\Omega|\max(1,|\bar{H}|))\); any rejected full step is automatically retried from \(u^n\) as two certified half-steps (\(k/2\)) and committed only if both half-steps pass.
\end{remark}

\begin{remark}[Global energy stability without a discrete maximum principle]
Although the continuous PDE preserves \(\phi\in[-1,1]\) (Proposition~\ref{prop:invariant-interval}), high-order spectral-element and second-order multistage schemes do not in general satisfy a strict pointwise discrete maximum principle without nonlinear limiting. Crucially, unlike splitting analyses that truncate \(F\) or \(h\) outside \([-1,1]\), our sextic majorizer \(G(s)\) establishes global convexity of \(\Ec,\Ee\) on all of \(\R\) (Lemma~\ref{lem:convexity}), so Theorems~\ref{thm:disc-exist} and~\ref{thm:disc-energy} hold unconditionally without assuming \(|\phi|\le 1\).
\end{remark}

% ============================================================================
\section{Moving-interface \texorpdfstring{\(hp\)}{hp} refinement strategy and analysis}\label{sec:hp}
% ============================================================================
\subsection{\texorpdfstring{\(2{:}1\)}{2:1}-balanced quadtree/octree conforming trace injection}\label{sec:hp-trace}
To avoid recursive chains of hanging-node dependencies and ill-conditioned constraint matrices, we maintain \emph{full \(2{:}1\) balance} across all face, edge, and corner neighbors of the tree partition \(\mathcal{T}_h\): any two leaf cells sharing a face, edge, or vertex differ in refinement level by at most one. This structural invariant yields two decisive benefits:
\begin{enumerate}
\item \textbf{Smooth geometric grading for the thermal boundary layer.} Outside the finest interfacial layer at level \(L_{\max}\), \(2{:}1\) balance automatically constructs a smooth octave-by-octave spatial transition into the far-field melt, preventing abrupt element-size jumps that degrade thermal-gradient accuracy and linear-solver conditioning.
\item \textbf{Explicit \(O(N)\) coarse-to-fine topological trace injection \(J\).} Under full \(2{:}1\) balance, every hanging vertex lies either at the midpoint of a unique maximal coarse edge or (in 3D) at the center of a unique maximal coarse face whose boundary entities are strictly coarser. In 2D, every maximal edge \(E\) carries the minimum degree \(q_E=\min\{r_K : K\text{ touches }E\}\) of its incident cells; in 3D, a brief fixed-point consistency pass ensures that every maximal boundary edge \(E\subset\partial F\) of a maximal face \(F\) satisfies \(q_E\le q_F\) (and similarly for face-midline edges). Ordering maximal edges and faces from coarsest to finest level, every hanging vertex, child edge, and quadrant child face is expressed in a single non-recursive pass by 1D/2D barycentric Lagrange interpolation from its strictly coarser master trace. Cell interior GLL nodes remain completely free.
\end{enumerate}
The resulting sparse \(O(N)\)-nonzero injection \(J\in\R^{N_{\mathrm{broken}}\times N_{\mathrm{conf}}}\) requires no global SVD or QR factorization. Because every constrained fine-side trace node is uniquely expressed by exact barycentric Lagrange interpolation from its strictly coarser master trace (with nested edge-to-face degree compatibility \(q_E\le q_F\) in 3D), while all master trace DOFs and element-interior GLL DOFs remain linearly independent, \(J\) has full column rank and maps \(\R^{N_{\mathrm{conf}}}\) onto a conforming \(hp\) subspace with continuous inter-element traces (matching the minimum-degree trace convention). In code, optional diagnostics assert constant reproduction (\(\|J\mathbf{1}_c-\mathbf{1}\|_\infty<10^{-13}\)) and random-vector inter-element trace continuity (\(\max|\text{jump}|<10^{-12}\), observed at \(\sim 10^{-16}\) across all verified meshes); see Figure~\ref{fig:hp-schematic}(a).

\subsection{Moving-interface \texorpdfstring{\([-0.98,0.98]\)}{[-0.98, 0.98]} \texorpdfstring{\(hp\)}{hp} adaptivity strategy}
\begin{figure}[htbp]
\centering
\includegraphics[width=\linewidth]{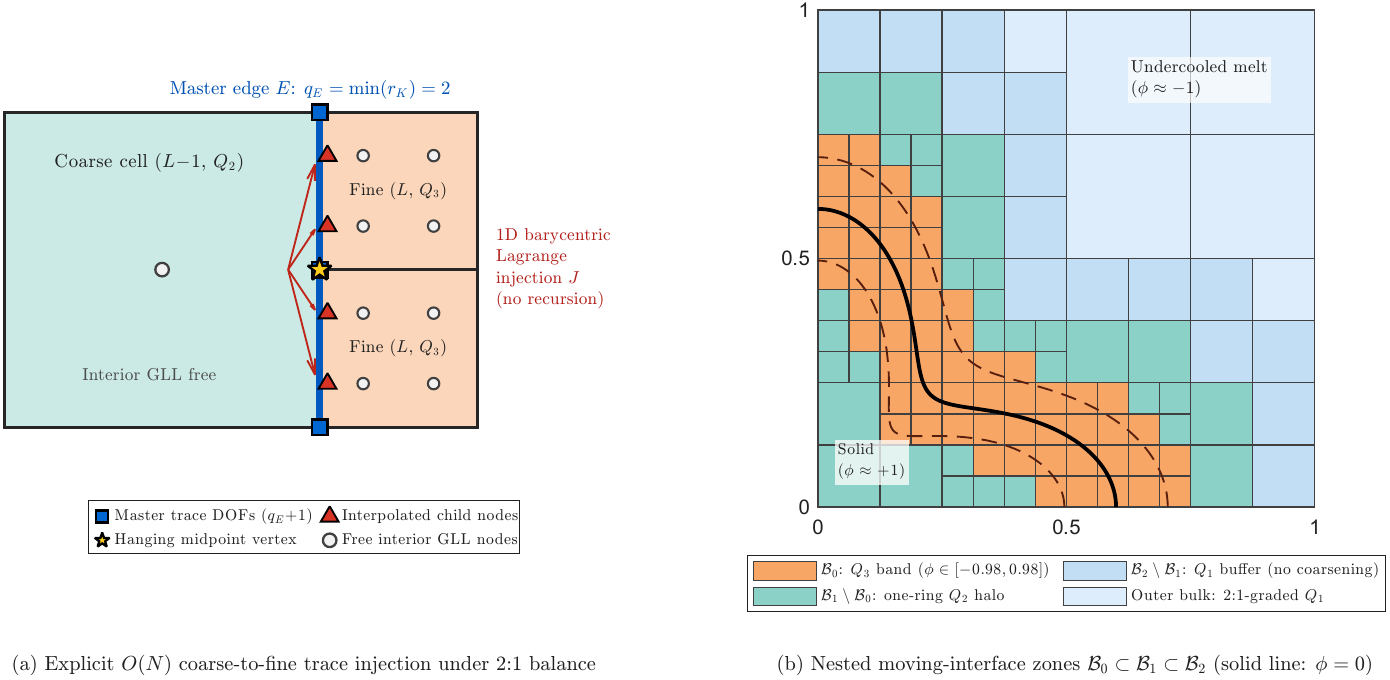}
\caption{Schematic of the \(2{:}1\)-balanced \(hp\) spectral-element discretization and moving-interface adaptivity. (a)~Coarse-to-fine conforming trace injection \(J\) across a \(2{:}1\) nonmatching interface: the hanging midpoint vertex (yellow star) and constrained fine-side GLL trace nodes (red triangles) are interpolated in a single non-recursive pass from the strictly coarser master trace DOFs (blue squares) at degree \(q_E=\min(r_{K^-},r_{K^+})\), while element interior GLL nodes (open circles) remain unconstrained. (b)~Nested adaptive zones \(\mathcal{B}_0\subset\mathcal{B}_1\subset\mathcal{B}_2\subset\mathcal{T}_h\) around a curved solid--liquid front: \(Q_3\) enrichment on the \([-0.98,0.98]\) interface band \(\mathcal{B}_0\), one-ring \(Q_2\) thermal halo \(\mathcal{B}_1\setminus\mathcal{B}_0\), two-ring coarsening buffer \(\mathcal{B}_2\setminus\mathcal{B}_1\), and \(2{:}1\)-graded \(Q_1\) bulk.}
\label{fig:hp-schematic}
\end{figure}

Because \(\phi\) varies rapidly only across the solid--liquid transition layer while remaining exponentially close to \(\pm 1\) in the solid core and undercooled melt, we drive spatial adaptation directly by the physical interface interval \([-0.98,0.98]\) (illustrated in Figure~\ref{fig:hp-schematic}(b)):
\begin{enumerate}
\item \textbf{Interface-band marking and \(h\)-refinement.} A cell \(K\in\mathcal{T}_h\) is marked as an \emph{interface cell} if its sampled polynomial phase values (evaluated at its nodal GLL/quadrature points) intersect \([-0.98,0.98]\), i.e., \(\max_{\bm{x}\in\mathcal{S}_K}\phi(\bm{x}) \ge -0.98\) and \(\min_{\bm{x}\in\mathcal{S}_K}\phi(\bm{x}) \le 0.98\). Every marked cell below the target interface level \(L_{\max}\) is recursively subdivided into \(2^d\) children, and the band test is re-evaluated at the sampled nodes of each child using the exact donor polynomial until all marked cells reach level \(L_{\max}\); denote the resulting set of interface cells by \(\mathcal{B}_0\subset\mathcal{T}_h\).
\item \textbf{Bulk \(h\)-coarsening.} Let \(\mathcal{B}_1\subset\mathcal{T}_h\) denote the union of \(\mathcal{B}_0\) with its one-ring neighbor layer, and let \(\mathcal{B}_2\subset\mathcal{T}_h\) be the two-ring protection zone obtained by extending \(\mathcal{B}_1\) by one safety ring to buffer the moving front between adaptation intervals (\(\mathcal{B}_0\subset\mathcal{B}_1\subset\mathcal{B}_2\)). Any complete family of \(2^d\) sibling cells above level \(L_{\min}\) lying strictly outside \(\mathcal{B}_2\) (all sampled values satisfying \(\phi>0.98\) or \(\phi<-0.98\)) is merged into its parent cell.
\item \textbf{\(2{:}1\) balance closure and zoned \(p\)-assignment.} Full \(2{:}1\) balance is restored on the updated tree. Polynomial degrees are then assigned by physical zone:
\begin{equation}
 r_K = \begin{cases}
  p_{\mathrm{high}} & \text{if } K\in\mathcal{B}_0\text{ (interface band }[-0.98,0.98]\text{)},\\
  p_{\mathrm{mid}}  & \text{if } K\in\mathcal{B}_1\setminus\mathcal{B}_0\text{ (one-ring thermal halo)},\\
  p_{\mathrm{low}}  & \text{otherwise (far-field pure solid/liquid bulk)},
 \end{cases}
 \label{eq:hp-zoning}
\end{equation}
with default hierarchy \((p_{\mathrm{low}},p_{\mathrm{mid}},p_{\mathrm{high}})=(1,2,3)\). Under the minimum-degree trace rule of Section~\ref{sec:hp-trace}, the intermediate \(Q_2\) halo insulates the \(Q_3\) band from the \(Q_1\) bulk: every interior edge or face strictly inside the \([-0.98,0.98]\) transition layer retains full \(Q_3\) conformity, while trace-degree transitions to \(Q_2\) and \(Q_1\) occur only in the exponentially flat tails (\(|\phi|>0.98\)) where the smoother thermal field governs spatial error.
\end{enumerate}

\subsection{Enthalpy-conserving, energy-nonincreasing inter-mesh transfer}
When the adaptive mesh changes from \(\mathcal{T}_{\mathrm{old}}\) (space \(V_{\mathrm{old}}\)) to \(\mathcal{T}_{\mathrm{new}}\) (space \(V_{\mathrm{new}}\)), both the approximation space and the quadrature partition change (note that because \(A(\nabla\phi)\) is non-polynomial, \(\Eorig^{\mathrm{old}}\) and \(\Eorig^{\mathrm{new}}\) denote the quadrature-defined original discrete energies on \(\mathcal{T}_{\mathrm{old}}\) and \(\mathcal{T}_{\mathrm{new}}\), respectively). To preserve exact enthalpy conservation and global monotonicity of this quadrature-defined discrete energy, we construct a trial state \((\widehat{\phi},T_{\mathrm{new}})\) in three steps followed by an explicit acceptance guard:
\begin{enumerate}
\item \textbf{Common-subdivision \(L^2\) mass projection.} Compute \(\widehat{\phi}\in V_{\mathrm{new}}\) and \(\widetilde{T}\in V_{\mathrm{new}}\) solving
\[
 (\widehat{\phi},v)_{M_{\mathrm{new}}} = (\phi_{\mathrm{old}},v)_{L^2(\Omega)}, \qquad
 (\widetilde{T},v)_{M_{\mathrm{new}}} = (T_{\mathrm{old}},v)_{L^2(\Omega)} \qquad \forall v\in V_{\mathrm{new}},
\]
evaluating the right-hand sides \emph{exactly} via Gauss quadrature on every nonempty intersection \(K_{\mathrm{new}}\cap K_{\mathrm{old}}\). Because \(\mathcal{T}_{\mathrm{old}}\) and \(\mathcal{T}_{\mathrm{new}}\) are both dyadic subtrees of the same root domain, every nonempty intersection coincides with the smaller of the two nested cells (\(K_{\mathrm{new}}\subseteq K_{\mathrm{old}}\) or \(K_{\mathrm{old}}\subset K_{\mathrm{new}}\)), avoiding polygon clipping and bounding the total number of intersections by \(O(N_{\mathrm{old}}+N_{\mathrm{new}})\); using dyadic leaf and ancestor-to-descendant hash maps, overlap lookup requires \(O(L_{\max}(N_{\mathrm{old}}+N_{\mathrm{new}}))\) operations, followed by sparse Cholesky mass solves on \(V_{\mathrm{new}}\).
\item \textbf{Exact enthalpy restoration.} Decompose \(\widetilde{T} = \widetilde{w} + \langle \widetilde{T}\rangle\mathbf{1}\) with \(\widetilde{w}\in V_{0,\mathrm{new}}\), and set the transferred mean temperature via the exact enthalpy lift:
\begin{equation}
 \bar{T}_{\mathrm{new}} = \bar{H} + K\langle h(\widehat{\phi})\rangle_{\nu_{\mathrm{new}}}.
 \label{eq:transfer-enthalpy}
\end{equation}
This guarantees \(\int_\Omega\bigl((\bar{T}_{\mathrm{new}}+\widetilde{w}) - K h(\widehat{\phi})\bigr)\,dx = |\Omega|\bar{H}\) identically.
\item \textbf{Energy-safe radial thermal projection and fallback acceptance guard.} Let \(\Eorig^{\mathrm{old}}\) be the donor energy and \(\Eorig^{\mathrm{raw}}=\Eorig^{\mathrm{new}}(\widehat{\phi},\bar{T}_{\mathrm{new}}+\widetilde{w})\). If \(\Eorig^{\mathrm{raw}} \le \Eorig^{\mathrm{old}}\), set \(w_{\mathrm{new}}=\widetilde{w}\). Otherwise, if \(\Eorig^{\mathrm{old}}\ge E_{\mathrm{base}}\) and \(\|\widetilde{w}\|_{M_{\mathrm{new}}}>0\), contract \(\widetilde{w}\) radially according to Proposition~\ref{prop:energy-transfer}. Finally, accept \((\mathcal{T}_{\mathrm{new}},\widehat{\phi},w_{\mathrm{new}})\) only if the trial state is finite and satisfies the discrete energy and enthalpy bounds; otherwise reject the candidate mesh and retain \((\mathcal{T}_{\mathrm{old}},\phi_{\mathrm{old}},w_{\mathrm{old}})\).
\end{enumerate}

\begin{proposition}[Energy-safe radial thermal projection]\label{prop:energy-transfer}
Fix \(\widehat{\phi}\in V_{\mathrm{new}}\) and \(\bar{T}_{\mathrm{new}}\) from \eqref{eq:transfer-enthalpy}, and define the baseline energy \(E_{\mathrm{base}} = \Eorig^{\mathrm{new}}(\widehat{\phi},\bar{T}_{\mathrm{new}}\mathbf{1})\). Whenever \(\Eorig^{\mathrm{old}}\ge E_{\mathrm{base}}\) and \(\|\widetilde{w}\|_{M_{\mathrm{new}}}>0\), the unique mass-norm closest zero-mean thermal fluctuation \(w_{\mathrm{new}}\in V_{0,\mathrm{new}}\) to \(\widetilde{w}\) satisfying \(\Eorig^{\mathrm{new}}(\widehat{\phi},\bar{T}_{\mathrm{new}}\mathbf{1}+w_{\mathrm{new}})\le \Eorig^{\mathrm{old}}\) is given by
\begin{equation}
 w_{\mathrm{new}} = \alpha_T \widetilde{w}, \qquad
 \alpha_T = \min\!\left(1,\, \sqrt{\frac{2(\Eorig^{\mathrm{old}} - E_{\mathrm{base}})}{c_T \|\widetilde{w}\|_{M_{\mathrm{new}}}^2}}\right) \in [0,1].
 \label{eq:radial-contraction}
\end{equation}
Moreover, the reconstructed field \(T_{\mathrm{new}}=\bar{T}_{\mathrm{new}}\mathbf{1}+w_{\mathrm{new}}\) conserves total physical enthalpy identically and satisfies \(\Eorig^{\mathrm{new}}(\widehat{\phi},T_{\mathrm{new}})\le \Eorig^{\mathrm{old}}\).
\end{proposition}
\begin{proof}
Because \(w\in V_{0,\mathrm{new}}\) is mass-orthogonal to constants,
\[
 \Eorig^{\mathrm{new}}(\widehat{\phi},\bar{T}_{\mathrm{new}}\mathbf{1}+w)
 = E_{\mathrm{base}} + \frac{c_T}{2}\|w\|_{M_{\mathrm{new}}}^2.
\]
Hence the energy constraint \(\Eorig^{\mathrm{new}}\le \Eorig^{\mathrm{old}}\) is the exact mass-norm ball \(\|w\|_{M_{\mathrm{new}}}\le R_{\mathrm{cap}} \coloneqq \sqrt{2(\Eorig^{\mathrm{old}}-E_{\mathrm{base}})/c_T}\). Projecting \(\widetilde{w}\) onto this centered Euclidean ball in the \(M_{\mathrm{new}}\)-inner product yields uniquely \eqref{eq:radial-contraction}. Because \(\langle\alpha_T\widetilde{w}\rangle=0\), the constant temperature mean \(\bar{T}_{\mathrm{new}}\) and thus the total enthalpy \(\mathcal{H}\) are unchanged.
\end{proof}

\begin{theorem}[Global quadrature-energy monotonicity across accepted moving \(hp\) meshes]\label{thm:global-energy}
Suppose each accepted time step satisfies the certified CSRK energy inequality \eqref{eq:disc-energy-ineq} (or \eqref{eq:two-comp-cert} with nonpositive right-hand side) and each moving-interface \(hp\) adaptation applies Proposition~\ref{prop:energy-transfer} with the fallback acceptance guard (retaining \(\mathcal{T}_{\mathrm{old}}\) if infeasible). Then the sequence of quadrature-defined original discrete energies \(\{\Eorig(u^n)\}_{n\ge 0}\) is globally nonincreasing across all accepted time steps and mesh adaptations, and the integrated physical enthalpy equals \(|\Omega|\bar{H}\) at every step.
\end{theorem}
\begin{proof}
At any step \(n\) without an accepted mesh change, \(\Eorig(u^{n+1})\le \Eorig(u^n)\) by Theorem~\ref{thm:disc-energy}. At an accepted adaptation event following step \(n\), Proposition~\ref{prop:energy-transfer} ensures \(\Eorig^{\mathrm{new}}(u_{\mathrm{tr}}^n)\le \Eorig^{\mathrm{old}}(u^n)\) (whereas a rejected candidate retains \((\mathcal{T}_{\mathrm{old}},u^n)\) with zero change), and Theorem~\ref{thm:disc-energy} on the active space gives \(\Eorig^{\mathrm{new}}(u^{n+1})\le \Eorig^{\mathrm{new}}(u_{\mathrm{tr}}^n)\). Chaining these inequalities establishes global monotonicity and exact enthalpy conservation.
\end{proof}

\begin{remark}[Feasibility margin and quantitative benignness of radial contraction]\label{rem:alpha-T}
Proposition~\ref{prop:energy-transfer} requires \(\Eorig^{\mathrm{old}}\ge E_{\mathrm{base}}\) and \(\|\widetilde{w}\|_{M_{\mathrm{new}}}>0\) when \(\Eorig^{\mathrm{raw}}>\Eorig^{\mathrm{old}}\); whenever either condition fails, the fallback acceptance guard rejects the candidate mesh and preserves \(\Eorig^{\mathrm{old}}\) trivially. Because our \(hp\) strategy avoids coarsening inside the marked interface zone \([-0.98,0.98]\) and coarsens only in the bulk where \(|\nabla\phi|\approx 0\) and \(F(\phi)\approx 0\), the phase-energy transfer perturbation is negligible relative to the thermal fluctuation energy \(\frac{c_T}{2}\|\widetilde{w}\|_{M_{\mathrm{new}}}^2\sim O(10^2\text{--}10^4)\), so \(\Eorig^{\mathrm{old}}\gg E_{\mathrm{base}}\) holds with a wide safety margin. Empirically, across all \(1{,}453\) accepted mesh adaptations recorded in the seven main trajectories of Examples~2 and~3 (excluding the auxiliary ablation studies, where the stripped narrow-band \([-0.90,0.90]\) 3D case records \(1-\alpha_T=1.52\times 10^{-3}\)), every archived time step strictly decreases \(\Eorig\) (\(\max\Delta\Eorig^{\mathrm{step}}<0\)), the radial contraction satisfies \(1-\alpha_T \le 2.40\times 10^{-8}\) in 2D and \(1-\alpha_T \le 3.08\times 10^{-6}\) in 3D---confirming that Proposition~\ref{prop:energy-transfer} acts purely as an \(O(10^{-8}\text{--}10^{-6})\) safeguard rather than artificial thermal damping---and all archived energy transfers lie within the floating-point rounding tolerance \(\Delta\Eorig^{\mathrm{tr}}\le 10^{-13}\max(1,|\Eorig^{\mathrm{old}}|)\), with peak positive IEEE double-precision increments of \(9.09\times 10^{-13}\) (2D \(m=4\)), \(7.28\times 10^{-12}\) (2D \(m=6\)), and \(2.91\times 10^{-11}\) (3D, \(\le 2.3\times 10^{-16}\) relative, \(1\)--\(2\) ulps). In the release codebase, the unified steppers and adaptation drivers additionally record explicit per-step retry flags, candidate-transfer rejections, iterative thermal residuals, and transfer enthalpy defects for new runs.
\end{remark}

\begin{algorithm}[htbp]
\caption{Adaptive \(2{:}1\)-balanced \(hp\) spectral-element CSRK time-stepping loop}
\label{alg:adaptive-csrk}
\begin{algorithmic}[1]
\State Construct initial \(2{:}1\)-balanced tree \(\mathcal{T}_0\) refined on sampled \(\phi_0\in[-0.98,0.98]\), assemble \(J,M_c,S_c\), project \((\phi^0,T^0)\), and record \(\bar{H}=\langle T^0-Kh(\phi^0)\rangle\).
\For{\(n = 0, 1, \dots, N_t-1\)}
  \If{mesh or step size \(k\) changed, or \(n=0\)}
    \State Factor/setup and cache SPD thermal Helmholtz solvers \(M_c + \frac{kD}{2}S_c\) and \(S_m M_c + kD S_c\) (\(m=1,\dots,5\)).
  \EndIf
  \State Solve predictor \eqref{eq:disc-pred} with frozen \(B_c(\phi^n)\) via Schur-eliminated damped Newton--PCG to obtain \(U_*=(\phi_*,w_*)\).
  \State Freeze midpoint coupling operator \(B_c(\phi_*)\) and solve stages \(m=1,\dots,5\) of \eqref{eq:disc-stage} to obtain \(u^{n+1}=U_5\).
  \State Check convergence, two-component certificate \eqref{eq:two-comp-cert}, \(\Eorig(u^{n+1})\le\Eorig(u^n)\), and enthalpy defect \(|\mathcal{H}(u^{n+1})-|\Omega|\bar{H}|\); if rejected, retry with two half-steps (\(k/2\)).
  \If{adaptation check triggered at step \(n+1\)}
    \State Mark interface cells \(\mathcal{B}_0\) where sampled \(\phi^{n+1}\) intersects \([-0.98,0.98]\), refine to \(L_{\max}\), coarsen bulk siblings outside \(\mathcal{B}_2\), enforce \(2{:}1\) balance, and assign \((p_{\mathrm{low}},p_{\mathrm{mid}},p_{\mathrm{high}})\) by \eqref{eq:hp-zoning}.
    \State Transfer \((\phi^{n+1},T^{n+1})\) via \(O(L_{\max}(N_{\mathrm{old}}+N_{\mathrm{new}}))\) common-subdivision \(L^2\) projection, enthalpy shift \eqref{eq:transfer-enthalpy}, and Proposition~\ref{prop:energy-transfer} radial projection if \(\Eorig^{\mathrm{raw}}>\Eorig^{\mathrm{old}}\); accept only if feasible and finite energy/enthalpy bounds hold, else retain \(\mathcal{T}_{\mathrm{old}}\).
  \EndIf
\EndFor
\end{algorithmic}
\end{algorithm}

% ============================================================================
\section{Numerical experiments}\label{sec:numerics}
% ============================================================================
We validate the proposed framework through three progressively demanding benchmark suites under insulated Neumann boundaries: temporal and spatial \(h\)/\(p\)-convergence tests together with a direct CSRK-versus-SAV--BDF2 benchmark (Example~1), strong-anisotropy planar fourfold, sixfold, and full-domain asymmetric rotated/multi-nuclei dendritic growth with a 2D independent-reference component ablation (Example~2), and three-dimensional cubic dendritic evolution with a 3D independent-reference ablation (Example~3). All time-dependent CSRK benchmarks implement the fully discrete five-stage scheme \eqref{eq:disc-pred}--\eqref{eq:disc-stage} with common Gauss quadrature \(\nu_K = 3r_K+1\), Schur-complement thermal elimination \eqref{eq:schur-thermal}, and the \(2{:}1\)-balanced moving-interface \([-0.98,0.98]\) \(hp\) strategy of Algorithm~\ref{alg:adaptive-csrk}; the stationary manufactured-solution spatial test in Example~1 uses a dedicated stationary coupled solver with four extra quadrature points (\(\nu_K = 3r_K+5\)), and the SAV--BDF2 baseline uses a separate linearized integrator. Nonmatching conforming spaces are verified via constant reproduction (\(\|J\mathbf{1}_c-\mathbf{1}\|_\infty<10^{-13}\)) and random-vector inter-element trace-jump assertions (threshold \(10^{-12}\), observed \(\sim 10^{-16}\)). Computations were executed in IEEE double precision in MATLAB R2025b on Intel Core i9-12900K workstations (128~GB RAM). All MATLAB source codes, experiment-specific benchmark drivers, figure generators, and archived summary datasets are publicly available at \url{https://github.com/minghuili123456789/anisotropic-dendrite-csrk-hp}.

\subsection{Example 1: Two-dimensional spatial and temporal accuracy verification}\label{sec:exp1}
\paragraph{Temporal accuracy on a nonmatching mixed-degree \(hp\) quadtree.}
On \(\Omega=[-1,1]^2\), we construct a \(2{:}1\)-balanced nonmatching \(hp\) quadtree consisting of seven cells: three level-1 cells of degrees \(r_K\in\{3,4,3\}\) and four level-2 children of degrees \(r_K\in\{2,3,3,2\}\) occupying the upper-right quadrant, giving \(N_{\mathrm{broken}}=107\) and \(N_{\mathrm{conf}}=69\) conforming unknowns per field. The trace injection satisfies \(\|J\mathbf{1}_c-\mathbf{1}\|_\infty = 0\) and \(\max|\text{jump}| = 2.22\times 10^{-16}\) on a random conforming vector. Using physical parameters \((\varrho,\varepsilon,\sigma,\lambda,D,K)=(10^3,0.015,0.02,400,0.0025,0.9)\) for both \(m=4\) and \(m=6\) (\(\theta_0=0\)) and smooth initial data
\[
 \phi_0(x,y) = 0.12 - 0.38x - 0.27y + 0.06xy, \qquad
 T_0(x,y) = -0.08 + 0.03x - 0.02y,
\]
we integrate the unforced coupled system on this fixed discrete space to \(T_f=0.008\) with step sizes \(k\in\{0.004,\allowbreak 0.002,\allowbreak 0.001,\allowbreak 0.0005,\allowbreak 0.00025\}\) and measure fixed-space \(L^2(\Omega)\) errors in \(\phi\), \(T\), and the nonlinear enthalpy field \(H=T-Kh(\phi)\) against a fine-step reference trajectory at \(k_{\mathrm{ref}}=3.125\times 10^{-5}\) (\(256\) steps, \(8\times\) finer than the smallest tested step). (Because \((H-H_{\mathrm{ref}})^2\) has polynomial degree up to \(10r_K\), we verified that evaluating \(\|H-H_{\mathrm{ref}}\|_{L^2}\) with exact \(5r_K+1\) Gauss points agrees with the default \(3r_K+1\) quadrature norm to 10 significant digits, relative difference \(\le 1.1\times 10^{-10}\).)

\begin{figure}[!htb]
\centering
\includegraphics[width=\linewidth]{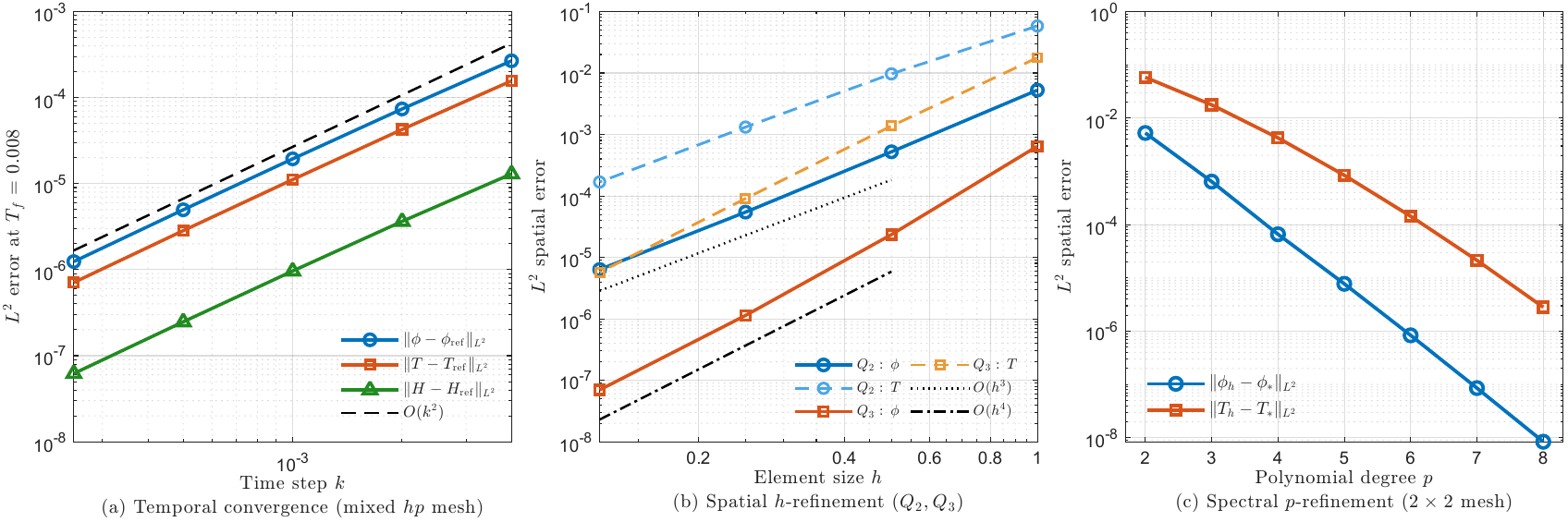}
\caption{Example~1: Two-dimensional accuracy verification. (a)~Second-order temporal convergence of the five-stage CSRK scheme on the nonmatching mixed-degree (\(Q_2/Q_3/Q_4\)) quadtree (\(m=4\)). (b)~Optimal spatial \(h\)-refinement rates \(O(h^3)\) for \(Q_2\) and \(O(h^4)\) for \(Q_3\). (c)~Exponential spectral \(p\)-refinement convergence on a fixed \(2\times 2\) mesh.}
\label{fig:exp1}
\end{figure}

Table~\ref{tab:exp1-time} and Figure~\ref{fig:exp1}(a) confirm clean second-order temporal convergence (asymptotic rates \(1.999\) in \(\phi\), \(2.001\) in \(T\), and \(1.996\) in \(H\)) for both fourfold and sixfold anisotropy. (The \(m=4\) and \(m=6\) errors agree to five significant digits because the stiff bulk reaction \(\varepsilon^{-2}F'(\phi)\sim O(4{,}400)\) dominates the modest \(\sigma=0.02\) anisotropy over \(T_f=0.008\), with anisotropy entering at the sixth digit in \(\phi,T\) and fifth digit in \(H\).) Across all trajectories, every step strictly decreases the quadrature-defined original discrete energy (\(\max_n\Delta\Eorig < 0\); e.g., at \(k=0.00025\), \(\Eorig\) decays monotonically from \(4{,}167.8883\) to \(4{,}001.7288\)), and the integrated enthalpy defect never exceeds \(1.11\times 10^{-16}\).

\begin{table}[htbp]
\centering
\caption{Example~1: Temporal \(L^2\) convergence, original-energy dissipation, and enthalpy conservation at \(T_f=0.008\) on the nonmatching mixed-degree (\(Q_2/Q_3/Q_4\)) quadtree.}
\label{tab:exp1-time}
\small
\begin{tabular}{@{}cccccccc@{}}
\toprule
\(k\) & \(\|\phi-\phi_{\mathrm{ref}}\|_{L^2}\) & Order & \(\|T-T_{\mathrm{ref}}\|_{L^2}\) & Order & \(\|H-H_{\mathrm{ref}}\|_{L^2}\) & Order & \(\max|\Delta\mathcal{H}|\) \\
\midrule
\multicolumn{8}{l}{\emph{Fourfold anisotropy (\(m=4,\ \sigma=0.02\)):}} \\
\(4.0\times 10^{-3}\) & \(2.6754\times 10^{-4}\) & --    & \(1.5649\times 10^{-4}\) & --    & \(1.2974\times 10^{-5}\) & --    & \(0.00\) \\
\(2.0\times 10^{-3}\) & \(7.3765\times 10^{-5}\) & 1.859 & \(4.2634\times 10^{-5}\) & 1.876 & \(3.6403\times 10^{-6}\) & 1.833 & \(0.00\) \\
\(1.0\times 10^{-3}\) & \(1.9397\times 10^{-5}\) & 1.927 & \(1.1140\times 10^{-5}\) & 1.936 & \(9.6440\times 10^{-7}\) & 1.916 & \(1.11\times 10^{-16}\) \\
\(5.0\times 10^{-4}\) & \(4.9633\times 10^{-6}\) & 1.966 & \(2.8411\times 10^{-6}\) & 1.971 & \(2.4759\times 10^{-7}\) & 1.962 & \(0.00\) \\
\(2.5\times 10^{-4}\) & \(1.2417\times 10^{-6}\) & 1.999 & \(7.0970\times 10^{-7}\) & 2.001 & \(6.2077\times 10^{-8}\) & 1.996 & \(1.11\times 10^{-16}\) \\
\midrule
\multicolumn{8}{l}{\emph{Sixfold anisotropy (\(m=6,\ \sigma=0.02\)):}} \\
\(4.0\times 10^{-3}\) & \(2.6754\times 10^{-4}\) & --    & \(1.5649\times 10^{-4}\) & --    & \(1.2974\times 10^{-5}\) & --    & \(0.00\) \\
\(2.0\times 10^{-3}\) & \(7.3765\times 10^{-5}\) & 1.859 & \(4.2634\times 10^{-5}\) & 1.876 & \(3.6404\times 10^{-6}\) & 1.833 & \(1.11\times 10^{-16}\) \\
\(1.0\times 10^{-3}\) & \(1.9397\times 10^{-5}\) & 1.927 & \(1.1140\times 10^{-5}\) & 1.936 & \(9.6442\times 10^{-7}\) & 1.916 & \(1.11\times 10^{-16}\) \\
\(5.0\times 10^{-4}\) & \(4.9633\times 10^{-6}\) & 1.966 & \(2.8411\times 10^{-6}\) & 1.971 & \(2.4760\times 10^{-7}\) & 1.962 & \(1.11\times 10^{-16}\) \\
\(2.5\times 10^{-4}\) & \(1.2417\times 10^{-6}\) & 1.999 & \(7.0970\times 10^{-7}\) & 2.001 & \(6.2078\times 10^{-8}\) & 1.996 & \(1.11\times 10^{-16}\) \\
\bottomrule
\end{tabular}
\end{table}

\paragraph{Spatial \(h\)-refinement and spectral \(p\)-refinement accuracy.}
To test spatial convergence independently of time stepping, we solve the stationary coupled system \(-\nabla\cdot\nabla_g A(\nabla\phi)+\varepsilon^{-2}F'(\phi)+c_T K h'(\phi)T=f_\phi\) and \(-D\Delta T=f_T\) on \(\Omega=[-1,1]^2\) with \((\varepsilon,\sigma,m,\theta_0,\lambda,D,K)=(0.25,0.03,4,0,2.0,0.5,0.8)\) and the non-polynomial smooth manufactured Neumann solution
\begin{align*}
 \phi_*(x,y) &= 0.72 + 0.08\sin\!\left(\frac{\pi x}{2}\right) + 0.06\sin\!\left(\frac{\pi y}{2}\right), \\
 T_*(x,y) &= 0.12\sin\!\left(\frac{\pi x}{2}\right)\sin\!\left(\frac{3\pi y}{2}\right) + 0.08\sin\!\left(\frac{3\pi x}{2}\right)\sin\!\left(\frac{\pi y}{2}\right),
\end{align*}
which satisfies \(\partial_{\bm{n}}\phi_*=\partial_{\bm{n}}T_*=0\) on \(\partial\Omega\), \(\langle T_*\rangle=0\), \(|\nabla\phi_*|>0\) throughout the interior \((-1,1)^2\) (vanishing only at the four boundary corners), and \(F''(\phi_*)>0\). Table~\ref{tab:exp1-space} and Figure~\ref{fig:exp1}(b)--(c) demonstrate optimal algebraic \(h\)-convergence rates \(O(h^{p+1})\) asymptotically on the finest levels for both \(Q_2\) (asymptotic orders \(3.101\) in \(\phi\) and \(2.969\) in \(T\), rising from \(2.589\) at the coarsest \(Q_2\) refinement) and \(Q_3\) (asymptotic orders \(4.026\) in \(\phi\) and \(3.980\) in \(T\)), together with rapid, monotonic spectral \(p\)-convergence on a fixed \(2\times 2\) mesh as \(p\) increases from \(2\) to \(8\) (\(\|\phi_h-\phi_*\|_{L^2}\) dropping nearly six orders of magnitude from \(5.26\times 10^{-3}\) to \(8.43\times 10^{-9}\) and \(\|T_h-T_*\|_{L^2}\) dropping over four orders of magnitude from \(5.82\times 10^{-2}\) to \(2.83\times 10^{-6}\)). Under matched degrees of freedom (\(N_{\mathrm{conf}}=289\)), spectral \(p\)-enrichment (\(p=8\)) outperforms \(h\)-refined \(Q_2\) by \(6{,}457\times\) in \(\phi\) and \(467\times\) in \(T\), and even surpasses \(h\)-refined \(Q_3\) at \(N_{\mathrm{conf}}=2{,}401\) (\(8.3\times\) more DOFs) by \(8.3\times\) in \(\phi\) and \(2.0\times\) in \(T\).

\begin{table}[!htb]
\centering
\caption{Example~1: Spatial \(h\)-refinement (\(Q_2, Q_3\)) and spectral \(p\)-refinement (\(2\times 2\) mesh, \(p=2,\dots,8\)) \(L^2\) errors for the smooth manufactured Neumann problem.}
\label{tab:exp1-space}
\small
\resizebox{\linewidth}{!}{%
\begin{tabular}{@{}cccccc@{\hspace{16pt}}cccc@{}}
\toprule
\multicolumn{6}{c}{\(h\)-refinement (\(Q_2\) and \(Q_3\))} & \multicolumn{4}{c}{\(p\)-refinement (\(2\times 2\) mesh)} \\
\cmidrule(r{14pt}){1-6}\cmidrule(l){7-10}
Space & \(h\) & \(N_{\mathrm{conf}}\) & \(\|\phi_h-\phi_*\|_{L^2}\) & Order & \(\|T_h-T_*\|_{L^2}\) (Ord.) & \(p\) & \(N_{\mathrm{conf}}\) & \(\|\phi_h-\phi_*\|_{L^2}\) & \(\|T_h-T_*\|_{L^2}\) \\
\midrule
\(Q_2\) & \(1.000\) & 25   & \(5.262\times 10^{-3}\) & --    & \(5.821\times 10^{-2}\) (--)    & 2 & 25  & \(5.262\times 10^{-3}\) & \(5.821\times 10^{-2}\) \\
\(Q_2\) & \(0.500\) & 81   & \(5.226\times 10^{-4}\) & 3.332 & \(9.672\times 10^{-3}\) (2.589) & 3 & 49  & \(6.431\times 10^{-4}\) & \(1.774\times 10^{-2}\) \\
\(Q_2\) & \(0.250\) & 289  & \(5.444\times 10^{-5}\) & 3.263 & \(1.319\times 10^{-3}\) (2.875) & 4 & 81  & \(6.608\times 10^{-5}\) & \(4.253\times 10^{-3}\) \\
\(Q_2\) & \(0.125\) & \(1{,}089\) & \(6.343\times 10^{-6}\) & 3.101 & \(1.684\times 10^{-4}\) (2.969) & 5 & 121 & \(7.759\times 10^{-6}\) & \(8.468\times 10^{-4}\) \\
\cmidrule(r{14pt}){1-6}
\(Q_3\) & \(1.000\) & 49   & \(6.431\times 10^{-4}\) & --    & \(1.774\times 10^{-2}\) (--)    & 6 & 169 & \(8.383\times 10^{-7}\) & \(1.441\times 10^{-4}\) \\
\(Q_3\) & \(0.500\) & 169  & \(2.349\times 10^{-5}\) & 4.775 & \(1.374\times 10^{-3}\) (3.691) & 7 & 225 & \(8.569\times 10^{-8}\) & \(2.142\times 10^{-5}\) \\
\(Q_3\) & \(0.250\) & 625  & \(1.136\times 10^{-6}\) & 4.371 & \(9.059\times 10^{-5}\) (3.922) & 8 & 289 & \(8.432\times 10^{-9}\) & \(2.825\times 10^{-6}\) \\
\(Q_3\) & \(0.125\) & \(2{,}401\) & \(6.970\times 10^{-8}\) & 4.026 & \(5.739\times 10^{-6}\) (3.980) &   &     & & \\
\bottomrule
\end{tabular}%
}
\end{table}

\paragraph{Direct same-mesh comparison with stabilized SAV--BDF2.}
To clarify the practical trade-offs between our original-energy-dissipative five-stage CSRK integrator and a representative linearized modified-energy auxiliary-variable implementation \cite{li2022,shen2018}, we conduct a controlled same-mesh benchmark on a fixed \(2{:}1\)-balanced \(hp\) quadtree covering the early growth region of the strong-anisotropy fourfold dendrite (\(m=4,\ \sigma=0.10,\ \mu_A=0.55,\ K=1.0\), \(97\) elements, \(N_{\mathrm{conf}}=637\) per field, \(\max|\text{jump}|=4.44\times 10^{-16}\), \(\Eorig(0)=4{,}793.34\)) over \(t\in[0,1.6]\) across seven step sizes \(k\in\{0.005,0.01,0.025,0.05,0.10,0.20,0.40\}\). We compare our CSRK scheme (free of empirical splitting/stabilization parameters, aside from physical, mesh, and solver tolerances) against a linearized second-order SAV--BDF2 discretization (quadratic stabilizer \(\frac{S_1}{2}\|\nabla\phi\|_{L^2}^2+\frac{S_2}{2\varepsilon^2}\|\phi\|_{L^2}^2\), shift \(C_0=25{,}000\), extrapolated latent-heat coupling \(h'(\phi^*)\), and one-step SAV--BDF1 startup) under three stabilization settings---unstabilized (\(S_1=S_2=0\)), moderate (\(S_1=S_2=2\)), and strong (\(S_1=S_2=4\))---evaluated against a fine CSRK temporal reference at \(k_{\mathrm{ref}}=0.001\) (\(\Eorig^{\mathrm{ref}}=4{,}671.37\), \(r_{\mathrm{tip}}^{\mathrm{ref}}=0.17255\)).

\begin{figure}[!htb]
\centering
\includegraphics[width=\linewidth]{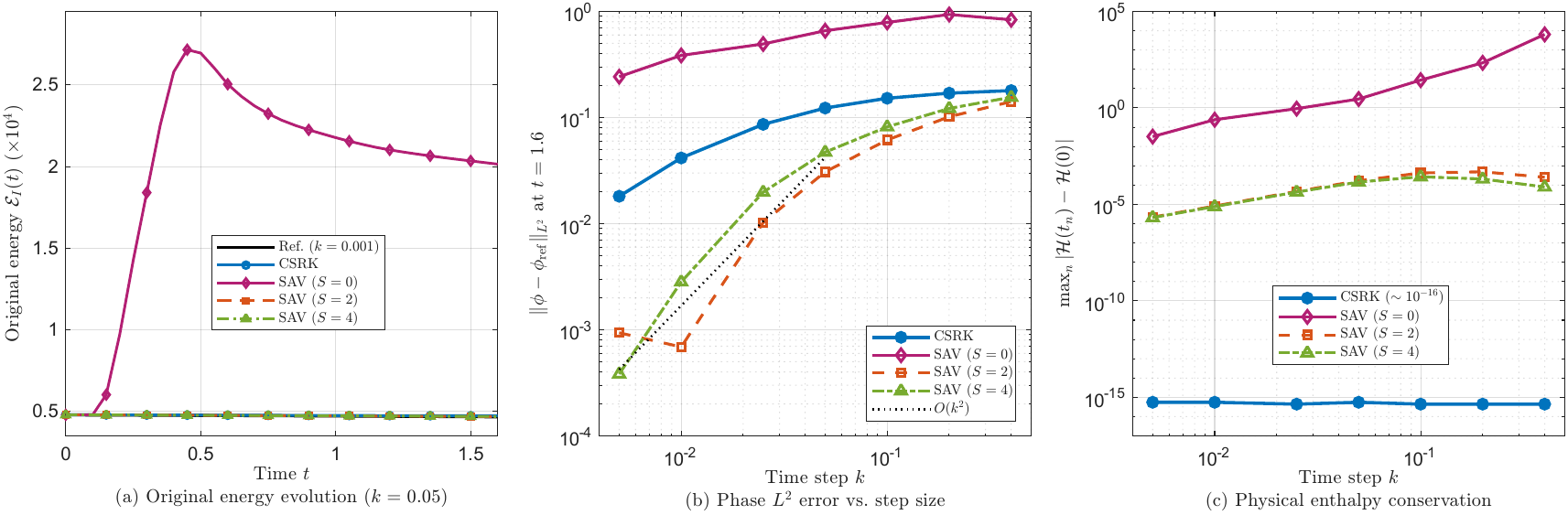}
\caption{Example~1: Direct same-mesh temporal comparison between the proposed five-stage CSRK scheme and the implemented linearized SAV--BDF2 scheme (\(S\in\{0,2,4\}\)) on the fixed \(2{:}1\)-balanced \(hp\) quadtree (\(t\in[0,1.6]\)). (a)~Original discrete energy \(\Eorig(t)\) at \(k=0.05\), exposing catastrophic original-energy growth in unstabilized SAV--BDF2 (\(S=0\)). (b)~Phase \(L^2\) error versus step size \(k\) relative to the \(k_{\mathrm{ref}}=0.001\) CSRK trajectory. (c)~Physical enthalpy conservation defect \(\max_n|\mathcal{H}(t_n)-\mathcal{H}(0)|\), illustrating \(10^{10}\times\)--\(10^{12}\times\) tighter conservation for CSRK.}
\label{fig:csrk-vs-sav}
\end{figure}

Table~\ref{tab:csrk-vs-sav} and Figure~\ref{fig:csrk-vs-sav} highlight three clear findings for this implementation and parameter set:
\begin{enumerate}
\item \textbf{Original-energy monotonicity vs.\ modified-energy drift.} Although SAV--BDF2 schemes are designed around a two-level \emph{modified} scalar energy, in our unstabilized SAV--BDF2 (\(S=0\)) run the \emph{original} physical energy \(\Eorig\) increases on every tested step size---from \(\max\Delta\Eorig^{\mathrm{step}}=+21.19\) at \(k=0.005\) to \(+4{,}536\) at \(k=0.05\) (driving \(\Eorig\) from \(4{,}793.3\) up to \(20{,}144.6\)) and \(+8.27\times 10^5\) at \(k=0.40\)---producing severely corrupted interfaces (\(\|\phi-\phi_{\mathrm{ref}}\|_{L^2}=0.243\)--\(0.936\)). By contrast, CSRK guarantees nonincreasing original discrete energy (observed \(\max\Delta\Eorig^{\mathrm{step}}<0\), zero energy increases) across the entire step-size range without empirical stabilization tuning.
\item \textbf{Machine-precision enthalpy conservation.} Because CSRK evolves the state directly in the enthalpy-conserving coordinate \eqref{eq:enthalpy-lift}, its physical enthalpy defect remains at \(\le 5.55\times 10^{-16}\) for all \(k\). Conversely, this linearized SAV--BDF2 implementation approximates the nonlinear latent-heat release as \(h'(\phi^*)\,\partial_t\phi\) using an extrapolated phase \(\phi^*\) rather than conserving \(T-Kh(\phi)\) identically, incurring \(10^{10}\times\)--\(10^{12}\times\) larger enthalpy drift even when stabilized (\(2.18\times 10^{-6}\) to \(4.89\times 10^{-4}\) for \(S=2,4\), and up to \(6.44\times 10^3\) for \(S=0\)).
\item \textbf{Honest cost and splitting-constant trade-off.} Because each SAV--BDF2 step requires only linear solves via the Sherman--Morrison scalar update, it is \(8\times\)--\(10\times\) faster per step than the five-stage nonlinear CSRK scheme on this small fixed mesh (\(0.48\)\,s vs.\ \(3.92\)\,s at \(k=0.05\)), and under tuned stabilization (\(S=2\)) exhibits a smaller temporal splitting prefactor at moderate \(k\) against the \(k_{\mathrm{ref}}=0.001\) reference (\(\|\phi-\phi_{\mathrm{ref}}\|_{L^2}=1.03\times 10^{-2}\) vs.\ \(8.63\times 10^{-2}\) at \(k=0.025\)). However, this SAV implementation depends sensitively on the empirical stabilizer \(S\) (failing badly at \(S=0\) and nearly doubling error from \(S=2\) to \(S=4\) at \(k=0.025\)), lacks a theoretical guarantee for \(\Eorig\), breaks exact enthalpy conservation, and its standard fixed-step two-level BDF2 modified-energy identity does not transfer directly across dynamic mesh changes without additional multistep restart or variable-step machinery. By contrast, our one-step CSRK integrator trades higher per-step nonlinear work (offset by \(7.4\times\)--\(75.3\times\) spatial \(hp\) DOF compression) and a larger global splitting constant for stabilization-free, certified original discrete energy monotonicity, machine-precision enthalpy conservation, and seamless one-step coupling with accepted moving \(hp\) meshes and variable step sizes (Theorem~\ref{thm:global-energy}).
\end{enumerate}

\begin{table}[!htb]
\centering
\caption{Example~1: Direct same-mesh temporal comparison between the proposed five-stage CSRK scheme and SAV--BDF2 (\(S\in\{0,2,4\}\)) on the fixed \(2{:}1\)-balanced \(hp\) quadtree (\(N_{\mathrm{conf}}=637,\ t=1.6,\ \Eorig(0)=4{,}793.34\)).}
\label{tab:csrk-vs-sav}
\footnotesize
\setlength{\tabcolsep}{4.5pt}
\resizebox{\linewidth}{!}{%
\begin{tabular}{@{}lccccccc@{}}
\toprule
Scheme & \(k\) & Wall (s) & \(\|\phi-\phi_{\mathrm{ref}}\|_{L^2}\) & \(\|T-T_{\mathrm{ref}}\|_{L^2}\) & \(\Eorig(1.6)\) & \(\max_n\Delta\Eorig^{\mathrm{step}}\) (inc.) & \(\max|\Delta\mathcal{H}|\) \\
\midrule
Proposed CSRK & \(0.005\) & \(31.98\) & \(1.812\times 10^{-2}\) & \(6.235\times 10^{-3}\) & \(4{,}678.09\) & \(-3.365\times 10^{-1}\) (0) & \(5.55\times 10^{-16}\) \\
Proposed CSRK & \(0.025\) & \(6.58\)  & \(8.629\times 10^{-2}\) & \(3.016\times 10^{-2}\) & \(4{,}707.47\) & \(-1.267\times 10^{0}\) (0)  & \(4.44\times 10^{-16}\) \\
Proposed CSRK & \(0.050\) & \(3.92\)  & \(1.229\times 10^{-1}\) & \(4.392\times 10^{-2}\) & \(4{,}727.42\) & \(-1.954\times 10^{0}\) (0)  & \(5.55\times 10^{-16}\) \\
Proposed CSRK & \(0.200\) & \(1.02\)  & \(1.698\times 10^{-1}\) & \(6.438\times 10^{-2}\) & \(4{,}762.55\) & \(-3.557\times 10^{0}\) (0)  & \(4.44\times 10^{-16}\) \\
\midrule
SAV--BDF2 (\(S=0\)) & \(0.005\) & \(3.84\) & \(2.425\times 10^{-1}\) & \(4.898\times 10^{-2}\) & \(4{,}582.29\) & \(+2.119\times 10^{1}\) (9)  & \(3.27\times 10^{-2}\) \\
SAV--BDF2 (\(S=0\)) & \(0.025\) & \(0.69\) & \(4.936\times 10^{-1}\) & \(1.528\times 10^{-1}\) & \(7{,}876.08\) & \(+1.530\times 10^{3}\) (7)  & \(9.15\times 10^{-1}\) \\
SAV--BDF2 (\(S=0\)) & \(0.050\) & \(0.35\) & \(6.578\times 10^{-1}\) & \(2.025\times 10^{-1}\) & \(20{,}144.6\) & \(+4.536\times 10^{3}\) (8)  & \(2.89\times 10^{0}\) \\
SAV--BDF2 (\(S=0\)) & \(0.200\) & \(0.10\) & \(9.359\times 10^{-1}\) & \(1.291\times 10^{-1}\) & \(1.68\times 10^{5}\) & \(+4.683\times 10^{4}\) (7) & \(2.17\times 10^{2}\) \\
\midrule
SAV--BDF2 (\(S=2\)) & \(0.005\) & \(3.71\) & \(9.377\times 10^{-4}\) & \(3.205\times 10^{-4}\) & \(4{,}671.06\) & \(-3.526\times 10^{-1}\) (0) & \(2.18\times 10^{-6}\) \\
SAV--BDF2 (\(S=2\)) & \(0.025\) & \(0.77\) & \(1.027\times 10^{-2}\) & \(3.614\times 10^{-3}\) & \(4{,}675.63\) & \(-1.570\times 10^{0}\) (0)  & \(4.72\times 10^{-5}\) \\
SAV--BDF2 (\(S=2\)) & \(0.050\) & \(0.48\) & \(3.079\times 10^{-2}\) & \(1.086\times 10^{-2}\) & \(4{,}684.16\) & \(-2.536\times 10^{0}\) (0)  & \(1.67\times 10^{-4}\) \\
SAV--BDF2 (\(S=2\)) & \(0.200\) & \(0.09\) & \(1.017\times 10^{-1}\) & \(3.951\times 10^{-2}\) & \(4{,}723.06\) & \(-5.846\times 10^{0}\) (0)  & \(4.89\times 10^{-4}\) \\
\midrule
SAV--BDF2 (\(S=4\)) & \(0.005\) & \(3.69\) & \(3.832\times 10^{-4}\) & \(1.293\times 10^{-4}\) & \(4{,}671.28\) & \(-3.521\times 10^{-1}\) (0) & \(2.07\times 10^{-6}\) \\
SAV--BDF2 (\(S=4\)) & \(0.025\) & \(0.77\) & \(1.982\times 10^{-2}\) & \(6.804\times 10^{-3}\) & \(4{,}679.27\) & \(-1.386\times 10^{0}\) (0)  & \(4.36\times 10^{-5}\) \\
SAV--BDF2 (\(S=4\)) & \(0.050\) & \(0.37\) & \(4.720\times 10^{-2}\) & \(1.641\times 10^{-2}\) & \(4{,}690.89\) & \(-2.164\times 10^{0}\) (0)  & \(1.44\times 10^{-4}\) \\
SAV--BDF2 (\(S=4\)) & \(0.200\) & \(0.10\) & \(1.219\times 10^{-1}\) & \(4.644\times 10^{-2}\) & \(4{,}734.22\) & \(-4.640\times 10^{0}\) (0)  & \(2.08\times 10^{-4}\) \\
\bottomrule
\end{tabular}%
}
\end{table}

\FloatBarrier
\subsection{Example 2: Two-dimensional moving-interface \texorpdfstring{\(hp\)}{hp} dendritic growth \texorpdfstring{(\(m=4\), \(m=6\), and asymmetric full-domain)}{(m=4, m=6, and asymmetric full-domain)}}\label{sec:exp2}
We next demonstrate the full moving-interface \([-0.98,0.98]\) \(hp\) strategy (Algorithm~\ref{alg:adaptive-csrk}) on strong-anisotropy (\(\mu_A>0\)) two-dimensional dendritic solidification benchmarks modeled after the physical regimes in \cite{lithesis,li2022}, here formulated strictly under homogeneous Neumann (insulated) boundary conditions rather than periodic boundaries. For the axis-aligned fourfold (\(m=4,\ \theta_0=0\)) and sixfold (\(m=6,\ \theta_0=\pi/2\)) crystals, \(\cos(m(\vartheta-\theta_0)) = \pm\cos(m\vartheta)\) is invariant under both \(x\mapsto -x\) and \(y\mapsto -y\), so solving on the first quadrant \([0,L]^2\) with homogeneous Neumann conditions produces a reflection-symmetric solution on \([-L,L]^2\) (mirrored across the axes for visualization); we subsequently remove all symmetry assumptions to test rotated (\(\theta_0=\pi/8\)) and multi-nuclei configurations directly on the full domain \([-1,1]^2\).

\paragraph{Fourfold needle dendrite and latent-heat sensitivity (\(m=4,\ \sigma=0.10,\ \mu_A=0.55\)).}
On \(\Omega=[-1,1]^2\) we set \((\varrho,\varepsilon,\sigma,\lambda,D)=(1000,\,0.015,\,0.10,\,400,\,2.5\times 10^{-3})\) with initial nucleus
\[
 \phi_0(\bm{x}) = \tanh\!\left(\frac{0.04 - |\bm{x}|}{\sqrt{2}\cdot 0.015}\right),
 \qquad T_0(\bm{x}) = -0.60\,\frac{1-\phi_0(\bm{x})}{2},
\]
and study three latent-heat parameters \(K\in\{0.8, 1.0, 1.2\}\). Since \(\sigma=0.10 > 1/15\), the fourfold surface-energy density is genuinely nonconvex and activates the global convex-splitting shift \(\mu_A = (1+0.10)(15\cdot 0.10 - 1) = 0.55 > 0\) from \eqref{eq:muA-2d}. We employ time step \(k=0.025\), a \(2{:}1\)-balanced quadtree with \(L_{\min}=2\), \(L_{\max}=6\) (\(h_{\min}=1/64=0.015625\)), polynomial degrees \((p_{\mathrm{low}},p_{\mathrm{mid}},p_{\mathrm{high}})=(1,2,3)\), one-ring \(Q_2\) halo, and adaptation every \(6\) steps. The equivalent uniform fine \(L_6\,Q_3\) discretization requires \(4{,}096\) elements and \(N_{\mathrm{unif}}=37{,}249\) conforming unknowns per quadrant (\(148{,}225\) on \([-1,1]^2\)).

\begin{figure}[htbp]
\centering
\includegraphics[width=\linewidth]{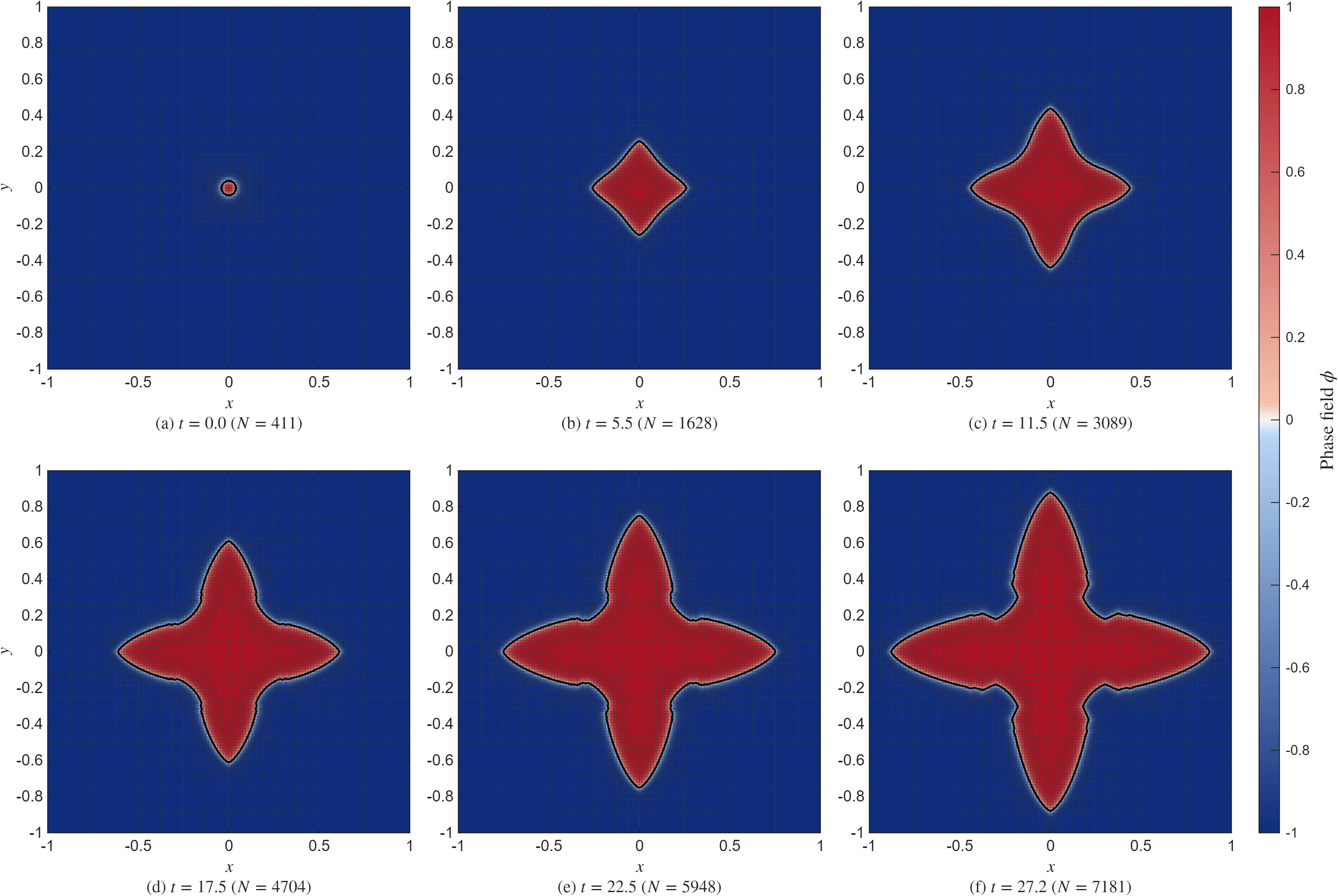}
\caption{Example~2 (\(m=4,\ \sigma=0.10,\ \mu_A=0.55,\ K=1.0\)): Formation process of the fourfold needle dendrite on the insulated square \([-1,1]^2\) from \(t=0\) to \(t=27.2\), showing the phase field \(\phi\), the \(\phi=0\) solid--liquid interface (black contour), and the moving \(2{:}1\)-balanced quadtree mesh. Quadrant conforming unknown counts \(N\) are shown below each panel (versus \(37{,}249\) on the equivalent uniform \(L_6\,Q_3\) mesh).}
\label{fig:exp2-m4-process}
\end{figure}

\begin{figure}[htbp]
\centering
\includegraphics[width=\linewidth]{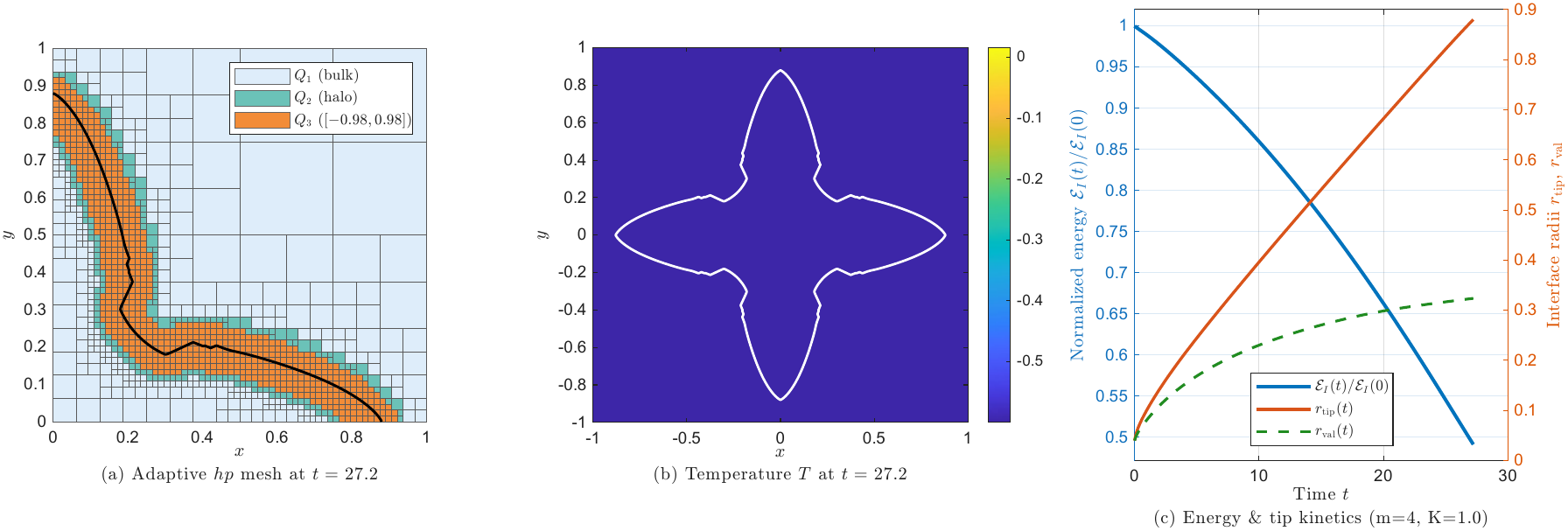}
\caption{Example~2 (\(m=4,\ K=1.0\)): Final-state \(hp\) adaptivity and thermodynamic diagnostics at \(t=27.225\). (a)~Quadrant \(2{:}1\)-balanced quadtree colored by polynomial degree (\(Q_1\) bulk, \(Q_2\) thermal halo, \(Q_3\) on \([-0.98,0.98]\)) with \(\phi=0\) contour. (b)~Full-domain temperature field \(T\) with \(\phi=0\) contour (white). (c)~Strict monotonic decay of the normalized original discrete energy \(\Eorig(t)/\Eorig(0)\) alongside tip and valley radii \(r_{\mathrm{tip}}(t), r_{\mathrm{val}}(t)\).}
\label{fig:exp2-m4-detail}
\end{figure}

Figure~\ref{fig:exp2-m4-process} illustrates the complete dendritic formation process for \(K=1.0\) across \(1{,}089\) steps to \(t=27.225\) (quadtree elements growing from \(115\) to \(1{,}204\) vs.\ \(4{,}096\) uniform): starting from the compact circular seed (\(r_{\mathrm{tip}}=0.040\), \(N_{\mathrm{conf}}=411\), a \(90.6\times\) initial DOF reduction), the crystal develops a faceted diamond (\(t=5.5\)), sharpens into four slender primary needle arms (\(t=11.5, 17.5\)), and by \(t=22.5\)--\(27.225\) reaches \(r_{\mathrm{tip}}=0.8800\) while valley growth stalls at \(r_{\mathrm{val}}=0.3238\) (tip-to-valley ratio \(2.718\)) with clearly visible nascent secondary sidebranching bulges along all four arm flanks. Figure~\ref{fig:exp2-m4-detail} displays the final \(2{:}1\)-balanced \(hp\) quadtree (orange \(Q_3\) ribbon tightly hugging both the needle tips and the nascent sidebranch bulges, framed by a teal \(Q_2\) halo and coarse light-blue \(Q_1\) bulk), the corresponding full-domain temperature field showing thin thermal boundary layers ahead of the four advancing tips and latent heat trapped inside the valleys, and the strictly monotonic decay of the original discrete energy \(\Eorig(0)=4{,}793.33\to\Eorig(27.225)=2{,}354.38\) (\(\max\Delta\Eorig^{\mathrm{step}}=-1.265\), max Newton residual \(9.97\times 10^{-11}\), \(\min\alpha_T=1-1.41\times 10^{-8}\) across \(180\) adaptations, and \(\max|\Delta\mathcal{H}|=3.33\times 10^{-16}\)).

\begin{figure}[htbp]
\centering
\includegraphics[width=\linewidth]{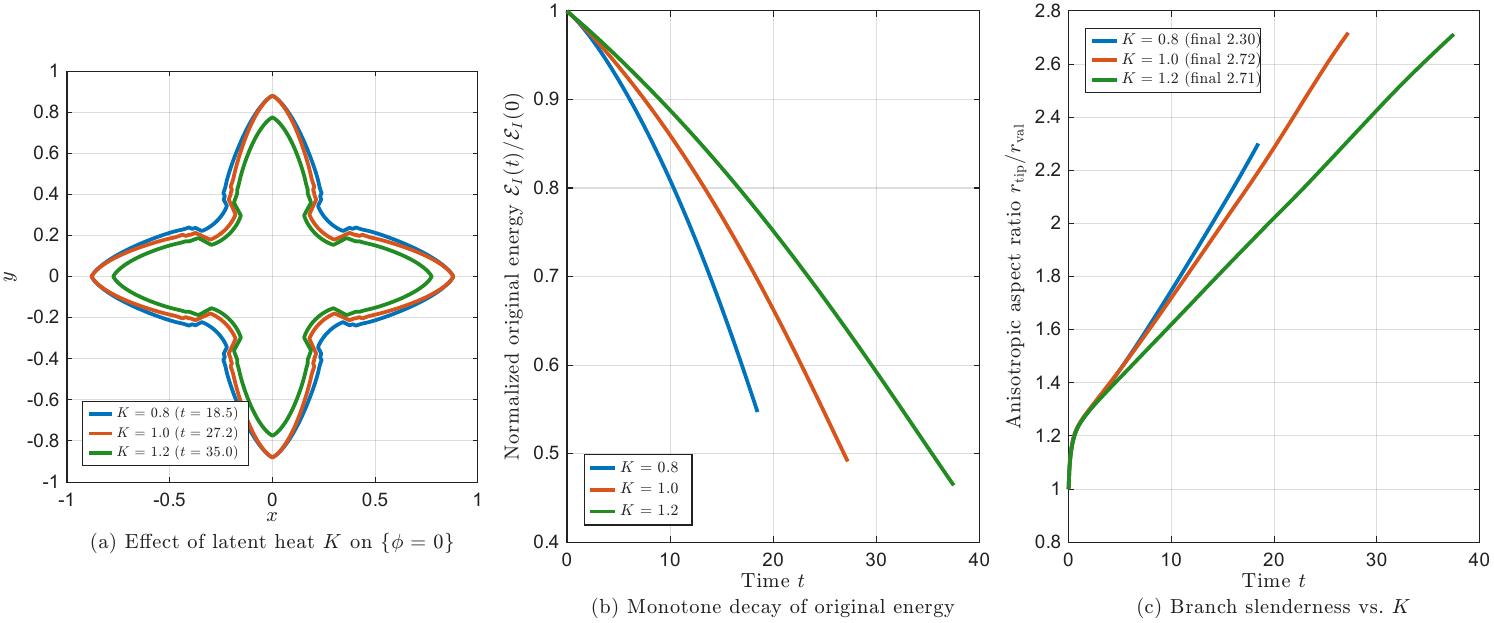}
\caption{Example~2 (\(m=4,\ \sigma=0.10\)): Effect of the latent-heat parameter \(K\in\{0.8,1.0,1.2\}\) under homogeneous Neumann boundary conditions. (a)~Overlaid \(\phi=0\) contours from the last saved field snapshot of each run (\(t=18.5\) for \(K=0.8\), \(t=27.2\) for \(K=1.0\), and \(t=35.0\) for \(K=1.2\)) illustrating branch narrowing and tip sharpening as \(K\) increases. (b)~Strictly monotonic decay of the original discrete energy \(\Eorig(t)/\Eorig(0)\) over the full trajectories (to \(t=37.5\) for \(K=1.2\)). (c)~Evolution of the anisotropic aspect ratio \(r_{\mathrm{tip}}/r_{\mathrm{val}}\).}
\label{fig:exp2-m4-K}
\end{figure}

Figure~\ref{fig:exp2-m4-K} compares \(K\in\{0.8,1.0,1.2\}\): increasing latent heat \(K\) slows primary growth (reaching \(r_{\mathrm{tip}}=0.880\) at \(t=18.475\) for \(K=0.8\) and \(t=27.225\) for \(K=1.0\), and \(r_{\mathrm{tip}}=0.815\) at \(t=37.50\) for \(K=1.2\)) while substantially narrowing the valleys (\(r_{\mathrm{val}}=0.3825,\,0.3238,\,0.3004\)) and sharpening the dendritic branches, reproducing the physical latent-heat dependence of \cite{lithesis} while rigorously dissipating the \emph{original} discrete energy \(\Eorig(t)\) in every case (all three runs satisfy \(\max|\Delta\mathcal{H}|=3.33\times 10^{-16}\)).

\begin{figure}[htbp]
\centering
\includegraphics[width=\linewidth]{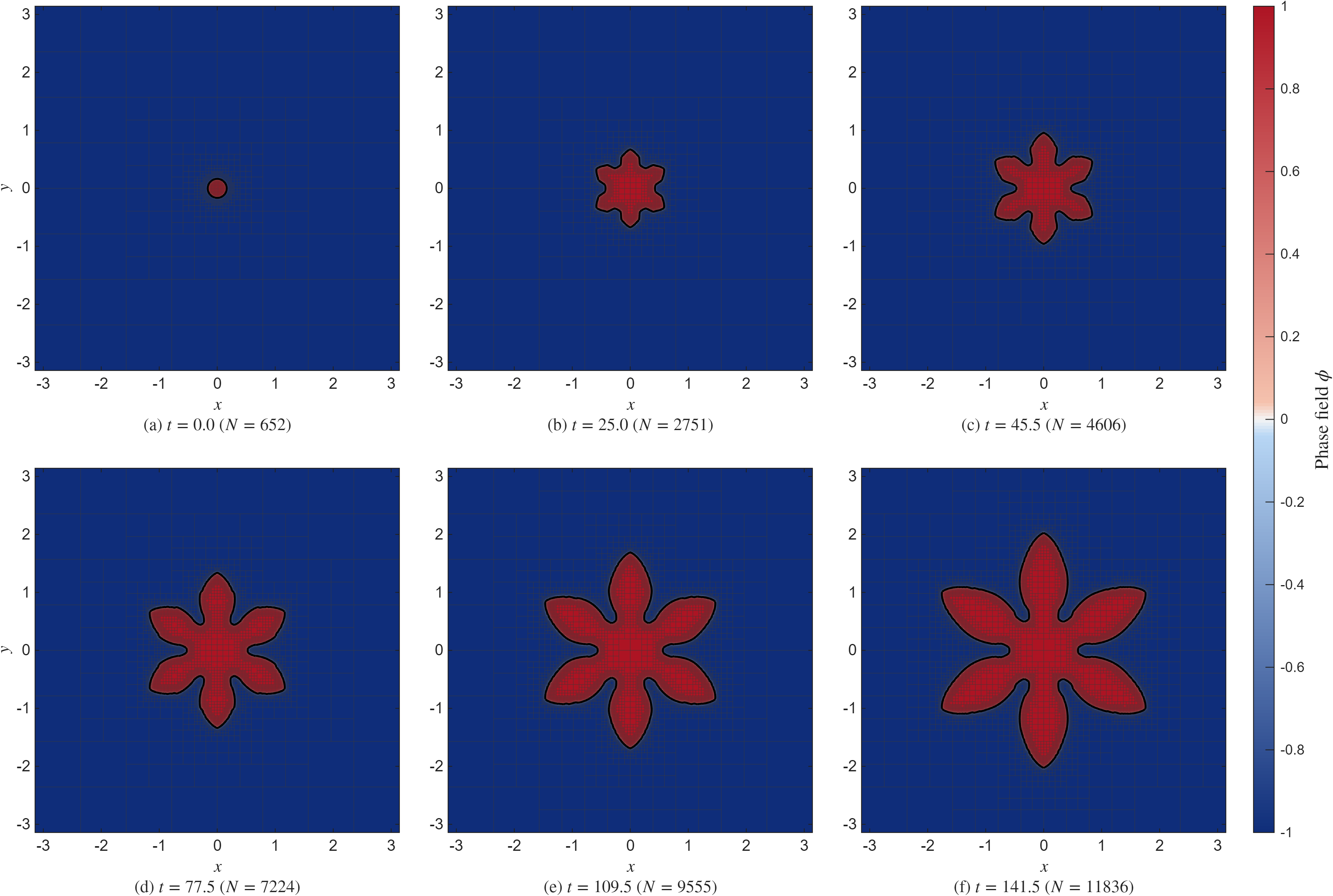}
\caption{Example~2 (\(m=6,\ \sigma=0.04,\ \mu_A=0.416,\ K=0.90\)): Formation process of the sixfold snowflake dendrite on the insulated square \([-\pi,\pi]^2\) from \(t=0\) to \(t=141.5\), showing \(\phi\), the \(\phi=0\) interface, and the moving \(2{:}1\)-balanced \(L_7\,Q_3\) quadtree mesh. Quadrant conforming unknown counts \(N\) are shown below each panel (versus \(148{,}225\) quadrant DOFs on the equivalent uniform \(L_7\,Q_3\) mesh).}
\label{fig:exp2-m6-process}
\end{figure}

\begin{figure}[htbp]
\centering
\includegraphics[width=\linewidth]{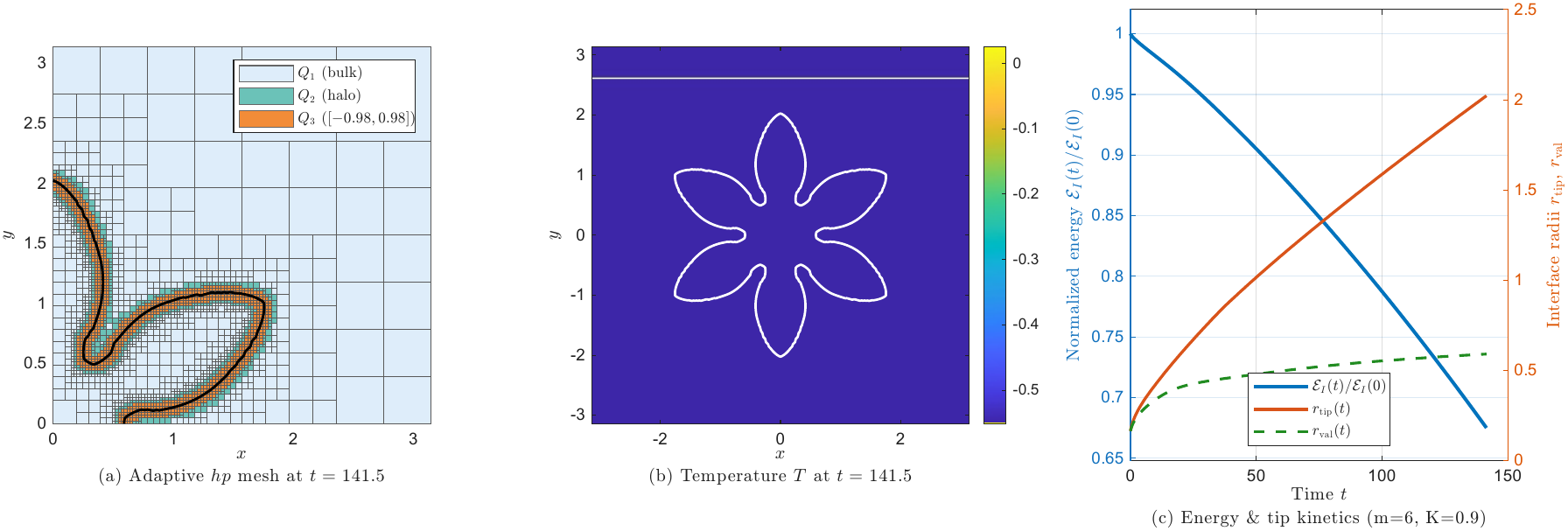}
\caption{Example~2 (\(m=6,\ K=0.90\)): Final-state \(hp\) adaptivity and diagnostics at \(t=141.5\). (a)~Degree-colored \(2{:}1\)-balanced quadtree tracking both the \(90^\circ\) and \(30^\circ\) snowflake arms. (b)~Full-domain temperature distribution \(T\). (c)~Monotonic decay of the original discrete energy \(\Eorig(t)/\Eorig(0)\) and tip/valley radius evolution.}
\label{fig:exp2-m6-detail}
\end{figure}

\paragraph{Sixfold snowflake dendritic growth (\(m=6,\ \sigma=0.04,\ \mu_A=0.416\)).}
On the larger insulated square \(\Omega=[-\pi,\pi]^2\), we adopt the snow-crystal parameters \((\varrho,\varepsilon,\sigma,\theta_0,\lambda,D,K)=(600,\,0.015,\,0.04,\,\pi/2,\,360,\,2.25\times 10^{-3},\,0.90)\) with initial seed
\[
 \phi_0(\bm{x}) = \tanh\!\left(\frac{0.027 - |\bm{x}|^2}{0.0054}\right),
 \qquad T_0(\bm{x}) = -0.55\,\frac{1-\phi_0(\bm{x})}{2},
\]
activating the nonconvex sixfold convex-splitting shift \(\mu_A = (1+0.04)(35\cdot 0.04 - 1) = 0.416 > 0\). To resolve the thin interfacial width \(\varepsilon/\pi \approx 4.8\times 10^{-3}\), we employ a seven-level \(2{:}1\)-balanced quadtree (\(L_{\min}=2\), \(L_{\max}=7\), \(h_{\min}=\pi/128\approx 0.0245\)) with \(Q_1/Q_2/Q_3\) degrees, a one-ring \(Q_2\) halo, and adaptation every \(6\) steps; the equivalent uniform \(L_7\,Q_3\) mesh requires \(16{,}384\) elements and \(N_{\mathrm{unif}}=148{,}225\) unknowns per quadrant (\(591{,}361\) on \([-\pi,\pi]^2\), \(1.5\times\) finer per coordinate than a \(512^2\) periodic Fourier grid). Advancing \(4{,}100\) steps (\(k=0.025\) for \(t\in[0,37.5]\) and \(k=0.04\) for \(t\in[37.5,141.5]\)) across \(680\) mesh adaptations (quadtree elements growing from \(187\) to \(2{,}398\) vs.\ \(16{,}384\) uniform) yields the six-arm snowflake evolution in Figure~\ref{fig:exp2-m6-process}: the six primary arms elongate dramatically from \(r_{\mathrm{tip}}=0.1644\) to \(r_{\mathrm{tip}}=2.0222\) while the valleys remain thermally arrested at \(r_{\mathrm{val}}=0.5907\), achieving an anisotropic aspect ratio of \(r_{\mathrm{tip}}/r_{\mathrm{val}}=3.423\). Figure~\ref{fig:exp2-m6-detail} confirms that the \(Q_3\) band tightly follows both the vertical (\(\theta=\pi/2\)) and oblique (\(\theta=\pi/6\)) arms inside the quadrant, while \(\Eorig\) decays strictly monotonically from \(39{,}730.09\) to \(26{,}820.99\) (\(\max\Delta\Eorig^{\mathrm{step}}=-1.610\), max Newton residual \(4.60\times 10^{-11}\), \(\min\alpha_T=1-1.94\times 10^{-8}\), \(\max|\Delta\mathcal{H}|=3.33\times 10^{-15}\)).

\begin{figure}[htbp]
\centering
\includegraphics[width=\linewidth]{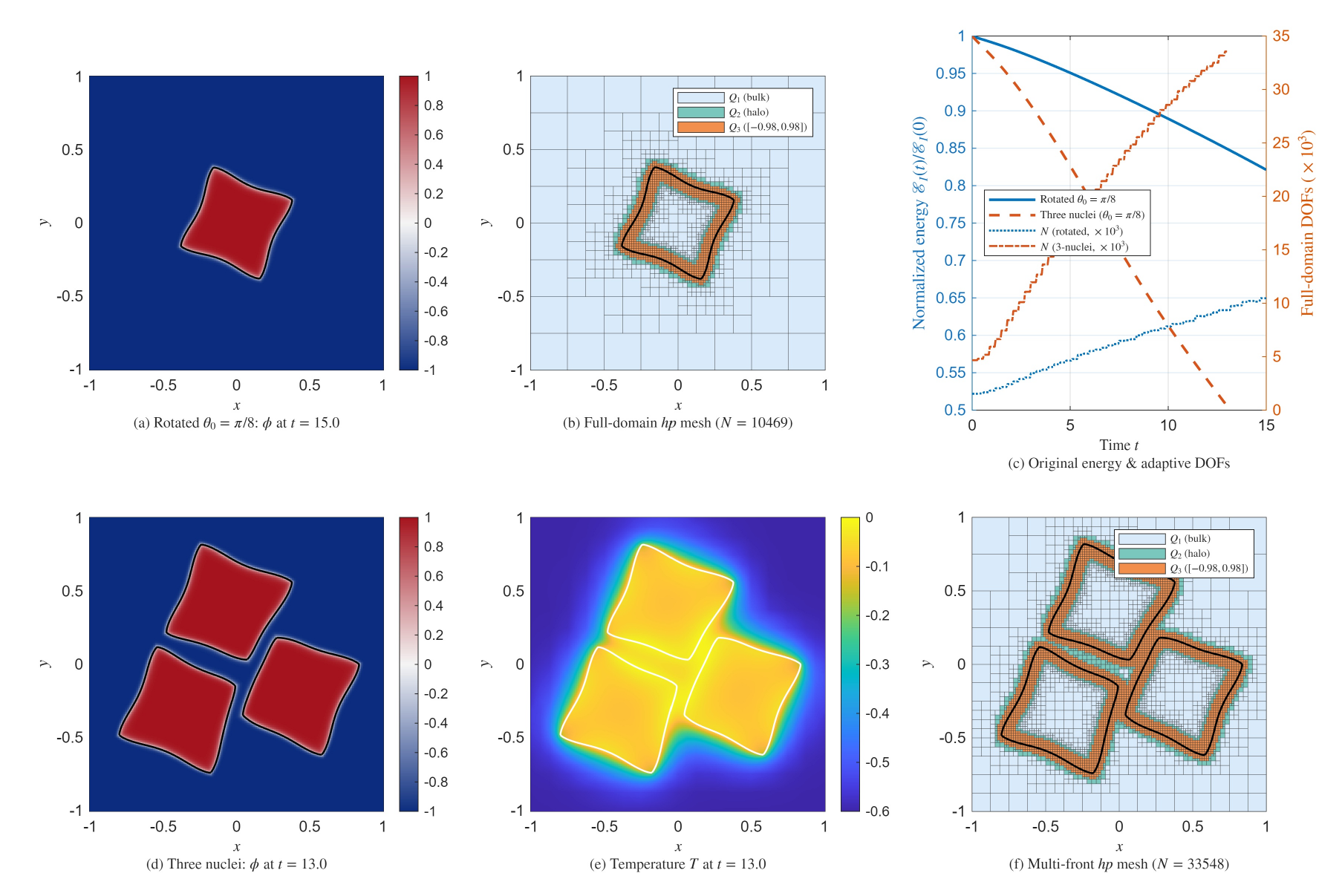}
\caption{Example~2 (Full-domain asymmetric 2D growth on \([-1,1]^2\) without symmetry reduction, \(m=4,\ \sigma=0.10,\ \mu_A=0.55,\ \theta_0=\pi/8\)): (a)--(b)~Rotated single crystal \(\phi\) and degree-colored \(2{:}1\)-balanced \(hp\) quadtree at \(t=15.0\) (\(K=1.0,\ N_{\mathrm{conf}}=10{,}469\), agreeing with the \(\theta_0=0\) reference to \(0.003\%\) in tip radius). (c)~Monotonic original-energy decay \(\Eorig(t)/\Eorig(0)\) and adaptive full-domain DOF counts for both asymmetric cases. (d)--(f)~Three asymmetrically placed rotated nuclei at \(t=13.0\) (\(K=0.8,\ N_{\mathrm{conf}}=33{,}548\)): phase field \(\phi\), temperature field \(T\) showing inter-crystal thermal impingement, and multi-front \(Q_3/Q_2/Q_1\) quadtree.}
\label{fig:exp2-asymmetric}
\end{figure}

\paragraph{Full-domain asymmetric growth: rotated crystal orientation and multi-nuclei competition.}
To confirm that neither the \(2{:}1\)-balanced quadtree nor the natural homogeneous Neumann formulation \eqref{eq:bc-neumann} relies on axis-aligned symmetry or quadrant reduction, we solve two genuinely non-symmetric strong-anisotropy (\(m=4,\ \sigma=0.10,\ \mu_A=0.55\)) problems directly on the full insulated square \(\Omega=[-1,1]^2\) (\(L_{\min}=2,\ L_{\max}=7\) on \([-1,1]^2\), matching the exact \(h_{\min}=1/64\) cell size and \(N_{\mathrm{unif}}=148{,}225\) full-domain DOFs of the uniform \(L_6\,Q_3\) quadrant mesh, with \(k=0.05\) and adaptation every \(6\) steps):
\begin{itemize}
\item \textbf{Case~A (Rotated crystal, \(\theta_0=\pi/8=22.5^\circ,\ K=1.0,\ t=15.0\)).} Rotating the preferred growth axes to \(22.5^\circ\) off the Cartesian mesh lines breaks reflection symmetry across both axes. Across \(300\) steps and \(49\) adaptations (Figure~\ref{fig:exp2-asymmetric}(a)--(c)), full-domain conforming unknowns grow from \(1{,}529\) (\(96.9\times\) initial compression) to \(10{,}469\) (\(14.2\times\) final, \(24.1\times\) average compression vs.\ \(148{,}225\) uniform) in \(366.3\)\,s (\(1.22\)\,s/step; adaptation overhead \(18.1\%\)), while \(\Eorig\) decays strictly monotonically from \(19{,}173.34\) to \(15{,}745.22\) (\(\max\Delta\Eorig^{\mathrm{step}}=-7.80\), \(\min\alpha_T=1-2.40\times 10^{-8}\), \(\max|\Delta\mathcal{H}|=2.89\times 10^{-15}\)). All four rotated arm tips reach identical lengths \(r_{\mathrm{tip}}=0.41043\) with inter-arm spread \(1.50\times 10^{-15}\). Compared against an axis-aligned (\(\theta_0=0\)) reference at identical \(k=0.05\) and \(t=15.0\) (\(r_{\mathrm{tip}}=0.41044,\ r_{\mathrm{val}}=0.26985\)), the rotated full-domain \(hp\) solution agrees to \(0.0031\%\) in tip radius and \(0.049\%\) in valley radius (\(r_{\mathrm{val}}=0.26999\)), confirming negligible directional mesh bias on the tested nonmatching \(2{:}1\) quadtree prior to boundary impingement.
\item \textbf{Case~B (Three asymmetrically placed rotated nuclei, \(\theta_0=\pi/8,\ K=0.8,\ t=13.0\)).} Three unequal circular seeds centered at \((-0.36,-0.30)\), \((0.40,-0.18)\), and \((-0.06,0.38)\) with radii \(0.045, 0.040, 0.042\) evolve simultaneously under full-domain insulated Neumann BCs across \(260\) steps and \(43\) adaptations (\(1{,}287.4\)\,s total, \(4.95\)\,s/step, adaptation share \(15.5\%\)). Figure~\ref{fig:exp2-asymmetric}(d)--(f) shows latent heat rejected into the narrow inter-crystal channels thermally arresting the mutually facing arms while outward-facing tips accelerate into the undercooled melt, as the \(Q_3/Q_2/Q_1\) mesh dynamically tracks all three interfaces simultaneously (\(N_{\mathrm{conf}}=4{,}673\to 33{,}548\), \(7.4\times\) average DOF compression) with strict original-energy monotonicity (\(\Eorig: 23{,}877.00\to 12{,}072.94,\ \max\Delta\Eorig^{\mathrm{step}}=-30.82\)), \(\min\alpha_T=1-2.09\times 10^{-8}\), and \(\max|\Delta\mathcal{H}|=2.44\times 10^{-15}\).
\end{itemize}

\paragraph{Computational efficiency, adaptation overhead, and uniform-grid baselines.}
On the same workstation and Newton tolerance (\(10^{-10}\)), we timed both matched uniform fine meshes over direct benchmark windows: uniform \(L_6\,Q_3\) (\(37{,}249\) quadrant DOFs, timed over \(10\) steps) requires \(7.325\)\,s/step on the \(m=4\) problem, and uniform \(L_7\,Q_3\) (\(148{,}225\) quadrant DOFs, timed over \(3\) steps) requires \(29.11\)\,s/step on the \(m=6\) problem (which would extrapolate to \(33.2\)\,h for \(4{,}100\) steps). Because the adaptive \(hp\) mesh preserves identical \(L_{\max}\,Q_3\) resolution across the entire \([-0.98,0.98]\) transition band and \(2{:}1\)-balanced \(Q_2\) thermal halo while coarsening only where \(|\phi\pm 1|<0.02\), it resolves the moving interface at the same local spectral accuracy. As summarized in Table~\ref{tab:exp2-efficiency}, the adaptive \(2{:}1\)-balanced \(hp\) solver achieves \(10.0\times\)--\(10.7\times\) average DOF reduction and \(12.5\times\)--\(13.3\times\) full-trajectory wall-clock speedup on \(m=4\) (\(19.0\times\) over the first \(700\) steps; mesh adaptation, balance closure, and \(L^2\) transfer account for \(20.4\%\)--\(21.2\%\) of total runtime), and \(25.1\times\) average (\(227.3\times\) initial) DOF reduction and \(29.7\times\) full-trajectory (\(79.0\times\) early-stage) wall-clock speedup on the high-resolution \(m=6\) snowflake benchmark (with all \(680\) adaptations accounting for \(23.2\%\) of the \(4{,}016.4\)\,s total wall time).

\begin{table}[htbp]
\centering
\caption{Example~2: Verified summary of strong-anisotropy (\(\mu_A>0\)) 2D moving-interface \(hp\) dendritic growth under homogeneous Neumann BCs and runtime comparison against short-window uniform fine spectral-element timing baselines on identical hardware. Reported conforming unknown counts \(N_{\mathrm{conf}}\) are per field; energies \(\Eorig\) and DOFs are evaluated on the computational quadrant \([0,L]^2\) for the four symmetric runs and on the full square \([-1,1]^2\) (\(N_{\mathrm{unif}}=148{,}225\)) for the two asymmetric runs.}
\label{tab:exp2-efficiency}
\footnotesize
\setlength{\tabcolsep}{4pt}
\resizebox{\linewidth}{!}{%
\begin{tabular}{@{}lccccccc@{}}
\toprule
Case (Domain, Grid) & Steps (\(t_{\mathrm{end}}\)) & \(N_{\mathrm{conf}}\) (init / avg / end) & Avg.\ DOF saving & \(r_{\mathrm{tip}}/r_{\mathrm{val}}\) & \(\Eorig(0)\to\Eorig(t_{\mathrm{end}})\) & \(\max|\Delta\mathcal{H}|\) & Avg.\ time/step \\
\midrule
\(m=4,\ K=0.8\) (\(L_6\,Q_3\)) & 739 (18.48) & \(411\) / \(3{,}472\) / \(6{,}838\) & \(10.73\times\) & 2.301 & \(5{,}990.7\to 3{,}277.1\) & \(3.33\times 10^{-16}\) & \(0.550\)\,s (\(13.3\times\)) \\
\(m=4,\ K=1.0\) (\(L_6\,Q_3\)) & \(1{,}089\) (27.23) & \(411\) / \(3{,}714\) / \(7{,}181\) & \(10.03\times\) & 2.718 & \(4{,}793.3\to 2{,}354.4\) & \(3.33\times 10^{-16}\) & \(0.586\)\,s (\(12.5\times\)) \\
\(m=4,\ K=1.2\) (\(L_6\,Q_3\)) & \(1{,}500\) (37.50) & \(411\) / \(3{,}584\) / \(6{,}682\) & \(10.39\times\) & 2.712 & \(3{,}995.1\to 1{,}854.3\) & \(3.33\times 10^{-16}\) & \(0.559\)\,s (\(13.1\times\)) \\
Uniform \(L_6\,Q_3\) (\(m=4\)) & -- & \(37{,}249\) (fixed) & \(1.00\times\) & -- & -- & \(1.8\times 10^{-15}\) & \(7.325\)\,s (\(1.0\times\)) \\
\midrule
\(m=6,\ K=0.9\) (\(L_7\,Q_3\)) & \(4{,}100\) (141.5) & \(652\) / \(5{,}911\) / \(11{,}836\) & \(25.08\times\) & 3.423 & \(39{,}730.1\to 26{,}821.0\) & \(3.33\times 10^{-15}\) & \(0.980\)\,s (\(29.7\times\)) \\
Uniform \(L_7\,Q_3\) (\(m=6\)) & -- & \(148{,}225\) (fixed) & \(1.00\times\) & -- & -- & \(2.3\times 10^{-14}\) & \(29.11\)\,s (\(1.0\times\)) \\
\midrule
Full \([-1,1]^2\), rot.\ \(\theta_0=\pi/8\) & 300 (15.00) & \(1{,}529\) / \(6{,}147\) / \(10{,}469\) & \(24.11\times\) & 1.520 & \(19{,}173.3\to 15{,}745.2\) & \(2.89\times 10^{-15}\) & \(1.221\)\,s \\
Full \([-1,1]^2\), 3 nuclei (\(\pi/8\)) & 260 (13.00) & \(4{,}673\) / \(19{,}982\) / \(33{,}548\) & \(7.42\times\) & -- & \(23{,}877.0\to 12{,}072.9\) & \(2.44\times 10^{-15}\) & \(4.952\)\,s \\
\bottomrule
\end{tabular}%
}
\end{table}

\paragraph{Equal-condition spatial fidelity and component ablation study.}
To quantify spatial accuracy against an \emph{independent} higher-order spatial reference under identical physical parameters, initial data, homogeneous Neumann boundary conditions, and CSRK step size (\(k=0.025\)), and to isolate the individual contributions of \(2{:}1\)-balanced quadtree \(h\)-refinement, spectral \(p\)-enrichment, the \([-0.98,0.98]\) interface band, and the \(Q_2\) thermal halo, we integrate the \(m=4,\ \sigma=0.10,\ K=1.0\) dendrite for \(60\) steps (\(t=1.50\), adapting every \(6\) steps for adaptive cases, during which the nucleus triples in radius from \(0.040\) to \(r_{\mathrm{tip}}=0.1227,\ r_{\mathrm{val}}=0.0964\); each candidate projects the analytical initial data onto its own mesh and preserves its corresponding discrete mean enthalpy \(\bar{H}\)). All ten candidate discretizations---including the matched uniform fine \(L_6\,Q_3\) grid itself (\(37{,}249\) quadrant DOFs per field, \(384.89\)\,s)---are evaluated via exact common-quadrature \(L^2(\Omega)\) errors (\(692{,}224\) Gauss points), relative original-energy errors, and tip/valley radius errors against an independent higher-order uniform \(L_6\,Q_4\) spatial reference (\(4{,}096\) elements, \(66{,}049\) quadrant DOFs per field, \(809.15\)\,s).

Table~\ref{tab:exp2-ablation} yields three clear conclusions:
\begin{enumerate}
\item \textbf{Substantial accuracy advantage under equal cost/DOF.} At virtually identical runtime (\(5.97\)\,s vs.\ \(5.71\)\,s) and half the degrees of freedom (\(536.1\) vs.\ \(1{,}089\)), our full \(2{:}1\)-balanced \(hp\) method reduces the phase \(L^2\) error by \(117\times\) and the tip-radius error by \(101\times\) compared with uniform \(L_4\,Q_2\); even against uniform \(L_4\,Q_3\) (\(4.5\times\) more DOFs and \(2.6\times\) slower), it reduces the phase and tip-radius errors by \(38\times\) and \(41\times\), respectively. Relative to the matched uniform \(L_6\,Q_3\) grid, it retains \(2.25\times 10^{-4}\) phase \(L^2\) fidelity, \(4.09\times 10^{-5}\) relative energy accuracy, and \(9.18\times 10^{-5}\) tip-radius accuracy with \(69.5\times\) fewer average DOFs and a \(64.5\times\) speedup (and \(123\times\) fewer DOFs / \(136\times\) speedup versus the independent \(L_6\,Q_4\) reference).
\item \textbf{Synergy of \(h\)-tree and high-\(p\) enrichment is the primary driver.} Neither component alone suffices: pure \(p\)-refinement on a fixed coarse \(L_4\) grid (\(\|\phi-\phi_{\mathrm{ref}}\|_{L^2}=8.61\times 10^{-3}\), \(38\times\) worse than full \(hp\)) fails to resolve the \(O(\varepsilon)\) interfacial layer, while pure \(h\)-refinement with bilinear \(Q_1\) elements on the \(L_6\) quadtree still suffers a \(1.06\times 10^{-2}\) phase error (\(47\times\) worse than full \(hp\)). Combining \(2{:}1\)-balanced \(h\)-refinement with \(Q_3\) interface enrichment drops the error by nearly two orders of magnitude.
\item \textbf{Observed trade-offs of the \([-0.98,0.98]\) band, \(Q_2\) halo, and bulk \(p\)-coarsening.} For a \(\tanh(s/\varepsilon)\) transition layer, the \([-0.98,0.98]\) band spans \(|s|\le\operatorname{atanh}(0.98)\,\varepsilon\approx 2.30\,\varepsilon\) and captures \(99.94\%\) of the interfacial gradient energy density \(\operatorname{sech}^4(s/\varepsilon)\) (whereas \([-0.90,0.90]\) truncates \(1.45\%\), \(24\times\) more tail energy). Consequently, widening the band from \([-0.90,0.90]\) to \([-0.98,0.98]\) cuts the phase error by \(43\%\) (\(3.92\times 10^{-4}\to 2.25\times 10^{-4}\)) at only \(16\%\) DOF and \(8\%\) wall-time overhead. Adding the one-ring \(Q_2\) halo at \(12\%\) DOF cost lowers phase, energy, and valley errors by \(30\%\), \(22\%\), and \(39\%\), respectively (while the thermal \(L^2\) error is essentially unchanged at \(5.71\times 10^{-4}\) and the small \(O(10^{-5})\) tip-radius error rises slightly from \(7.76\times 10^{-5}\) to \(9.18\times 10^{-5}\)). Widening further to \([-0.999,0.999]\) gives another \(2.1\times\) phase-error reduction at \(32\%\) additional DOFs, placing \([-0.98,0.98]\) at a favorable accuracy--cost knee among the tested configurations. Finally, keeping \(Q_3\) everywhere on the \(2{:}1\)-balanced quadtree (\emph{\(h\)-only \(Q_3\)}) reproduces uniform \(L_6\,Q_3\) to within \(6\%\) in \(\phi\) (\(3.23\times 10^{-5}\) vs.\ \(3.04\times 10^{-5}\)) and to two significant digits in energy and tip/valley radii at \(38.0\times\) fewer DOFs and \(46.5\times\) speedup (with \(3.18\times\) larger bulk thermal error \(1.88\times 10^{-5}\) vs.\ \(5.92\times 10^{-6}\) due to outer \(h\)-coarsening). The remaining \(O(10^{-4})\) gap of the zoned \(Q_3/Q_2/Q_1\) scheme reflects primarily bulk \(p\)-coarsening in the wider thermal diffusion layer (together with mixed-degree trace transitions), offering an extra \(\sim 1.8\times\) DOF reduction in 2D and up to \(27\times\) asymptotic per-bulk-cell DOF compression in 3D while preserving \(O(10^{-4})\) fidelity.
\end{enumerate}

\begin{table}[htbp]
\centering
\caption{Example~2 (\(m=4,\ \sigma=0.10,\ K=1.0,\ k=0.025,\ t=1.50\)): Equal-condition spatial accuracy and component ablation study measured by exact common-quadrature errors (\(692{,}224\) Gauss points) against the independent higher-order uniform \(L_6\,Q_4\) spatial reference (\(66{,}049\) quadrant DOFs per field, \(809.15\)\,s).}
\label{tab:exp2-ablation}
\footnotesize
\setlength{\tabcolsep}{4pt}
\resizebox{\linewidth}{!}{%
\begin{tabular}{@{}lccccccc@{}}
\toprule
Spatial discretization & Avg.\ DOFs & Wall (s) & \(\|\phi-\phi_{\mathrm{ref}}\|_{L^2}\) & \(\|T-T_{\mathrm{ref}}\|_{L^2}\) & \(|\Delta\Eorig|/\Eorig^{\mathrm{ref}}\) & \(|r_{\mathrm{tip}}-r_{\mathrm{tip}}^{\mathrm{ref}}|\) & \(|r_{\mathrm{val}}-r_{\mathrm{val}}^{\mathrm{ref}}|\) \\
\midrule
Uniform \(L_4\,Q_2\) (coarse baseline) & \(1{,}089\) & \(5.71\) & \(2.626\times 10^{-2}\) & \(7.661\times 10^{-3}\) & \(1.852\times 10^{-3}\) & \(9.235\times 10^{-3}\) & \(4.155\times 10^{-3}\) \\
Uniform \(L_4\,Q_3\) (coarse baseline) & \(2{,}401\) & \(15.34\) & \(8.594\times 10^{-3}\) & \(2.400\times 10^{-3}\) & \(5.785\times 10^{-4}\) & \(3.742\times 10^{-3}\) & \(4.087\times 10^{-4}\) \\
\midrule
\(p\)-only (\(L_4\), zoned \(Q_3/Q_2/Q_1\)) & \(348.2\) & \(4.55\) & \(8.613\times 10^{-3}\) & \(2.402\times 10^{-3}\) & \(5.817\times 10^{-4}\) & \(3.748\times 10^{-3}\) & \(4.223\times 10^{-4}\) \\
\(h\)-only \(Q_1\) (\(2{:}1\) \(L_6\) tree, \(p\equiv 1\)) & \(120.6\) & \(2.27\) & \(1.058\times 10^{-2}\) & \(2.995\times 10^{-3}\) & \(7.754\times 10^{-4}\) & \(3.417\times 10^{-3}\) & \(2.417\times 10^{-3}\) \\
\(hp\), narrower band \([-0.90,0.90]\) & \(461.4\) & \(5.53\) & \(3.922\times 10^{-4}\) & \(5.608\times 10^{-4}\) & \(4.908\times 10^{-5}\) & \(1.177\times 10^{-4}\) & \(7.775\times 10^{-5}\) \\
\(hp\), no \(Q_2\) halo (\(\text{halo}=0\)) & \(478.4\) & \(5.20\) & \(3.218\times 10^{-4}\) & \(5.706\times 10^{-4}\) & \(5.267\times 10^{-5}\) & \(7.759\times 10^{-5}\) & \(9.036\times 10^{-5}\) \\
\textbf{Full \(hp\) (\([-0.98,0.98]\), \(Q_3/Q_2/Q_1\))} & \(\mathbf{536.1}\) & \(\mathbf{5.97}\) & \(\mathbf{2.252\times 10^{-4}}\) & \(\mathbf{5.710\times 10^{-4}}\) & \(\mathbf{4.089\times 10^{-5}}\) & \(\mathbf{9.178\times 10^{-5}}\) & \(\mathbf{5.497\times 10^{-5}}\) \\
\(hp\), wider band \([-0.999,0.999]\) & \(705.4\) & \(7.02\) & \(1.084\times 10^{-4}\) & \(3.157\times 10^{-4}\) & \(2.161\times 10^{-5}\) & \(5.287\times 10^{-5}\) & \(1.675\times 10^{-5}\) \\
\(h\)-only \(Q_3\) (\(2{:}1\) \(L_6\) tree, \(p\equiv 3\)) & \(981.2\) & \(8.29\) & \(3.226\times 10^{-5}\) & \(1.880\times 10^{-5}\) & \(1.237\times 10^{-6}\) & \(9.576\times 10^{-6}\) & \(2.111\times 10^{-6}\) \\
\midrule
Uniform \(L_6\,Q_3\) (matched fine grid) & \(37{,}249\) & \(384.89\) & \(3.045\times 10^{-5}\) & \(5.919\times 10^{-6}\) & \(1.248\times 10^{-6}\) & \(9.825\times 10^{-6}\) & \(2.158\times 10^{-6}\) \\
Uniform \(L_6\,Q_4\) (independent ref.) & \(66{,}049\) & \(809.15\) & \(0\) (ref) & \(0\) (ref) & \(0\) (ref) & \(0\) (ref) & \(0\) (ref) \\
\bottomrule
\end{tabular}%
}
\end{table}

\subsection{Example 3: Three-dimensional cubic-anisotropy \texorpdfstring{\(hp\)}{hp} dendritic growth}\label{sec:exp3}
Finally, we test the full three-dimensional matrix-free sum-factorized \(hp\)-SEM solver and \(2{:}1\)-balanced octree adaptivity on three-dimensional axis-aligned (\(\mathcal{R}=I\)) cubic dendritic solidification in the insulated cube \(\Omega=[-\pi,\pi]^3\) with homogeneous Neumann boundary conditions (replacing the periodic boundary setup of \cite{lithesis}). By reflection symmetry of \eqref{eq:aniso-3d} across the three coordinate planes, solving on the octant \([0,\pi]^3\) with homogeneous Neumann conditions on all six faces yields a reflection-symmetric solution on the full cube \([-\pi,\pi]^3\).

We set the physical parameters to the strong nonconvex cubic-anisotropy regime of \cite{lithesis}:
\begin{align*}
 (\varrho,\varepsilon,\delta,\lambda,D,K) &= (1000,\,0.015,\,0.10,\,400,\,2.5\times 10^{-3},\,0.90), \\
 \phi_0(\bm{x}) &= \tanh\!\left(\frac{0.027-|\bm{x}|^2}{0.0054}\right),
 \qquad T_0(\bm{x}) = -0.55\,\frac{1-\phi_0(\bm{x})}{2}.
\end{align*}
Because \(\delta=0.10 > 1/15\), the 3D cubic surface-energy density is nonconvex and activates the global convex-splitting stabilization \(\mu_A^{(3\mathrm{D})}=(1+0.10)(15\cdot 0.10 - 1)=0.55>0\) in \eqref{eq:muA-3d}. We configure the \(2{:}1\)-balanced octree with \(L_{\min}=1\), \(L_{\max}=5\) (\(h_{\min}=\pi/32\approx 0.0982\)), polynomial degrees \((p_{\mathrm{low}},p_{\mathrm{mid}},p_{\mathrm{high}})=(1,2,3)\) on the moving \([-0.98,0.98]\) transition band with a one-ring \(Q_2\) halo, and initial analytical seed pre-refinement achieving exact octree conformity (\(\max\text{jump}=8.88\times 10^{-16}\)). The equivalent uniform \(L_5\,Q_3\) spectral-element mesh requires \(32{,}768\) elements and \(N_{\mathrm{unif}}=912{,}673\) conforming unknowns per field per octant (\(7{,}189{,}057\) unknowns per field on the full cube, \(1.5\times\) finer per coordinate than a \(128^3\) periodic Fourier grid).

\begin{figure}[htbp]
\centering
\includegraphics[width=\linewidth]{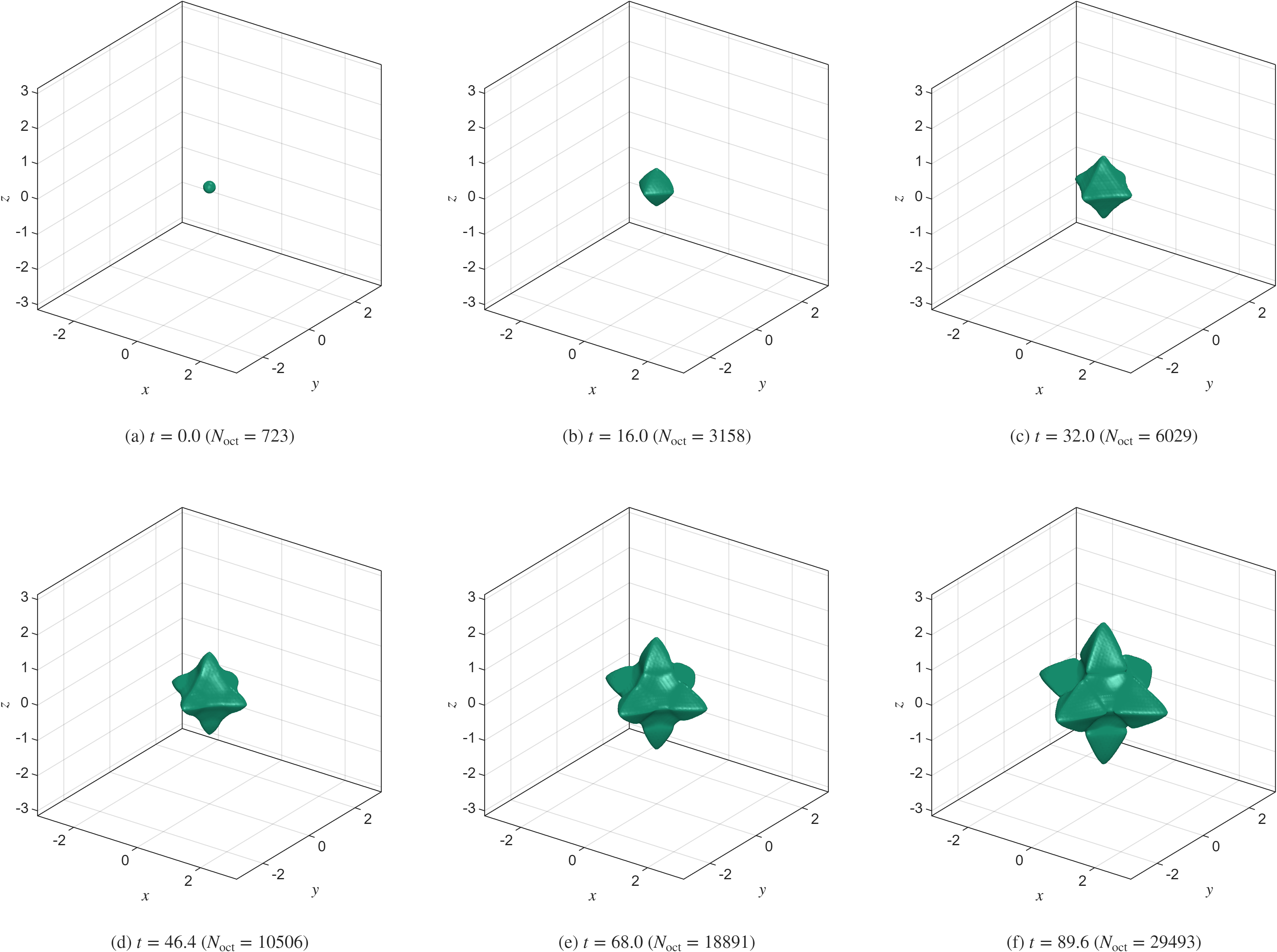}
\caption{Example~3 (3D cubic anisotropy, \(\delta=0.10,\ \mu_A=0.55,\ K=0.90\)): Full-cube dendritic formation process (\(\phi=0\) isosurfaces on \([-\pi,\pi]^3\) under homogeneous Neumann BCs) across six saved field snapshots from \(t=0\) to \(t=89.6\) (diagnostic trajectory continues to \(t=94.4\) in Figure~\ref{fig:exp3-3d-detail} and Table~\ref{tab:exp3-3d}). Octant conforming unknown counts \(N_{\mathrm{oct}}\) are listed below each panel (compared with \(912{,}673\) octant DOFs on the equivalent uniform \(L_5\,Q_3\) mesh).}
\label{fig:exp3-3d-process}
\end{figure}

We advance \(840\) five-stage CSRK steps (\(k=0.10\) for \(320\) steps over \(t\in[0,32.0]\), checking adaptation every \(6\) steps, followed by \(k=0.12\) for \(520\) steps over \(t\in[32.0,94.40]\), checking adaptation every \(5\) steps, excluding segment-final steps) with \(137\) actual dynamic \(2{:}1\)-balanced octree adaptations (octant elements growing from \(120\) to \(2{,}584\) vs.\ \(32{,}768\) uniform). Figure~\ref{fig:exp3-3d-process} presents six full-cube \(\phi=0\) isosurfaces tracing the three-dimensional dendritic formation process: from the small spherical nucleus (\(R_0=\sqrt{0.027}\approx 0.1643\), \(r_{\mathrm{ax}}=0.1635\), \(N_{\mathrm{oct}}=723\), a \(1{,}262\times\) initial DOF reduction), six prominent cubic arms emerge along the \(\pm x, \pm y, \pm z\) crystal axes, sharpen markedly, and develop pronounced axial ribs and concave \([111]\) valleys, reaching \(r_{\mathrm{ax}}=2.1126\) versus \(r_{\mathrm{diag}}=0.7635\) at \(t=94.40\) (3D tip-to-valley ratio \(r_{\mathrm{ax}}/r_{\mathrm{diag}}=2.767\)).

\begin{figure}[htbp]
\centering
\includegraphics[width=\linewidth]{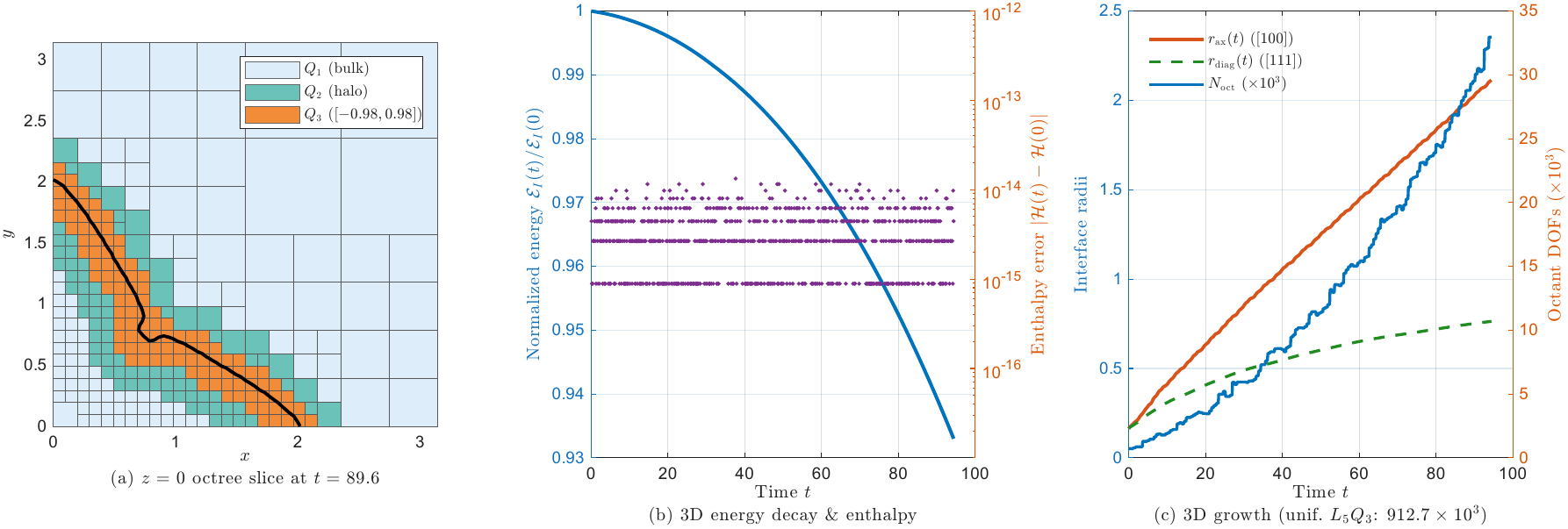}
\caption{Example~3 (3D): \(2{:}1\)-balanced octree adaptivity and thermodynamic diagnostics. (a)~\(z=0\) slice of the octree at \(t=89.6\) colored by polynomial degree (\(Q_1\) bulk, \(Q_2\) halo, \(Q_3\) on \([-0.98,0.98]\)) with \(\phi=0\) contour. (b)~Strictly monotonic decay of the original discrete energy \(\Eorig(t)/\Eorig(0)\) and fixed-mesh time-step integrated enthalpy error \(|\mathcal{H}(t)-\mathcal{H}(0)|\le 1.33\times 10^{-14}\) over the full trajectory \(t\in[0,94.4]\). (c)~Evolution of axial \([100]\) and diagonal \([111]\) radii alongside the adaptive octant DOF count \(N_{\mathrm{oct}}(t)\).}
\label{fig:exp3-3d-detail}
\end{figure}

Figure~\ref{fig:exp3-3d-detail} summarizes the 3D octree grading and structure-preserving diagnostics (octant energies and per-field octant DOFs): panel~(a) displays the \(z=0\) slice of the displayed \(2{:}1\)-balanced octree snapshot (\(t=89.6\)) colored by polynomial degree, confirming crisp localization of \(Q_3\) elements to the curved interfacial band and smooth \(2{:}1\) coarsening into the interior solid and outer melt; panel~(b) verifies monotonic decay of the octant quadrature-defined original discrete energy from \(\Eorig(0)=138{,}944.71\) to \(\Eorig(94.40)=129{,}635.25\) (every time step decreases energy with \(\max\Delta\Eorig^{\mathrm{step}}=-1.427<0\), max Newton residual \(5.17\times 10^{-11}\), \(\min\alpha_T=1-3.08\times 10^{-6}\) and peak transfer floating-point rounding \(\max\Delta\Eorig^{\mathrm{tr}}=+2.91\times 10^{-11}\le 2.3\times 10^{-16}\Eorig\) across all \(137\) octree transfers, remaining nonincreasing within the stated IEEE tolerance) together with fixed-mesh time-step integrated enthalpy conservation within \(1.33\times 10^{-14}\) throughout (an end-to-end \(840\)-step re-execution of the frozen release driver reproduces the initial/final energies, axial and diagonal radii, step/transfer energy bounds, DOF counts, and all \(137\) adaptation events to full \(15\)-digit IEEE precision while recording all \(840\) steps with \(0\) step retries, \(0\) rejected adaptations, peak thermal residual \(2.99\times 10^{-13}\) across the \(>30{,}000\)-DOF iterative tail, peak transfer enthalpy defect \(1.15\times 10^{-14}\) across all \(137\) octree transfers, and an exact Step-\(840\) terminal snapshot); panel~(c) tracks the rapid divergence between axial and body-diagonal radii alongside the compact octree DOF count \(N_{\mathrm{oct}}\), which averages only \(12{,}124\) unknowns per field (\(75.28\times\) average DOF compression relative to the matched theoretical uniform \(L_5\,Q_3\) count) and remains at \(32{,}961\) (\(27.69\times\) compression) even when the mature six-arm dendrite spans over two-thirds of the cube.

The entire \(840\)-step 3D trajectory completes in \(7{,}401.4\)\,s (\(8.81\)\,s/step average, of which all \(137\) octree adaptations and energy-safe transfers account for \(730.7\)\,s, or \(9.9\%\)). As an indicative short-window uniform timing reference on our i9-12900K workstation, even a uniform \(L_4\,Q_3\) mesh with \(8\times\) fewer elements (\(2\times\) coarser mesh size \(h\); \(4{,}096\) octant cells, \(117{,}649\) octant unknowns per field) requires \(345.4\)\,s per CSRK step over a two-step timing window (\(351.6\)\,s and \(339.3\)\,s), compared with \(8.81\)\,s/step averaged over the \(840\)-step finer \(L_5\,Q_3\) adaptive run (Table~\ref{tab:exp3-3d}; such per-step wall-clock ratios depend on hardware, thread load, and the two-step versus full-trajectory averaging window and are reported only as an order-of-magnitude timing reference rather than an equal-accuracy speedup).

\begin{table}[htbp]
\centering
\caption{Example~3: Verified performance and structure-preservation summary of the 3D \(2{:}1\)-balanced octree \(hp\)-SEM + five-stage CSRK solver on the symmetry octant \([0,\pi]^3\) (per-field octant DOFs \(N_{\mathrm{oct}}\), octant energies \(\Eorig\), and peak fixed-mesh time-step integrated enthalpy defect \(\max|\Delta\mathcal{H}|\)) compared with a two-step uniform \(L_4\,Q_3\) timing baseline (\(8\times\) fewer elements, \(2\times\) coarser \(h\)) and the matched theoretical \(L_5\,Q_3\) DOF count.}
\label{tab:exp3-3d}
\footnotesize
\setlength{\tabcolsep}{4pt}
\resizebox{\linewidth}{!}{%
\begin{tabular}{@{}lccccccc@{}}
\toprule
Discretization & Steps (\(t_{\mathrm{end}}\)) & \(N_{\mathrm{oct}}\) (init / avg / end) & DOF compression & \(r_{\mathrm{ax}}/r_{\mathrm{diag}}\) & \(\Eorig(0)\to\Eorig(t_{\mathrm{end}})\) & \(\max|\Delta\mathcal{H}|\) & Avg.\ time/step \\
\midrule
Adaptive \(L_5\,Q_3\) (\(hp\)) & 840 (94.4) & \(723\) / \(12{,}124\) / \(32{,}961\) & \(1{,}262\times\) / \(75.3\times\) / \(27.7\times\) & 2.767 & \(138{,}944.7\to 129{,}635.3\) & \(1.33\times 10^{-14}\) & \(8.81\)\,s \\
Uniform \(L_4\,Q_3\) (\(8\times\) fewer cells) & 2 (bench) & \(117{,}649\) (fixed) & -- & -- & -- & \(4.0\times 10^{-15}\) & \(345.4\)\,s (2-step ref.) \\
Uniform \(L_5\,Q_3\) (matched DOFs) & -- & \(912{,}673\) (fixed) & \(1.0\times\) & -- & -- & -- & intractable \\
\bottomrule
\end{tabular}%
}
\end{table}

\begin{figure}[htbp]
\centering
\includegraphics[width=\linewidth]{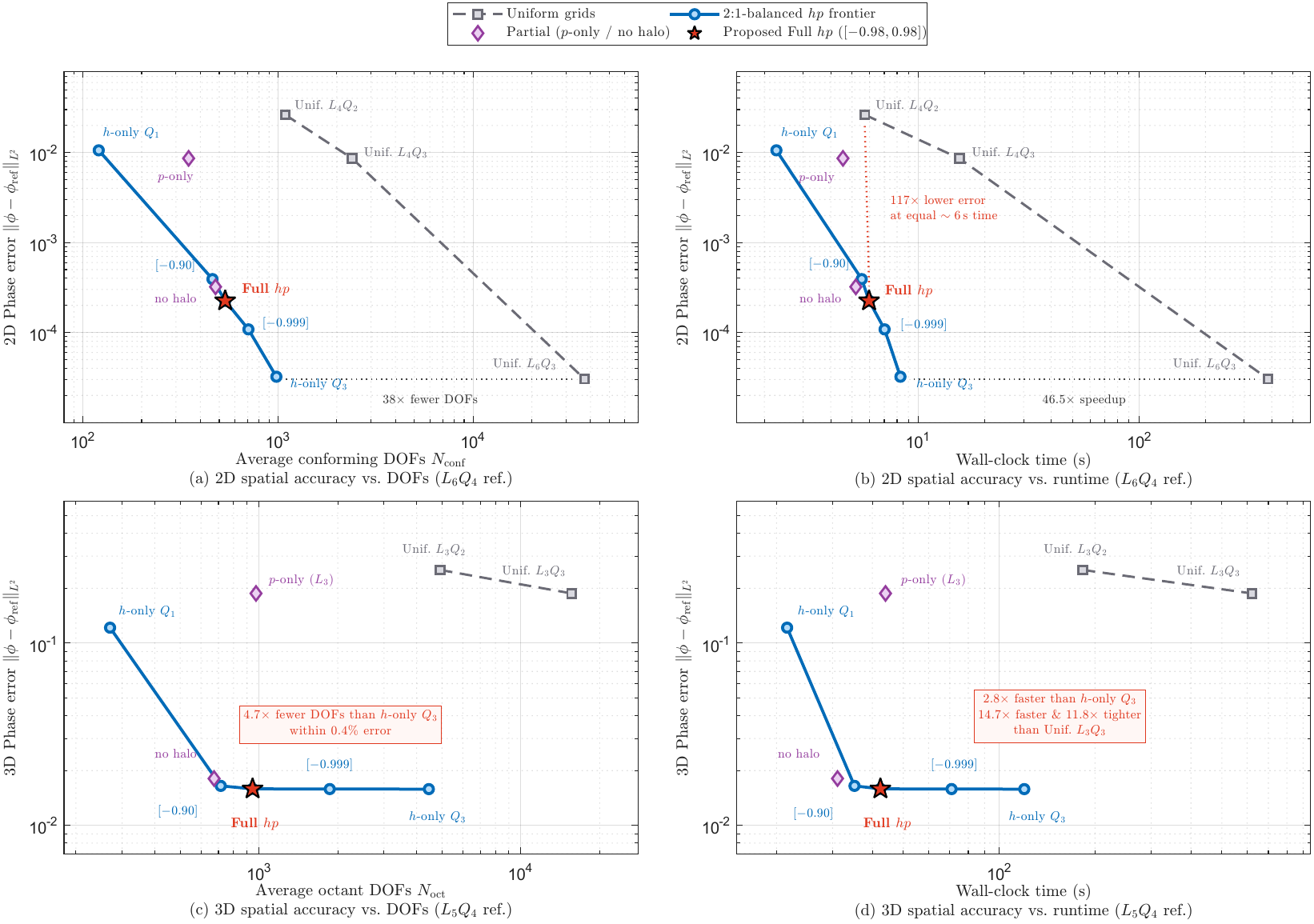}
\caption{Unified 2D and 3D accuracy--cost Pareto comparison from the equal-condition independent higher-order \(Q_4\) spatial ablation studies (Tables~\ref{tab:exp2-ablation} and~\ref{tab:exp3-ablation}). (a)--(b)~Two-dimensional phase \(L^2\) error versus average conforming DOFs per field and wall-clock time (\(L_6\,Q_4\) reference, \(t=1.50\)). (c)--(d)~Three-dimensional phase \(L^2\) error versus average octant DOFs per field and wall-clock time (\(L_5\,Q_4\) reference, \(t=5.00\)), exposing the sharp 3D knee where the proposed zoned Full \(hp\) scheme (\([-0.98,0.98]\), \(Q_3/Q_2/Q_1\), red star) stays within \(0.35\%\) phase \(L^2\) error of \(h\)-only \(Q_3\) at \(4.7\times\) fewer DOFs and \(2.8\times\) lower runtime.}
\label{fig:ablation-pareto}
\end{figure}

\paragraph{3D equal-condition spatial accuracy and component ablation study.}
To verify whether the 2D ablation conclusions carry over---and sharpen---in three dimensions, we integrate the exact Example~3 cubic-anisotropy solidification problem (\(\delta=0.10,\ \mu_A=0.55,\ K=0.90\)) on \([0,\pi]^3\) for \(50\) five-stage CSRK steps (\(k=0.10,\ t=5.00\), adapting every \(5\) steps for adaptive cases, during which the nucleus grows from \(r_{\mathrm{ax}}=0.1635\) to \(r_{\mathrm{ax}}=0.2955,\ r_{\mathrm{diag}}=0.2397\), developing clear cubic anisotropy \(r_{\mathrm{ax}}/r_{\mathrm{diag}}=1.233\)). Because uniform \(L_5\,Q_3\) (\(912{,}673\) octant DOFs per field) is intractable, we compute an independent higher-order 3D spatial reference on a strictly finer-degree \(2{:}1\)-balanced \(L_5\,Q_4\) octree covering the wider band \([-0.999,0.999]\) with \(p\equiv 4\) globally (\(14{,}533.9\) average octant DOFs per field, \(432.10\)\,s), and evaluate exact common-subdivision 3D \(L^2(\Omega)\) errors, relative energy errors, and axial/diagonal radius errors across nine candidate discretizations (Table~\ref{tab:exp3-ablation} and the unified 2D/3D accuracy--cost Pareto frontiers in Figure~\ref{fig:ablation-pareto}).

The 3D ablation reinforces every element of the proposed spatial strategy even more decisively than in 2D:
\begin{enumerate}
\item \textbf{Zoned \(hp\) closely matches \(h\)-only \(Q_3\) accuracy at \(4.7\times\) fewer DOFs and \(2.8\times\) speedup.} Whereas bulk \(p\)-coarsening yielded a \(1.8\times\) DOF saving in the small 2D benchmark, in 3D the full zoned \(hp\) scheme (\([-0.98,0.98]\), \(Q_3/Q_2/Q_1\)) reproduces uniform-degree \(h\)-only \(Q_3\) to within \(0.35\%\) in \(\|\phi-\phi_{\mathrm{ref}}\|_{L^2}\) (\(1.589\times 10^{-2}\) vs.\ \(1.583\times 10^{-2}\), i.e., \(2.85\times 10^{-3}\) relative \(L^2\) error on \([0,\pi]^3\)), \(0.21\%\) in \(T\) (\(9.47\times 10^{-4}\) relative), and \(0.85\%\) in energy (and \(7.1\%\) lower axial-radius error) while reducing average DOFs by \(4.71\times\) (\(944.9\) vs.\ \(4{,}452.3\)) and wall time by \(2.83\times\) (\(42.38\)\,s vs.\ \(119.89\)\,s). Because the local \(L_5\,Q_3\) interface error already dominates the outer bulk thermal \(Q_1\) error in 3D, retaining \(Q_3\) in the bulk buys virtually no accuracy at nearly fivefold DOF cost.
\item \textbf{Combined octree \(h\)-branching and \(Q_3\) enrichment outperforms uniform 3D grids by over an order of magnitude.} Relative to uniform \(L_3\,Q_3\), full \(hp\) uses \(16.5\times\) fewer DOFs and runs \(14.7\times\) faster while improving accuracy by \(11.8\times\) in \(\phi\), \(19.6\times\) in \(T\), \(22.2\times\) in energy, and \(27.6\times\) in \(r_{\mathrm{ax}}\). At virtually identical cost to \(p\)-only (\(L_3\)) (\(944.9\) vs.\ \(974.1\) DOFs, \(42.38\)\,s vs.\ \(43.97\)\,s), adding \(2{:}1\)-balanced octree \(h\)-branching reduces \(\phi\), \(T\), and \(r_{\mathrm{ax}}\) errors by \(11.8\times\), \(19.5\times\), and \(27.5\times\); conversely, stripping spectral enrichment down to \(h\)-only \(Q_1\) degrades errors by \(7.7\times\)--\(16.2\times\).
\item \textbf{Observed impact of the \([-0.98,0.98]\) band and \(Q_2\) thermal halo in 3D.} Removing the one-ring \(Q_2\) halo degrades the 3D thermal \(L^2\) error by \(41\%\), diagonal radius error by \(38\%\), energy error by \(29\%\), and phase error by \(14\%\). Narrowing the band to \([-0.90,0.90]\) worsens \(|r_{\mathrm{ax}}-r_{\mathrm{ax}}^{\mathrm{ref}}|\) by \(28\%\) and triggers a \(4{,}860\times\) larger inter-mesh thermal contraction (\(1-\alpha_T=1.52\times 10^{-3}\) vs.\ \(3.13\times 10^{-7}\)). Conversely, widening to \([-0.999,0.999]\) nearly doubles the 3D DOF count (\(+97\%\)) and adds \(67\%\) runtime for only \(0.3\%\) phase-error reduction, making \([-0.98,0.98]\) with a one-ring \(Q_2\) halo the best accuracy--cost trade-off among the tested 3D configurations.
\end{enumerate}

\begin{table}[htbp]
\centering
\caption{Example~3 (3D cubic anisotropy, \(\delta=0.10,\ K=0.90,\ k=0.10,\ 50\text{ steps},\ t=5.00\)): Equal-condition 3D spatial accuracy and component ablation study measured by exact common-subdivision 3D \(L^2\) integration against the independent higher-order \(2{:}1\)-balanced \(L_5\,Q_4\) octree spatial reference (\([-0.999,0.999]\), \(p\equiv 4\), \(14{,}533.9\) avg.\ octant DOFs per field, \(432.10\)\,s).}
\label{tab:exp3-ablation}
\footnotesize
\setlength{\tabcolsep}{4pt}
\resizebox{\linewidth}{!}{%
\begin{tabular}{@{}lccccccc@{}}
\toprule
3D spatial discretization & Avg.\ DOFs & Wall (s) & \(\|\phi-\phi_{\mathrm{ref}}\|_{L^2}\) & \(\|T-T_{\mathrm{ref}}\|_{L^2}\) & \(|\Delta\Eorig|/\Eorig^{\mathrm{ref}}\) & \(|r_{\mathrm{ax}}-r_{\mathrm{ax}}^{\mathrm{ref}}|\) & \(|r_{\mathrm{diag}}-r_{\mathrm{diag}}^{\mathrm{ref}}|\) \\
\midrule
Uniform \(L_3\,Q_2\) (coarse baseline) & \(4{,}913\) & \(183.00\) & \(2.520\times 10^{-1}\) & \(8.954\times 10^{-2}\) & \(1.154\times 10^{-3}\) & \(2.255\times 10^{-1}\) & \(1.032\times 10^{-1}\) \\
Uniform \(L_3\,Q_3\) (coarse baseline) & \(15{,}625\) & \(621.64\) & \(1.874\times 10^{-1}\) & \(5.672\times 10^{-2}\) & \(7.454\times 10^{-4}\) & \(1.665\times 10^{-1}\) & \(5.015\times 10^{-2}\) \\
\midrule
\(p\)-only (\(L_3\), zoned \(Q_3/Q_2/Q_1\)) & \(974.1\) & \(43.97\) & \(1.873\times 10^{-1}\) & \(5.667\times 10^{-2}\) & \(7.450\times 10^{-4}\) & \(1.658\times 10^{-1}\) & \(5.020\times 10^{-2}\) \\
\(h\)-only \(Q_1\) (\(2{:}1\) \(L_5\) tree, \(p\equiv 1\)) & \(269.9\) & \(21.59\) & \(1.216\times 10^{-1}\) & \(3.204\times 10^{-2}\) & \(3.882\times 10^{-4}\) & \(9.766\times 10^{-2}\) & \(3.583\times 10^{-2}\) \\
\(hp\), narrower band \([-0.90,0.90]\) & \(714.3\) & \(35.14\) & \(1.649\times 10^{-2}\) & \(3.023\times 10^{-3}\) & \(3.531\times 10^{-5}\) & \(7.701\times 10^{-3}\) & \(3.483\times 10^{-3}\) \\
\(hp\), no \(Q_2\) halo (\(\text{halo}=0\)) & \(673.6\) & \(31.08\) & \(1.809\times 10^{-2}\) & \(4.082\times 10^{-3}\) & \(4.326\times 10^{-5}\) & \(6.000\times 10^{-3}\) & \(4.594\times 10^{-3}\) \\
\textbf{Full \(hp\) (\([-0.98,0.98]\), \(Q_3/Q_2/Q_1\))} & \(\mathbf{944.9}\) & \(\mathbf{42.38}\) & \(\mathbf{1.589\times 10^{-2}}\) & \(\mathbf{2.899\times 10^{-3}}\) & \(\mathbf{3.353\times 10^{-5}}\) & \(\mathbf{6.038\times 10^{-3}}\) & \(\mathbf{3.337\times 10^{-3}}\) \\
\(hp\), wider band \([-0.999,0.999]\) & \(1{,}862.6\) & \(70.89\) & \(1.584\times 10^{-2}\) & \(2.894\times 10^{-3}\) & \(3.328\times 10^{-5}\) & \(6.576\times 10^{-3}\) & \(3.305\times 10^{-3}\) \\
\(h\)-only \(Q_3\) (\(2{:}1\) \(L_5\) tree, \(p\equiv 3\)) & \(4{,}452.3\) & \(119.89\) & \(1.583\times 10^{-2}\) & \(2.893\times 10^{-3}\) & \(3.325\times 10^{-5}\) & \(6.499\times 10^{-3}\) & \(3.306\times 10^{-3}\) \\
\midrule
\(2{:}1\) \(L_5\,Q_4\) (\([-0.999,0.999]\), ref.) & \(14{,}533.9\) & \(432.10\) & \(0\) (ref) & \(0\) (ref) & \(0\) (ref) & \(0\) (ref) & \(0\) (ref) \\
\bottomrule
\end{tabular}%
}
\end{table}

\FloatBarrier
% ============================================================================
\section{Conclusions}\label{sec:conclusions}
% ============================================================================
We have presented a unified, thermodynamically consistent numerical framework for two- and three-dimensional anisotropic dendritic solidification that combines enthalpy-conserving gradient-flow coordinates, second-order five-stage convex splitting, and \(2{:}1\)-balanced nonmatching \(hp\) spectral-element adaptivity. By bounding the curvature of the quintic latent-heat enthalpy lift with a sextic majorizer and incorporating an explicit anisotropy stabilization constant \(\mu_A\), we established a global convex splitting of the \emph{original} physical free energy and proved local second-order temporal consistency on smooth fixed-space trajectories, unique solvability, exact physical enthalpy conservation, and unconditional nonincrease of the original quadrature-defined discrete energy for every positive step size at both the time-semidiscrete and fully discrete levels, reinforced by a computable two-component a posteriori energy certificate covering both inexact phase Newton steps and iterative 3D thermal solves.

On the spatial side, full \(2{:}1\) balance on nonmatching quadtrees and octrees enables an explicit \(O(N)\)-nonzero conforming trace injection with continuous inter-element traces (verified numerically to \(\sim 10^{-16}\)) and smooth geometric grading across the thermal boundary layer. Focusing finest elements and \(Q_3\) enrichment on sampled \(\phi\in[-0.98,0.98]\) with a one-ring \(Q_2\) thermal halo, bulk \(Q_1\) coarsening, \(O(L_{\max}(N_{\mathrm{old}}+N_{\mathrm{new}}))\) dyadic common-subdivision transfer, and an enthalpy-preserving, energy-safe radial thermal projection with an explicit fallback acceptance guard guarantees global quadrature-energy monotonicity across all accepted updates. Exact Schur-complement elimination and matrix-free 3D sum-factorized phase Hessians resolve full 2D (\(m=4,6\)), full-domain asymmetric rotated and multi-nuclei, and 3D cubic dendritic trajectories with \(7.4\times\)--\(75.3\times\) average DOF compression, \(0.003\%\) tip-radius agreement between rotated and axis-aligned configurations at \(t=15.0\), and substantial equal-condition runtime savings (\(64.5\times\) full-run speedup over uniform \(L_6\,Q_3\) in the 2D ablation, and \(2.8\times\) speedup over \(h\)-adaptive \(Q_3\) in the 3D ablation). Direct same-mesh temporal benchmarking against a linearized SAV--BDF2 implementation clarifies the complementary trade-off between linear modified-energy integrators and our stabilization-free original-energy-dissipative CSRK scheme, while equal-condition 2D and 3D ablations against independent higher-order \(Q_4\) spatial references quantify the benefits of balanced tree \(h\)-branching, \(Q_3\) enrichment, the \([-0.98,0.98]\) band, and the \(Q_2\) halo, showing that bulk \(p\)-coarsening in 3D cuts DOFs by \(4.7\times\) and runtime by \(2.8\times\) relative to \(h\)-adaptive \(Q_3\) octrees at \(0.35\%\) phase \(L^2\) error difference (\(t=5.00\)).

Several directions merit future investigation. First, combining our \(2{:}1\)-balanced octree \(hp\)-SEM discretization with distributed-memory forest-of-octrees partitioning \cite{burstedde2011} and algebraic multigrid preconditioning will enable full-domain 3D simulations with deep secondary and tertiary sidebranching and arbitrary crystal orientations without octant symmetry. Second, leveraging the intermediate CSRK stages for a posteriori adaptive time-step control can further accelerate long-time coarsening and late-stage solidification while retaining unconditional original-energy dissipation. Finally, extending the enthalpy-conserving gradient-flow splitting to thermo-solutal alloy solidification and Navier--Stokes melt convection represents a natural next step.

\paragraph{Code and data availability.}
The complete MATLAB implementation of the 2D/3D five-stage CSRK solver, \(2{:}1\)-balanced quadtree/octree \(hp\)-SEM discretization, experiment-named drivers reproducing Examples~1--3, figure scripts, and archived benchmark summary data supporting all tables and figures in this manuscript are openly available at \url{https://github.com/minghuili123456789/anisotropic-dendrite-csrk-hp}.


\begin{thebibliography}{41}
\setlength{\itemsep}{1.5pt plus 0.5pt}
\setlength{\parskip}{0pt}
\bibitem{babuska1994}
I.~Babu\v{s}ka and M.~Suri, The \(p\) and \(h\)-\(p\) versions of the finite element method, basic principles and properties, \emph{SIAM Rev.} 36 (1994), 578--632.

\bibitem{boettinger2002}
W.~J.~Boettinger, J.~A.~Warren, C.~Beckermann, and A.~Karma, Phase-field simulation of solidification, \emph{Annu.\ Rev.\ Mater.\ Res.} 32 (2002), 163--194.

\bibitem{burstedde2011}
C.~Burstedde, L.~C.~Wilcox, and O.~Ghattas, \texttt{p4est}: Scalable algorithms for parallel adaptive mesh refinement on forests of octrees, \emph{SIAM J.\ Sci.\ Comput.} 33 (2011), 1103--1133.

\bibitem{caginalp1986}
G.~Caginalp, An analysis of a phase field model of a free boundary, \emph{Arch.\ Rational Mech.\ Anal.} 92 (1986), 205--245.

\bibitem{demkowicz2006}
L.~Demkowicz, \emph{Computing with \(hp\)-Adaptive Finite Elements, Vol.~1: One and Two Dimensional Elliptic and Maxwell Problems}, Chapman \& Hall/CRC, Boca Raton, 2006.

\bibitem{deville2002}
M.~O.~Deville, P.~F.~Fischer, and E.~H.~Mund, \emph{High-Order Methods for Incompressible Fluid Flow}, Cambridge University Press, Cambridge, 2002.

\bibitem{eggleston2001}
J.~J.~Eggleston, G.~B.~McFadden, and P.~W.~Voorhees, A phase-field model for highly anisotropic interfacial energy, \emph{Physica D} 150 (2001), 91--103.

\bibitem{eyre1998}
D.~J.~Eyre, Unconditionally gradient stable time marching the Cahn--Hilliard equation, \emph{MRS Online Proc.\ Libr.} 529 (1998), 39--46.

\bibitem{greenwood2018}
M.~Greenwood, K.~N.~Shampur, N.~Ofori-Opoku, T.~Pinomaa, L.~Wang, S.~Gurevich, and N.~Provatas, Quantitative 3D phase field modelling of solidification using next-generation adaptive mesh refinement, \emph{Comput.\ Mater.\ Sci.} 142 (2018), 153--171.

\bibitem{guo2024}
Y.~Guo, M.~Azaiez, and C.~Xu, An efficient numerical method for the anisotropic phase field dendritic crystal growth model, \emph{Commun.\ Nonlinear Sci.\ Numer.\ Simul.} 131 (2024), 107858.

\bibitem{jeong2001}
J.-H.~Jeong, N.~Goldenfeld, and J.~A.~Dantzig, Phase field model for three-dimensional dendritic growth with fluid flow, \emph{Phys.\ Rev.\ E} 64 (2001), 041602.

\bibitem{karma1996}
A.~Karma and W.-J.~Rappel, Quantitative phase-field modeling of dendritic growth in two and three dimensions, \emph{Phys.\ Rev.\ E} 57 (1998), 4323--4349.

\bibitem{karniadakis2005}
G.~E.~Karniadakis and S.~J.~Sherwin, \emph{Spectral/\(hp\) Element Methods for Computational Fluid Dynamics}, 2nd ed., Oxford University Press, Oxford, 2005.

\bibitem{kobayashi}
R.~Kobayashi, Modeling and numerical simulations of dendritic crystal growth, \emph{Physica D} 63 (1993), 410--423.

\bibitem{lan2003}
C.~W.~Lan, Y.~C.~Chang, and C.~J.~Shih, Adaptive phase field simulation of non-isothermal free dendritic growth of a binary alloy, \emph{Acta Mater.} 51 (2003), 1857--1869.

\bibitem{lithesis}
M.~Li, \emph{Highly Efficient Numerical Methods for Phase-Field Models of Multiphase Systems}, Doctoral Dissertation, Xiamen University, 2022.

\bibitem{li2022}
M.~Li, M.~Azaiez, and C.~Xu, New efficient time-stepping schemes for the anisotropic phase-field dendritic crystal growth model, \emph{Comput.\ Math.\ Appl.} 109 (2022), 204--215.

\bibitem{likim2012}
Y.~Li and J.~Kim, Phase-field simulations of crystal growth with adaptive mesh refinement, \emph{Int.\ J.\ Heat Mass Transfer} 55 (2012), 7926--7932.

\bibitem{likim2011}
Y.~Li, H.-G.~Lee, and J.~Kim, A fast, robust, and accurate operator splitting method for phase-field simulations of crystal growth, \emph{J.\ Cryst.\ Growth} 321 (2011), 176--182.

\bibitem{likim2023}
Y.~Li, K.~Qin, Q.~Xia, and J.~Kim, A second-order unconditionally stable method for the anisotropic dendritic crystal growth model with an orientation-field, \emph{Appl.\ Numer.\ Math.} 184 (2023), 512--526.

\bibitem{mullins1963}
W.~W.~Mullins and R.~F.~Sekerka, Morphological stability of a particle growing by diffusion or heat flow, \emph{J.\ Appl.\ Phys.} 34 (1963), 323--329.

\bibitem{penrose1990}
O.~Penrose and P.~C.~Fife, Thermodynamically consistent models of phase-field type for the kinetics of phase transitions, \emph{Physica D} 43 (1990), 44--62.

\bibitem{provatas1999}
N.~Provatas, N.~Goldenfeld, and J.~Dantzig, Adaptive mesh refinement computation of solidification microstructures using dynamic data structures, \emph{J.\ Comput.\ Phys.} 148 (1999), 265--290.

\bibitem{rosam2007}
J.~Rosam, P.~K.~Jimack, and A.~Mullis, A fully implicit, fully adaptive time and space discretisation method for phase-field simulation of binary alloy solidification, \emph{J.\ Comput.\ Phys.} 225 (2007), 1271--1287.

\bibitem{schmidt1996}
A.~Schmidt, Computation of three dimensional dendrites with finite elements, \emph{J.\ Comput.\ Phys.} 125 (1996), 293--312.

\bibitem{shentangwang}
J.~Shen, T.~Tang, and L.-L.~Wang, \emph{Spectral Methods: Algorithms, Analysis and Applications}, Springer, Berlin, 2011.

\bibitem{shen2012}
J.~Shen, C.~Wang, X.~Wang, and S.~M.~Wise, Second-order convex splitting schemes for gradient flows with Ehrlich--Schwoebel type energy: Application to thin film epitaxy, \emph{SIAM J.\ Numer.\ Anal.} 50 (2012), 105--125.

\bibitem{shen2018}
J.~Shen, J.~Xu, and J.~Yang, The scalar auxiliary variable (SAV) approach for gradient flows, \emph{J.\ Comput.\ Phys.} 353 (2018), 407--416.

\bibitem{wang1993}
S.-L.~Wang, R.~F.~Sekerka, A.~A.~Wheeler, B.~T.~Murray, S.~R.~Coriell, R.~J.~Braun, and G.~B.~McFadden, Thermodynamically-consistent phase-field models for solidification, \emph{Physica D} 69 (1993), 189--200.

\bibitem{wang2022}
W.~Wang, M.~Azaiez, and C.~Xu, An unconditionally stable fast high order method for thermal phase change models, \emph{Comput.\ Fluids} 237 (2022), 105306.

\bibitem{wangxiao2022}
Y.~Wang, X.~Xiao, and X.~Feng, An accurate and parallel method with post-processing boundedness control for solving the anisotropic phase-field dendritic crystal growth model, \emph{Commun.\ Nonlinear Sci.\ Numer.\ Simul.} 115 (2022), 106717.

\bibitem{yamanaka2011}
A.~Yamanaka, T.~Aoki, S.~Ogawa, and T.~Takaki, GPU-accelerated phase-field simulation of dendritic solidification in a binary alloy, \emph{J.\ Cryst.\ Growth} 318 (2011), 40--45.

\bibitem{yang2016}
X.~Yang, Linear, first and second-order, unconditionally energy stable numerical schemes for the phase field model of homopolymer blends, \emph{J.\ Comput.\ Phys.} 327 (2016), 294--316.

\bibitem{yang2019}
X.~Yang, Efficient linear, stabilized, second-order time marching schemes for an anisotropic phase field dendritic crystal growth model, \emph{Comput.\ Methods Appl.\ Mech.\ Engrg.} 347 (2019), 316--339.

\bibitem{yang2020}
X.~Yang, Efficient and energy stable scheme for an anisotropic phase-field dendritic crystal growth model using the scalar auxiliary variable (SAV) approach, \emph{J.\ Math.\ Study} 53 (2020), 212--236.

\bibitem{yang2021a}
X.~Yang, Fully-discrete spectral-Galerkin scheme with decoupled structure and second-order time accuracy for the anisotropic phase-field dendritic crystal growth model, \emph{Int.\ J.\ Heat Mass Transfer} 180 (2021), 121750.

\bibitem{yang2021b}
X.~Yang, On a novel full decoupling, linear, second-order accurate, and unconditionally energy stable numerical scheme for the anisotropic phase-field dendritic crystal growth model, \emph{Int.\ J.\ Numer.\ Methods Eng.} 122 (2021), 4129--4153.

\bibitem{zeg}
A.~Zaitzeff, S.~Esedo\u{g}lu, and K.~Garikipati, High order, semi-implicit, energy stable schemes for gradient flows, \emph{J.\ Comput.\ Phys.} 447 (2021), 110688.

\bibitem{zhang2019}
J.~Zhang, C.~Chen, and X.~Yang, A novel decoupled and stable scheme for an anisotropic phase-field dendritic crystal growth model, \emph{Appl.\ Math.\ Lett.} 95 (2019), 122--129.

\bibitem{zhangshen2022}
Y.~Zhang and J.~Shen, A generalized SAV approach with relaxation for dissipative systems, \emph{J.\ Comput.\ Phys.} 464 (2022), 111311.

\bibitem{zhao2017}
J.~Zhao, Q.~Wang, and X.~Yang, Numerical approximations for a phase field dendritic crystal growth model based on the invariant energy quadratization approach, \emph{Int.\ J.\ Numer.\ Methods Eng.} 110 (2017), 279--300.
\end{thebibliography}
\end{document}